\documentclass[3p,times]{elsarticle}

\journal{arXiv}

\usepackage{amssymb,amsmath,mathtools}
\usepackage{amsthm}
        \theoremstyle{plain}

        \newdefinition{rem}{Remark}
\usepackage[section]{algorithm}
\usepackage{algpseudocode}
\usepackage{etoolbox}
\renewcommand{\vec}[1]{\boldsymbol{#1}}
\renewcommand{\Vec}[1]{\mathbf{#1}}

\newcommand{\tensorOne}[1]{\boldsymbol#1}
\newcommand{\tensorTwo}[1]{\boldsymbol#1}
\newcommand{\funSpace}[1]{\mathcal{#1}}
\newcommand{\vecFunSpace}[1]{\boldsymbol{\mathcal{#1}}}
\newcommand{\tensorFour}[1]{\textbf{\sffamily #1}}
\newcommand{\tC}{\widetilde{C}}
\newcommand{\tS}{\widetilde{S}}
\newcommand{\dt}{\Delta t}
\usepackage{subcaption}
\usepackage{arydshln}
\usepackage{booktabs}
\usepackage{multirow}

\usepackage[nodots]{numcompress}

\usepackage{xcolor}

\usepackage[colorlinks=true, allcolors=blue]{hyperref}
\biboptions{sort&compress}
\usepackage{todonotes}

\begin{document}

\begin{frontmatter}

\title{Algebraically stabilized Lagrange multiplier method for frictional contact mechanics with hydraulically active fractures}


\author[Stanford]{Andrea Franceschini\corref{mycorrespondingauthor}}
\ead{af1990@stanford.edu}

\author[LLNL]{Nicola Castelletto}
\ead{castelletto1@llnl.gov}

\author[LLNL]{Joshua A. White}
\ead{jawhite@llnl.gov}

\author[Stanford]{Hamdi A. Tchelepi}
\ead{tchelepi@stanford.edu}

\address[Stanford]{Energy Resources Engineering, Stanford University,
  Stanford, United States}
\address[LLNL]{Atmospheric, Earth, and Energy Division, Lawrence Livermore National
  Laboratory, United States}
\cortext[mycorrespondingauthor]{Corresponding author}

\begin{abstract}
Accurate numerical simulation of coupled fracture/fault deformation and fluid flow is crucial to the performance and safety assessment of many subsurface systems.
In this work, we consider the discretization and enforcement of contact conditions at such surfaces.
The bulk rock deformation is simulated using low-order continuous finite elements, while frictional contact conditions are imposed by means of a Lagrange multiplier method.
We employ a cell-centered finite-volume scheme to solve the fracture fluid mass balance equation.
From a modeling perspective, a convenient choice is to use a single grid for both mechanical and flow processes, with piecewise-constant interpolation of Lagrange multipliers, i.e., contact tractions and fluid pressure.
Unfortunately, this combination of displacement and multiplier variables is not uniformly inf-sup stable, and therefore requires a stabilization technique.
Starting from a macroelement analysis, we develop two algebraic stabilization approaches and compare them in terms of robustness and convergence rate.
The proposed approaches are validated against challenging analytical two- and three-dimensional benchmarks to demonstrate accuracy and robustness.
These benchmarks include both pure contact mechanics problems and well as problems with tightly-coupled fracture flow.
\end{abstract}

\begin{keyword}
Contact mechanics \sep
Lagrange multipliers \sep
Darcy fracture flow \sep
stabilization
\MSC[2010] 65N08 \sep 65N12 \sep 65N30
\end{keyword}

\end{frontmatter}

\allowdisplaybreaks

\section{Introduction}
\label{sec:intro}

To accurately simulate the geomechanical response of a subsurface system,
such as an aquifer or reservoir, it is often important to model faults and fractures \cite{fossen2016structural}. Phenomena such as
micro-seismicity \cite{urbancic1993microseismicity}, fluid leakage
\cite{dockrill2010structural}, fault reactivation \cite{morris1996slip}, and fracture
propagation \cite{zhang2007effects} depend strongly on coupled fluid-structure interaction. As a result, it
often necessary to explicitly model complex hydromechanical behavior at these surfaces
\cite{ferronato2008numerical}. Experimental data
confirms a strong dependence of fracture properties, like conductivity, on 
contact conditions \cite{witherspoon1980validity,zimmerman1996hydraulic,zhang2019experimental}.  The core of the modeling challenge is dealing with a lubricated frictional contact problem \cite{Wri06}. Specifically, we have fluid pressure acting as a normal
forcing term on the surfaces of the fracture, while the conductivity of the fracture is a
strong function of the effective aperture.  This establishes a two-way and highly nonlinear coupling.  A different but related challenge is to accurately
model fracture propagation based on the resulting stress and strain fields near the fracture tip \cite{dahi2011numerical}.

Many approaches have been developed over the years to solve the contact mechanics problem, and an extensive review is beyond the scope of this work. The largest
difference between methods lies in how the discontinuity is fundamentally represented. A first class of models uses a Discrete Fracture Model (DFM), whereby a
conformal computational grid is introduced, and the conservation equations are discretized
using this grid. These methods have the advantage that standard discretization techniques
may be directly applied. A potential drawback is that a conformal mesh for complex fracture
networks may require distorted or excessively-refined meshes.  With DFM-based models, the
discontinuity can be explicitly
discretized using zero-thickness interface elements
\cite{goodman1968model,ferronato2008numerical, GarKarTch16,Set_etal17}. In some models, thin layers
of continuum finite elements with plastic behavior are used instead to mimic the fracture rheology
\cite{rutqvist2002modeling,rutqvist2008coupled,pan2014approach,lee2019tough}.
A second class employs an embedding strategy, in which a continuum discretization
scheme is enriched to capture discontinuities cutting through continuum elements.
For example, Embedded Discrete Fracture Models (EDFM) introduce additional, local
degrees of freedom to capture opening and sliding modes
\cite{shakiba2015using,ren2016fully,wong2019investigation,wu2019integrating} while Extended Finite Elements Methods (XFEM) introduce global
degrees of freedom for this purpose
\cite{deb2009extended,zhang2011extended,mohammadi2012xfem,ren2016fully,flemisch2016review,
vahab2017numerical,khoei2018enriched}.
A third class discards explicit discontinuity surfaces altogether, and instead represents
a fracture via a regularized (smooth) field using a continuum discretization.  This latter
class includes Phase Field and damage-mechanics-based approaches \cite{verhoosel2013phase,
amiri2014phase,wheeler2014augmented,geelen2019phase}.

In this work, we represent fractures explicitly using a conforming discretization. For these
elements, the Karush-Kuhn-Tucker (KKT) conditions \cite{SimHug98}
for normal impenetrability and frictional compatibility must be enforced.
The two most common classes of methods used to
fulfill these conditions are penalty \cite{peric1992computational,zang2011contact,GarKarTch16,
Set_etal17} and Lagrange multiplier methods \cite{hild2010stabilized,JhaJua14,FraFerJanTea16,
Ber_etal19,koppel2019stabilized}.
In penalty methods, constraints are satisfied in an approximate way using stiff
``springs`` connecting the two surfaces of a fracture, and no additional degree of freedom
are introduced. However, penalty methods suffer from several
drawbacks, including  ill-conditioning of the resulting linear systems
and a strong dependence of solution quality on the penalty parameters
\cite{zavarise1998method,FerJanPin12}.
In Lagrange multiplier-based methods,
the KKT conditions are instead imposed directly. This advantage, however,
comes at the cost of adding new global primary variables.
Moreover, the resulting matrices have a
generalized saddle-point structure that requires specialized solvers for efficiency and
scalability \cite{benzi2004preconditioner,BenGolLie05,franceschini2019block}.
Other widely used techniques include Nitsche's method \cite{wriggers2008formulation,
chouly2013nitsche,annavarapu2014nitsche,capatina2016nitsche}, regularized and augmented
Lagrange multipliers \cite{lin1994using,cavalieri2013augmented,pietrzak1999large,
hirmand2015augmented}, and mortar methods \cite{belgacem1997mortar,flemisch2005new,
seitz2016isogeometric,popp2014dual}.
In particular, mortar methods were originally developed to allow matching of
different regions in a domain decomposition framework.  They provide useful flexibility in dealing with large deformations and non-conforming discretization of the
contact area \cite{wohlmuth2000mortar,laursen2006mortar}.

In this paper, we focus entirely on a Lagrange multiplier approach to impose the normal and frictional compatibility conditions on the fracture surfaces.
The displacement field in the bulk is approximated by lowest-order continuous finite elements.
Using lubrication theory \cite{witherspoon1980validity}, the fracture flow discretization
relies on a cell-centered finite volume approach using a two-point flux approximation
(TPFA) scheme \cite{EymGalHer00} for the numerical flux.
We assume that the rock matrix is impermeable, though an extension of the proposed strategy to handle a poroelastic bulk is clear.
Motivated by the coupling between fracture contact and flow processes, we use piecewise-constant equal-order interpolation of Lagrange multipliers, i.e. contact tractions, and fluid pressure.
This is the natural and convenient choice given the mixed finite-element / finite-volume scheme adopted.  Unfortunately, this combination of displacement and contact traction/pressure variables is not uniformly inf-sup stable \cite[Section 3.1]{wohlmuth2011variationally} and requires the implementation of a stabilization methodology.

To rectify this deficiency, possible remedies include enriching the displacement space
with locally supported bubble functions \cite{brezzi1984stabilization,BreMar01,HauLeT07,
caylak2012stabilization} or using a coarser mesh
for the traction/pressure space \cite{Woh01}.  Here, we instead introduce a stabilizing
modification to the mass balance and constraint equations.
We start by analyzing the pure contact mechanics problem, and explore several stabilization alternatives.
First, we derive a local traction jump stabilization that depends on a stabilization
parameter $\alpha$ defined locally at the macroelement level, following the macroelement
analysis \cite{SilKec90,bressan2013isogeometric,boffi2013mixed,he2005stabilized}
originally proposed for the Stokes equation.
Second, to circumvent difficulties associated with the definition of $\alpha$ in the presence of distorted grids and severe material heterogeneity, we develop an algebraic variant of the local traction jump stabilization technique and then generalize it to an algebraic global traction jump stabilization.
Finally, we extend the stabilization procedure to problems that include fluid flow.

The paper is organized as follows. In Section \ref{sec:gen_fram}, we present the general
problem statement, providing the strong form of the governing equations.
In Section \ref{sec:numerical_model}, the discretized form and
its linearization are provided. We also address the nonlinear optimization problem associated
with the imposition of the KKT conditions.
In Section \ref{sec:stabilization}, we describe three related stabilization strategies.  We discuss the relative advantages of each, before settling on a global algebraic stabilization approach as our recommended strategy.  We test the performance of this method on mechanical contact problems in Section \ref{sec:anal_cont_mech} and hydromechanical contact problems in Section \ref{sec:anal_flow}.  We then conclude the paper with a few remarks regarding future work.

\section{Problem statement}
\label{sec:gen_fram}

A fracture in a continuous medium is modeled as an internal discontinuity where additional
constraints must be enforced. These constraints depend on the local state of the discontinuity---that is, whether the fracture is
opening, sliding, or sticking. Inside the fracture,
a single-phase fluid is present and a mass balance equation---expressed at the Darcy scale---is
solved for pressure.  The momentum and mass balance equations are then coupled through the fluid pressure and fracture aperture. For simplicity of exposition, we assume that the matrix is impermeable.  The extension to a poroelastic medium, however, is clearly of interest in practical applications and could be included in a more general formulation \cite{coussy2004poromechanics}.  Similarly, extensions to multiphase flow and non-isothermal conditions may be required for some reservoir applications.

\begin{figure}
\hfill
\begin{subfigure}{0.29\linewidth}
  \centering
  \includegraphics[width=\linewidth]{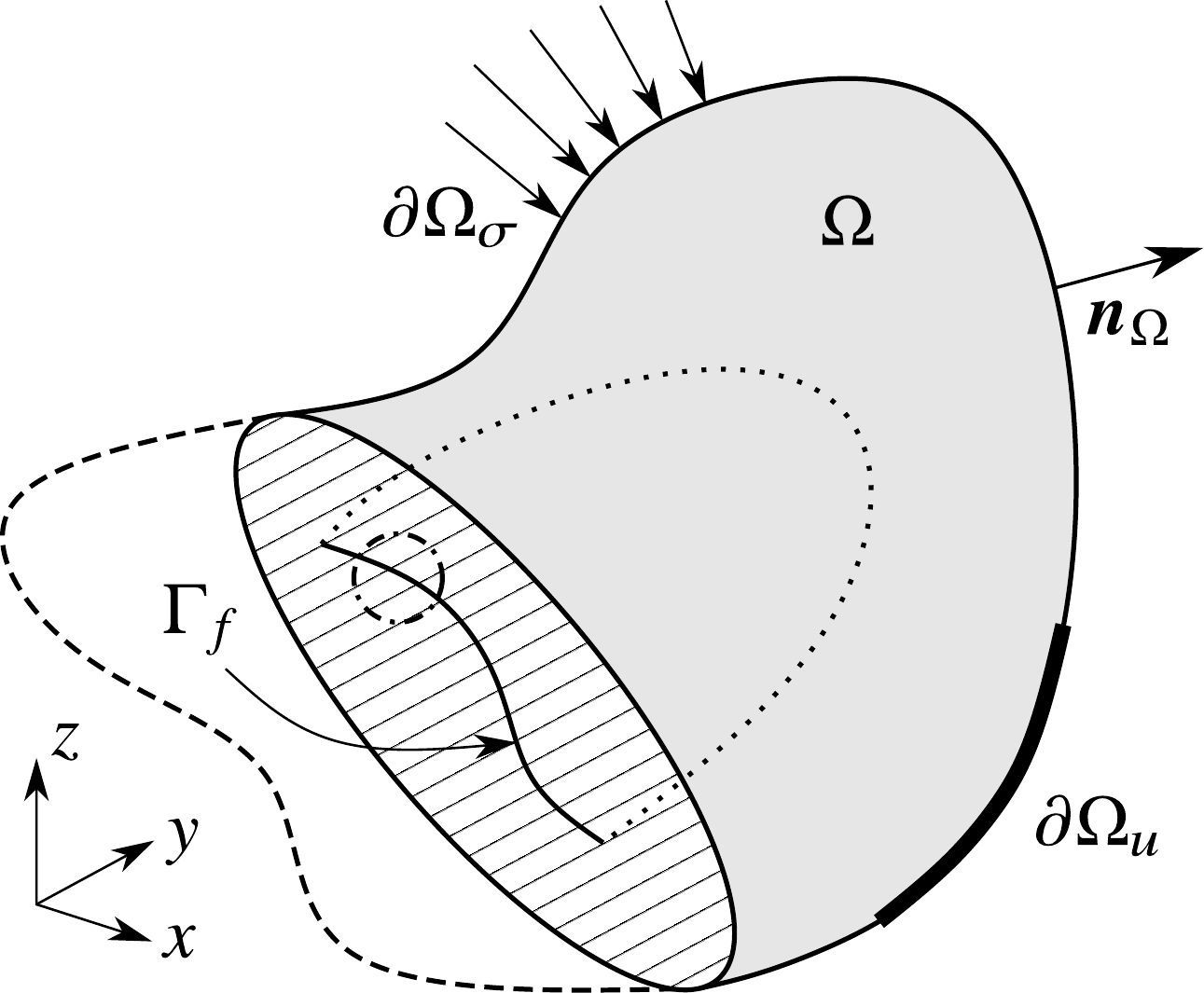}
  \caption{}
  \label{fig:conc_scheme1}
\end{subfigure}
\hfill
\begin{subfigure}{0.29\linewidth}
  \centering
  \includegraphics[width=\linewidth]{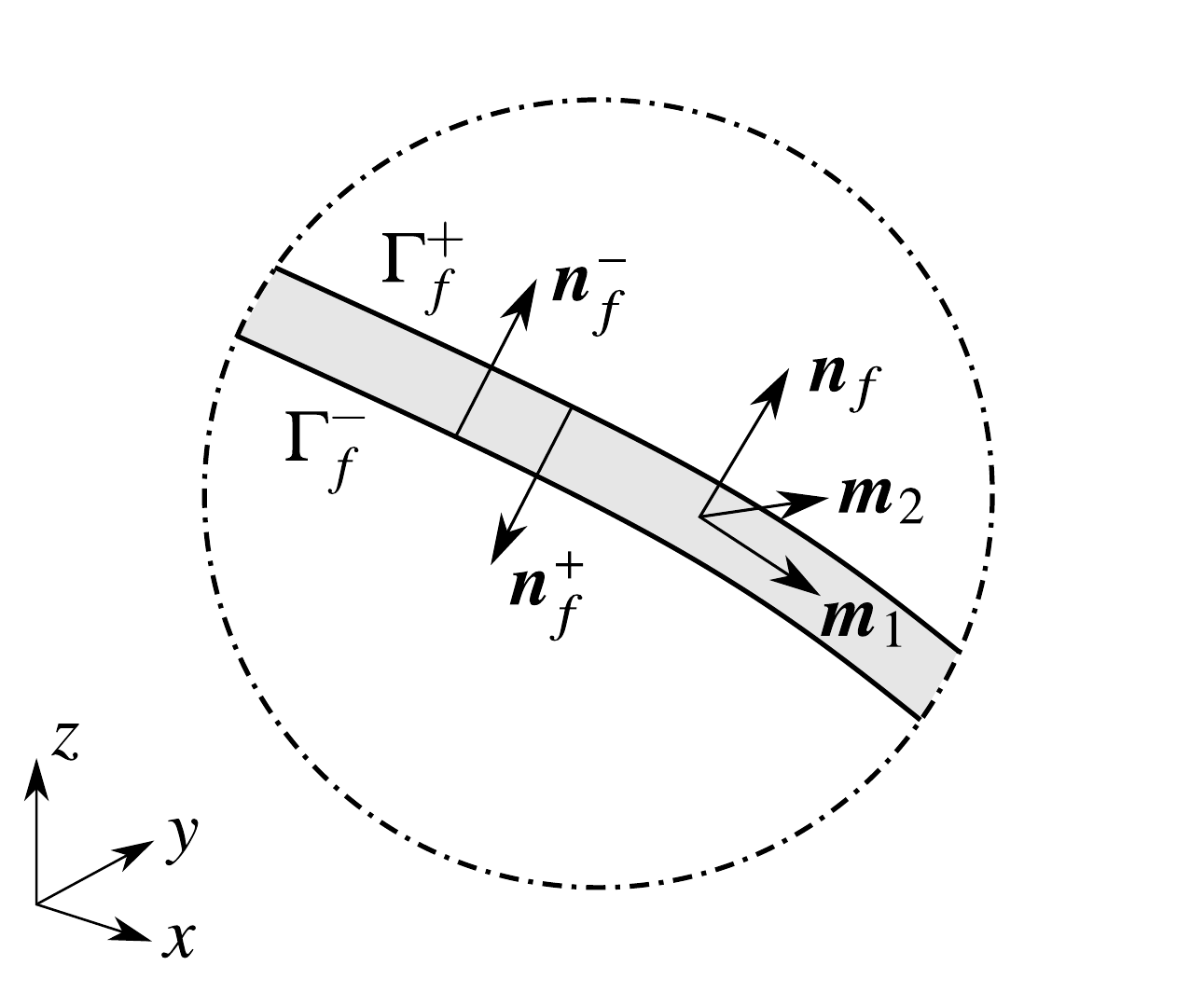}
  \caption{}
  \label{fig:conc_scheme2}
\end{subfigure}
\hfill
\begin{subfigure}{0.29\linewidth}
  \centering
  \includegraphics[width=\linewidth]{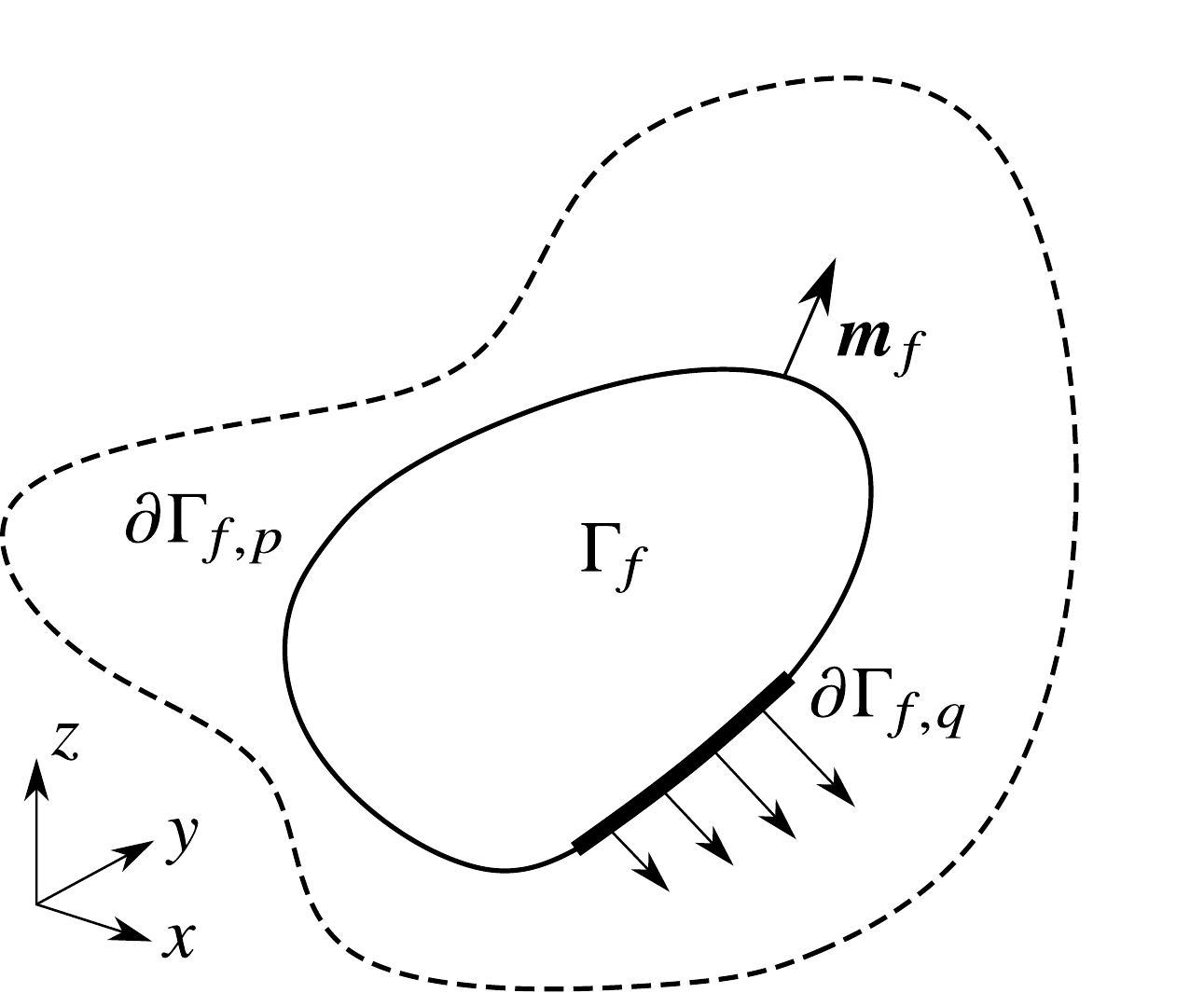}
  \caption{}
  \label{fig:conc_scheme3}
\end{subfigure}
\hfill\null
  \caption{Conceptual scheme for fracture modeling (a); local reference frame on the
    discontinuity surface (b); and the scheme for the flow domain (c).}
  \label{fig:conc_scheme}
\end{figure}

We consider an elastic, closed, polyhedral domain $\overline{\Omega} = \Omega \cup \partial\Omega$ in $\mathbb{R}^3$, with $\Omega$ an open set, $\partial \Omega$ its boundary, and $\vec{n}_{\Omega}$ its outer normal vector.
As usual, the boundary is subdivided into two non-overlapping subsets where
$\partial\Omega = \overline{\partial\Omega_u \cup \partial\Omega_{\sigma}}$ such that
$\partial\Omega_u \cap \partial\Omega_{\sigma} = \varnothing$, where Dirichlet and Neumann boundary conditions for displacement and traction are applied (see Fig. \ref{fig:conc_scheme1}).
From a mathematical standpoint, a fracture is described as an internal boundary $\Gamma_f$ embedded in $\Omega$, consisting of two overlapping surfaces $\Gamma_f^-$ and $\Gamma_f^+$ as shown in Fig. \ref{fig:conc_scheme2}.
The fracture orientation is characterized by a unit vector orthogonal to the fracture plane.
By convention, we choose $\vec{n}_f = \vec{n}_f^- = -\vec{n}_f^+$. 
On this lower-dimensional domain, the pressure field is defined.
Let $\overline{\Gamma}_f = \Gamma_f \cup \partial \Gamma_f$ denote the closed domain occupied by the fracture, with $\Gamma_f$ a two-dimensional (2D) open surface and $\partial \Gamma_f$ a one-dimensional (1D) curve defining its boundary.
The fracture boundary is subdivided into two non-overlapping subsets, $\partial \Gamma_f =
\overline{\partial \Gamma_{f,p} \cup \partial \Gamma_{f,q}}$ such that $\partial \Gamma_{f,p} \cap \partial \Gamma_{f,q} =
\varnothing$, corresponding to the position of Dirichlet and Neumann boundary conditions
for pressure and flux fields, respectively, as shown in Fig. \ref{fig:conc_scheme3}.
We define $\vec{m}_{f}$ as the outer normal vector for $\Gamma_f$.
Finally, let $\mathbb{T} = (0, t_{\max})$ denote the time domain of interest.

We assume quasi-static conditions and infinitesimal strains.
The fluid is taken to be incompressible, and we neglect body force and buoyancy effects.
The strong form of initial boundary value problem (IBVP) can then be stated as \cite{KikOde88,Lau03,Wri06}:
\begin{subequations}
\begin{align}
  \intertext{Find $\vec{u}:
\overline{\Omega} \times \mathbb{T} \rightarrow \mathbb{R}^3$ and $p: \overline{\Gamma}_f
\times \mathbb{T} \rightarrow \mathbb{R}$ such that}
    - \text{div} \;\tensorTwo{\sigma} (\vec{u}) &= \vec{0} & &\mbox{ in } \Omega \times \mathbb{T}
& &\mbox{(linear momentum balance)}, \label{eq:momentumBalanceS}\\
    \dot{g}_N(\vec{u}) + \text{div} \; \vec{q} (\tensorOne{u},p) &= q_s & &\mbox{ in } \Gamma_f \times \mathbb{T} & &\mbox{(mass balance)},
\label{eq:massBalanceS}\\
    \tensorTwo{\sigma}(\tensorOne{u}) \cdot \vec{n}_f - p \vec{n}_f &= \vec{0} & &\mbox{ on } \Gamma_f
\times \mathbb{T} & &\mbox{(traction balance on the fracture)},
\label{eq:momentumBalanceS_F}\\
    \vec{u} &= \bar{\vec{u}} & &\mbox{ on } \partial \Omega_{u} \times \mathbb{T} &
&\mbox{(prescribed boundary displacement)}, \label{eq:momentumBalanceS_DIR}\\
    \tensorTwo{\sigma}(\tensorOne{u}) \cdot \vec{n} &= \bar{\vec{t}} & &\mbox{ on }
\partial \Omega_{\sigma} \times \mathbb{T} & &\mbox{(prescribed boundary traction)},
\label{eq:momentumBalanceS_NEU}\\
    p &= \bar{p} & &\mbox{ on } \partial \Gamma_{f,p} \times \mathbb{T} & &\mbox{(prescribed
boundary pressure)}, \label{eq:massBalanceS_DIR}\\
    \vec{q} \cdot \vec{m}_f &= \bar{q} & &\mbox{ on } \partial \Gamma_{f,q} \times
\mathbb{T} & &\mbox{(prescribed boundary flux)}, \label{eq:massBalanceS_NEU}\\
    \vec{u}|_{t=0} &= \vec{u}_0 & &\mbox{ in } \overline{\Omega} & &\mbox{(initial
displacement)}, \label{eq:ICu}\\
    p|_{t=0} &= p_0 & &\mbox{ in } \overline{\Gamma}_f & &\mbox{(initial pressure)},
\label{eq:ICp}\\
   \intertext{subject to the constraints}
    t_N = \vec{t} \cdot \vec{n}_f &\le 0 & &\mbox{ on } \Gamma_f \times \mathbb{T} &
&\mbox{(normal contact conditions)}, \label{eq:normal_contact_KKT_1} \\
    g_N = \llbracket \vec{u} \rrbracket \cdot \vec{n}_f &\ge 0 & &\mbox{ on } \Gamma_f
\times \mathbb{T}, \label{eq:normal_contact_KKT_2} \\
    t_N g_N &= 0 & &\mbox{ on } \Gamma_f \times \mathbb{T},
\label{eq:normal_contact_KKT_3} \\
    \nonumber \\
    \left\| \vec{t}_T \right\|_2 - \tau_{\max}(t_N) & \le 0 & &\mbox{
on } \Gamma_f \times \mathbb{T}, & & \mbox{(Coulomb frictional contact conditions)},
\label{eq:frictional_contact_KKT_1_new} \\
    \dot{\vec{g}}_T \cdot \vec{t}_T - \tau_{\max}(t_N) \left\| \dot{\vec{g}}_T \right\|_2  &= 0 & &\mbox{
on } \Gamma_f \times \mathbb{T}. \label{eq:frictional_contact_KKT_2_new}
\end{align}
\label{eq:IBVP}\null
\end{subequations}
Here, known boundary and initial conditions are given as 
$\bar{\vec{u}}: \partial \Omega_{u} \times \mathbb{T} \rightarrow \mathbb{R}^3$,
$\bar{\vec{t}}: \partial \Omega_{\sigma} \times \mathbb{T} \rightarrow \mathbb{R}^3$,
$\vec{u}_0: \overline{\Omega} \rightarrow \mathbb{R}^3$, $q_s: \Gamma_f \times
\mathbb{T} \rightarrow \mathbb{R}$, $\bar{p}: \partial \Gamma_{f,p} \times \mathbb{T} \rightarrow
\mathbb{R}$, $\bar{q}: \partial \Gamma_{f,q} \times \mathbb{T} \rightarrow \mathbb{R}$, and
$p_0: \overline{\Gamma}_f \times \mathbb{T} \rightarrow \mathbb{R}$.
The following symbols, variables, and constitutive relationships are also introduced:
\begin{itemize}
  \item $\tensorTwo{\sigma}(\vec{u}) = \tensorFour{C}: \nabla^s \vec{u}$ is the Cauchy stress tensor, which is related to the displacement vector $\vec{u}$ by the fourth-order elasticity tensor $\tensorFour{C}$, with $\nabla^s$ the symmetric gradient operator;
  \item $\vec{q} (\vec{u},p) = -\frac{\displaystyle C_f(\vec{u})}{\displaystyle  \mu} \nabla p$ is the fluid volumetric flux in the fracture domain---assuming laminar flow and validity of Darcy's law \cite{witherspoon1980validity}---with $\nabla p$ the fluid pressure gradient, $\mu$ the fluid viscosity (constant), and $C_f$ the isotropic fracture hydraulic conductivity modeled as in \cite{GarKarTch16}:
\begin{equation}
  C_f = C_{f,0} + \frac{g_N^3}{12}.
  \label{eq:cond_def}
\end{equation}
Note that $C_{f,0}$ captures the conductivity associated with two irregular surfaces that are
in contact \cite{kamenov2013laboratory}.
This allows fluid to flow and pressure to propagate even if the
fracture is nominally closed. From a physical viewpoint,
the volume between two rough surfaces in contact is non-zero and fluid may infiltrate between asperities. For simplicity, here we assume a constant closed conductivity, though for some applications a normal-stress dependent model may be preferred.
  \item $\vec{t} = \tensorTwo{\sigma} \cdot \vec{n}_f^- = - \tensorTwo{\sigma} \cdot
\vec{n}_f^+ = (t_N \vec{n}_f + \vec{t}_T)$ is the traction vector over $\Gamma_f$, with
$t_N$ and $\vec{t}_T = (t_{m_1} \vec{m}_1 + t_{m_2} \vec{m}_2)$ its normal and tangential
component, respectively, with respect to the local reference system shown in Fig. \ref{fig:conc_scheme2};
  \item $\tau_{\text{max}} = c - t_N \tan(\theta)$ is the limit value for the modulus of $\vec{t}_T$ according to the Coulomb criterion, with $c$ and $\theta$ the cohesion and friction angle, respectively;
  \item $\llbracket \ast \rrbracket$ denotes the jump of a quantity across $\Gamma_f$, with
$\llbracket \vec{u} \rrbracket = ( \left.\vec{u}\right|_{ \Gamma_f^+ } -
\left.\vec{u}\right|_{ \Gamma_f^- } ) = (g_N \vec{n}_f + \vec{g}_T)$ the relative
displacement across $\Gamma_f$, where $g_N$ and $\vec{g}_T$ are normal and tangential
components, respectively, and $\left.\vec{u}\right|_{ \Gamma_f^+ }$ and
$\left.\vec{u}\right|_{ \Gamma_f^- }$ are the restrictions of $\vec{u}$ on $\Gamma_f^+$ and
$\Gamma_f^-$.
\end{itemize}
For additional details regarding the governing formulation, we refer the reader to
\cite{KikOde88,Lau03,Wri06}.

\begin{rem}
In the literature, other constitutive laws are available that relate opening
and fracture conductivity. For the laminar regime, they are based on
improvements of the original Poiseuille law \cite{johnson2016handbook}, describing flow between parallel smooth surfaces. Corrections are introduced to address surface irregularities. In \eqref{eq:cond_def}, $C_{f,0}$ is used to quantify the surface
roughness, while other authors used a multiplicative coefficient $1/f$
\cite{witherspoon1980validity,schrefler2006adaptive} or a hydraulic aperture, different
from the nominal aperture $g_N$ \cite{zimmerman1996hydraulic}. For the turbulent regime,
the constitutive behavior can be based on more complex hydrodynamic relationships, like the
Darcy-Weisback equation and Moody's diagram \cite{johnson2016handbook}. See also
\cite{vahab2017numerical} for additional discussion on this topic.
\end{rem}

Throughout this work, we will use subscripts $N$ and $T$ to identify the normal and tangential components of a vector quantity with respect to the discontinuity $\Gamma_f$.
Specifically, given the vector $\tensorOne{v} \in \mathbb{R}^3$, we have
\begin{align}
  \tensorOne{v} &= v_N \tensorOne{n}_f + \tensorOne{v}_T, &
  v_N &= \tensorOne{v} \cdot \tensorOne{n}_f, &
  \tensorOne{v}_T &= (\tensorTwo{1} - \vec{n}_f \otimes \vec{n}_f) \cdot \tensorOne{v},
  \label{eq:vec_NTcomp}
\end{align}
with $\otimes$ the dyadic product.
Note that we consider the static Coulomb law.
Therefore, in a discrete setting the tangential velocity $\dot{\tensorOne{g}}_T$ in \eqref{eq:frictional_contact_KKT_2_new} can be replaced with the tangential displacement increment \cite{wohlmuth2011variationally}.
This has to be done at every timestep if an implicit time-marching scheme is used.
From now on, $\dot{\vec{g}}_T$ will be replaced by $\Delta_n \vec{g}_T$, where $\Delta_n$ denotes the backward difference operator such that $\Delta_n (\ast) = (\ast)_{n} - (\ast)_{n-1}$ with the subscript $n$ denoting the current discrete time level $t_n$.
The Coulomb frictional contact conditions Eqs. \eqref{eq:frictional_contact_KKT_1_new}--\eqref{eq:frictional_contact_KKT_2_new} can then be rewritten as:
\begin{subequations}
\begin{align}
  || \vec{t}_{T,n} ||_2 < \tau_{\max}
  &\Longrightarrow
  \Delta_n \vec{g}_T = \vec{0},
  \label{eq:frictional_contact_1} \\
  || \vec{t}_{T,n} ||_2 = \tau_{\max}
  &\Longrightarrow
  \tensorOne{t}_{T,n} =
  \tensorOne{t}_{T}^{*}(\tensorOne{u}_n, t_{N,n}) =
  \tau_{\max}(t_{N,n}) \frac{\Delta_n {\vec{g}}_{T} }{|| \Delta_n {\vec{g}}_{T} ||_2}.
  \label{eq:frictional_contact_2}
\end{align}
\end{subequations}

In our framework, the fracture is explicitly modeled according to a Discrete Fracture
Model \cite{karimi2003efficient,GarKarTch16}, with $\Gamma_f$ encompassing the whole region
where opening or contact may take place at any $t \in \mathbb{T}$ \cite{Lau03}.
That is, we assume $\Gamma_f$ is fixed and does not propagate.  The only unknown is then its partitioning into stick,
slip and open patches at a given point in time, $\Gamma_f = \Gamma_f^{\text{stick}} \cup \Gamma_f^{\text{slip}} \cup
\Gamma_f^{\text{open}}$.
These three modes are associate with behavior regimes, namely:
\begin{itemize}
  \item \textbf{Stick mode} on $\Gamma_f^{\text{stick}}$: the discontinuity is fully
closed and compressed with the Coulomb criterion satisfied, i.e., $t_N < 0$ (Eq.
\eqref{eq:normal_contact_KKT_1}) and $\left\| \vec{t}_T \right\|_2 < \tau_{\max}(t_N)$ (Eq. \eqref{eq:frictional_contact_KKT_1_new}).
The three components of the traction are unknown and such that no relative movement is
allowed between $\Gamma_f^+$ and $\Gamma_f^-$. The conductivity is constant and equal to
$C_{f,0}$. This implies linear, steady-state flow behavior within the stick portion.
  \item \textbf{Slip mode} on $\Gamma_f^{\text{slip}}$: the discontinuity is compressed, but a
slip displacement $\vec{g}_T$ between $\Gamma_f^+$ and $\Gamma_f^-$ is allowed for. Only
the normal traction component $t_N$ is unknown, the tangential component having magnitude
$|| \vec{t}_T^* ||_2 = \tau_{\text{max}}(t_N)$ and direction collinear with
$\dot{\vec{g}}_T$.  The conductivity behavior is the same as for the stick mode, as no
slip-induced dilation is modeled. The flow follows, as before, a linear, steady
state model.
  \item \textbf{Open mode} on $\Gamma_f^{\text{open}}$: $\Gamma_f^+$ and $\Gamma_f^-$ are
not in contact and a free relative displacement $\llbracket \vec{u} \rrbracket$ is
allowed for. The traction is known and equal to the zero vector in $\mathbb{R}^3$. In this
case, the conductivity is related to the opening (see Eq. \eqref{eq:cond_def}) and the
fluid behavior can be described as a nonlinear transient flow, with linear storage and
nonlinear conductivity.
\end{itemize}
The numerical strategy to resolve this partitioning is described in the next section.

\section{Numerical model}
\label{sec:numerical_model}

\subsection{Discrete formulation}
\label{sec:variational_form_new}

We solve numerically the IBVP \eqref{eq:IBVP} using a mixed finite-element/finite-volume approach for the spatial discretization combined with a fully-implicit time-marching scheme.
In particular, the contact mechanics problem is addressed using the saddle point formulation \cite{KikOde88,Lau03,Wri06} where traction vectors acting on $\Gamma_f$ are introduced as additional primary variables serving as Lagrange multipliers to enforce normal and frictional contact constraints.
The simulation of the fluid flow through the discontinuity relies on a finite volume method \cite{EymGalHer00}.

We introduce a triangulation $\mathcal{T}$ of the domain consisting of nonoverlapping hexahedral cells that conforms to the discontinuity surface, i.e. $\overline{\Omega} = \bigcup_{\tau \in \mathcal{T}} \overline{\tau}$.
Let us define $\mathcal{F}_f$ as the set of quadrilateral faces $\varphi$ in the triangulation such that $\overline{\Gamma}_f = \bigcup_{\varphi \in \mathcal{F}_f} \overline{\varphi}$.
The edge normal unit vector to face $\varphi_K$ in the discontinuity plane is denoted by $\tensorOne{m}_{K}$, with $K$ the face global index.
Let $\mathcal{E}_f = (\mathcal{E}_{\text{int}} \cup \mathcal{E}_p \cup \mathcal{E}_q)$ be the set of edges belonging to faces defining $\overline{\Gamma}_f$, with $\mathcal{E}_{\text{int}}$ (respectively $\mathcal{E}_p$, $\mathcal{E}_q$) the set of edges included in $\Gamma_f$ (respectively $\partial \Gamma_{f,p}$, $\partial \Gamma_{f,q}$).
An edge in $\mathcal{E}_{\text{int}}$ shared by faces $\varphi_K$ and $\varphi_L$ is denoted as $\varepsilon_{K,L}$, with face global indices such that $K<L$.
Similarly, an edge in $\mathcal{E}_p \cup \mathcal{E}_q$ belongs to a unique face $\varphi_K$ and is simply denoted as $\varepsilon_{K}$.
A unique orientation on $\Gamma_f$ is associated with each edge in $\mathcal{E}_f$ through a unit vector $\tensorOne{m}_{\varepsilon}$.
We set $\tensorOne{m}_{\varepsilon} = \tensorOne{m}_{K}$ for any edge $\varepsilon_{K,L} \in \mathcal{E}_{\text{int}}$ and $\varepsilon_{K} \in \mathcal{E}_{p} \cup \mathcal{E}_{q}$.

The following discrete approximation spaces are used for the displacement and traction field \cite{KikOde88,wohlmuth2011variationally}
\begin{subequations} \label{eq:func_discrete_spaces}
\begin{align}
  \vecFunSpace{U}^h &:=
  \left\lbrace
    \tensorTwo{\eta} \in [C^0(\overline{\Omega})]^3 :
	  \tensorTwo{\eta} = \bar{\vec{u}} \text{ on } \partial \Omega_u, 
  	\tensorTwo{\eta}\vert_{\tau} \in \vec{Q}_1(\tau) \; \forall \tau \in \mathcal{T}
  \right\rbrace \subset \vecFunSpace{U}, \\
  \vecFunSpace{U}_0^h &:=
  \left\lbrace
    \tensorTwo{\eta} \in [C^0(\overline{\Omega})]^3 :
    \tensorTwo{\eta} = \vec{0} \text{ on } \partial \Omega_u,
    \tensorTwo{\eta}\vert_{\tau} \in \vec{Q}_1(\tau) \; \forall \tau \in \mathcal{T}
  \right\rbrace \subset \vecFunSpace{U}_0, \\
  \vecFunSpace{M}^h(t^h_{N,n}) &:=
    \left\lbrace
    \tensorOne{\mu} \in [L^2(\Gamma_f)]^3 :
    \tensorOne{\mu}\vert_{\varphi} \in \vec{P}_0(\varphi) \; \forall \varphi \in \mathcal{F}_f,
    \mu_N \le 0,
    ||\tensorOne{\mu}_T||_2 \le \tau_{\max}(t^h_{N,n})
  \right\rbrace,	\subset \vecFunSpace{M}, 
\end{align}
\end{subequations}
with $\vecFunSpace{U}$ and $\vecFunSpace{U}_0$ the space of vector functions in $[L^2(\Omega)]^3$ with first derivatives in $L^2(\Omega)$ having trace on $\partial \Omega_u$ equal to $\bar{\vec{u}}$ and $\vec{0}$, respectively, and $\vecFunSpace{M}$ the dual space of the trace space $\vecFunSpace{W}$ of $\vecFunSpace{U}_0$ restricted to $\Gamma_f$.
Here, $C^0(\overline{\Omega})$ and $L^2(\Omega)$ denote the space of continuous and square Lebesgue-integrable functions on $\overline{\Omega}$ and $\Omega$, $\vec{Q}_1(\tau) = [Q_1(\tau)]^3$, $Q_1(\tau) = \text{span}\{1, x, y, z, xy, xz, yz, xyz \}$---namely, the mapping to $\tau$ of the space of trilinear polynomials on the unit cube $[0,1]^3$ in $\mathbb{R}^3$---, and $\vec{P}_0 (\tau) = [P_0(\tau)]^3$, with $P_0(\tau)$ the space of piecewise constant functions.
Hereafter, for a given domain $D$, the compact notation $(*,\star)_D$ is used to denote the $L^2$-inner product of scalar, vector, or rank-2 tensor functions in $L^2(D)$, $[L^2(D)]^3$, or $[L^2(D)]^{3 \times 3}$, as appropriate---for example, $( \nabla^s \tensorOne{ \eta }, \tensorTwo{\sigma} )_{\Omega} = \int_{\Omega} \nabla^{s} \tensorTwo{\eta} : \tensorTwo{\sigma} \, \mathrm{d}V$.

The pressure field is approximated in the space of piecewise constant functions on faces in $\mathcal{F}_f$, namely
\begin{align} \label{eq:FV_func_discrete_spaces}
  \mathcal{P}^h &:=
  \left\lbrace
    \chi \in L^2(\Gamma_f) :
    \chi\vert_{\varphi} \in P_0 (\varphi) \; \forall \varphi \in \mathcal{F}_f
  \right\rbrace.
\end{align}
In this work the discretization of the mass balance equation \eqref{eq:massBalanceS} is based on the classical two point flux approximation (TPFA) scheme.
To allow for a unified presentation of the coupled mixed finite-element/finite-volume model, the TPFA scheme is written in weak form \cite{EymGalHer00,EymGalHer07,Age_etal10}.
Therefore, for the space $\mathcal{P}^h$ we introduce the following weighted inner product  
\begin{align}
  [\chi, p^h]_{\mathcal{F}_f} = 
  \smashoperator{\sum_{\varepsilon = \varepsilon_{K,L} \in \mathcal{E}_{\text{int}} }}
  ( \chi\vert_{ \varphi_L } - \chi\vert_{ \varphi_K } ) \Upsilon_{KL} ( p^h\vert_{ \varphi_L } - p^h\vert_{ \varphi_K } )
  + \smashoperator{\sum_{\varepsilon = \varepsilon_{K} \in \mathcal{E}_p }}
  \chi\vert_{ \varphi_K } \Upsilon_K p^h\vert_{ \varphi_K },
  \label{eq:FV_inner_product}
\end{align}
where $\Upsilon_{KL} = (\Upsilon_{K} \Upsilon_{L}) / (\Upsilon_{K} + \Upsilon_{L})$ is the harmonic average of positive one-sided transmissibility $\Upsilon_{K}$ and $\Upsilon_{L}$ associated to face $\varphi_K$ and $\varphi_L$, respectively.
In particular, the one-sided transmissibility is expressed as the product of a nonlinear ($C_{f, \, \beta }$) and a constant ($\overline{\Upsilon}_{\beta}$) term, respectively
\begin{align}
  \Upsilon_{\beta} &= C_{f, \, \beta } \overline{\Upsilon}_{\beta},
  &
  C_{f, \, \beta }
  &= C_{f,0}
  + \frac{1}{| \varphi_{\beta} |} \int_{\varphi_{\beta} } \frac{ g_{N,n}^3 }{12} \, \mathrm{d}A,
  &
  \overline{\Upsilon}_{\beta} &= \frac{|\varepsilon|}{\mu} \frac{(\tensorOne{x}_{\varepsilon} - \tensorOne{x}_{\beta}) \cdot \tensorOne{m}_{\beta}}{||\tensorOne{x}_{\varepsilon} - \tensorOne{x}_{\beta} ||_2^2},
  &
  \beta &= \{K, L \},
  \label{eq:T_def}
\end{align}
with $|\varepsilon|$ the edge lenght, $\tensorOne{x}_{\varepsilon}$ a collocation point associated with each edge in $\mathcal{E}^f$, and $\tensorOne{x}_{\beta}$ and $| \varphi_{\beta} |$ the barycenter and area of
face $\varphi_{\beta}$.
Notice that $C_{f, \, \beta }$ is the mean value of the fracture conductivity over $\varphi_{\beta}$.
The dependence of the hydraulic conductivity on normal displacement jump makes the transition from a closed to open state even more challenging.
\begin{rem}
In \eqref{eq:FV_inner_product}, the term $\hat{q}_{K,L} = - \Upsilon_{KL} ( p^h\vert_{ \varphi_L } - p^h\vert_{ \varphi_K } )$ represents the numerical flux approximation based on the TPFA scheme \cite{EymGalHer00,barth2018finite} to the volume flow rate per unit length through $\varepsilon = \varepsilon_{K,L} \in \mathcal{E}_{\text{int}}$, i.e $\hat{q}_{K,L} \approx \int_{\varepsilon} \tensorOne{q} \cdot \tensorOne{m}_{\varepsilon} \;\text{d}l$.
Similarly, for a boundary edge $\varepsilon = \varepsilon_{K} \in \mathcal{E}_p$ in the presence of homogeneous pressure conditions, such an approximation is given by $\hat{q}_{K} = - \Upsilon_{K} p^h\vert_{ \varphi_K }$.
\end{rem}
\begin{rem}
In the derivation of transmissibility coefficients, the collocation point $\tensorOne{x}_{\varepsilon}$ in \eqref{eq:T_def} is used to impose point-wise pressure continuity between faces sharing an edge.
Different strategies can be used to choose $\tensorOne{x}_{\varepsilon}$.
Following \cite{KarDur16}, in our implementation we select $\tensorOne{x}_{\varepsilon}$ as the intersection of the edge $\varepsilon_{K,L}$ and the line connecting the barycenters of faces $\varphi_K$ and $\varphi_L$.
For a boundary edge $\varepsilon_K \in \mathcal{E}_p \cup \mathcal{E}_q$, $\tensorOne{x}_{\varepsilon}$ is chosen as the  orthogonal projection of $\tensorOne{x}_K$.
\end{rem}
\begin{rem}
Because of its robustness, simplicty and monotonicity-preserving properties, the classical TPFA scheme is the  method-of-choice in industrial reservoir simulation.
However, TPFA may lead to an inconsistent numerical flux on arbitrary grids \cite{EymGalHer00,barth2018finite}, for which more sophisticated discretization methods like multipoint and/or nonlinear schemes should be considered, see~\cite{droniou2014finite,TerMalTch17,Sch_etal18} and references therein.  Note that the stabilization scheme described below will work without modification for these alternative flux approximations.
\end{rem}
Discretizing the time interval $\mathbb{T}$ into $n_{\mathbb{T}}$ subintervals of size $\Delta_n t = (t_n - t_{n-1})$, the mesh-dependent fully discrete weak form of \eqref{eq:IBVP} is: find $\{ \vec{u}^h_n, \tensorOne{t}^h_n, p^h_n \} \in \vecFunSpace{U}^h \times \vecFunSpace{M}^h(t^h_{N,n}) \times \funSpace{P}^h$ such that for all $\{ \tensorOne{\eta}, \tensorOne{\mu}, \chi \} \in \vecFunSpace{U}_0^h \times \vecFunSpace{M}^h(t^h_{N,n}) \times \funSpace{P}^h$
\begin{subequations} \label{eq:weak_form_discr}
\begin{align}
  \mathcal{R}_{u} =
  & ( \nabla^s \tensorOne{ \eta }, \tensorTwo{\sigma}_n )_{\Omega}
  + ( \llbracket \tensorOne{ \eta } \rrbracket, \tensorOne{ t }_{n}^h - p_n^h \tensorOne{n}_f )_{\Gamma_f}
  - ( \tensorOne{\eta}, \bar{\vec{t}}_n )_{\partial \Omega_{\sigma}}
  = 0,
	\label{eq:weak_form_mom_h} \\
	\mathcal{R}_t =
  & ( t_{N,n}^h - \mu_N, g_{N,n} )_{\Gamma_f}
  + (\vec{t}_{T,n}^h - \tensorTwo{\mu}_T, \Delta_n {\vec{g}}_{T} )_{\Gamma_f}
  \ge 0,
  \label{eq:weak_lam_h}\\
  \mathcal{R}_{p} =
  & \left( \chi, \frac{\Delta_n g_{N} }{\Delta_n t} \right)_{\Gamma_f}
  + [ \chi, p_n^h ]_{\mathcal{F}_f} 
  - F_{\mathcal{F}_f}(\chi)
  + G_{\mathcal{F}_f}(\chi)
  - ( \chi, q_{s,n} )_{\Gamma_f} = 0.
  \label{eq:weak_form_mass_h}
\end{align}
\end{subequations}
with $n \in \{1, 2, \ldots, n_{\mathbb{T}} \}$, and $F_{\mathcal{F}_f}(\chi)$ and $G_{\mathcal{F}_f}(\chi)$ two functionals that incorporate boundary conditions \eqref{eq:massBalanceS_DIR}-\eqref{eq:massBalanceS_NEU}
\begin{align}
  F_{\mathcal{F}_f}(\chi) &= \smashoperator{\sum_{\varepsilon = \varepsilon_{K} \in \mathcal{E}_p }} \chi\vert_{\varphi_K } \Upsilon_K \frac{1}{|\varepsilon|}\int_{\varepsilon} \bar{p}_{n} \, \mathrm{d}l,
  &
  G_{\mathcal{F}_f}(\chi) &= \smashoperator{\sum_{\varepsilon = \varepsilon_{K} \in \mathcal{E}_q }} \chi\vert_{\varphi_K } \int_{\varepsilon} \bar{q}_{n} \, \mathrm{d}l.
\end{align}

Let $\{ \tensorOne{\eta}_i \}_{i \in \mathcal{N}_u \cup \overline{\mathcal{N}}_u}$ be the standard vector nodal basis functions for the global finite element space of continuous piecewise-$[Q_1]^3$ functions associated with $\mathcal{T}$, with $\mathcal{N}_u$ and $\overline{\mathcal{N}}_u$ the set of indices of basis function vanishing on $\partial \Omega_u$ and having support on $\partial \Omega_u$, respectively.
Note that $\text{card}(\mathcal{N}_u) + \text{card}(\overline{\mathcal{N}}_u)$ is equal to three times the number of vertices in $\mathcal{T}$.
Let $\{ \chi_k \}_{k \in \mathcal{N}_p}$ be the basis for $\funSpace{P}^h$, with $\chi_k$ the characteristic function of the $k$th face in $\mathcal{F}_f$ such that $\chi_k(\tensorOne{x}) = 1$, if $\tensorOne{x} \in \varphi_k$, $\chi_k(\tensorOne{x}) = 0$, if $\tensorOne{x} \notin \varphi_k$, and $\mathcal{N}_p = \{1,\ldots, \text{card}(\mathcal{F}_f) \}$.
Let $\{ \tensorOne{\mu}_j \}_{j \in \mathcal{N}_t}$ be the piecewise-constant vector basis for $\funSpace{M}^h(t_N)$, with $\mathcal{N}_t = \{1,\ldots, 3\cdot\text{card}(\mathcal{F}_f) \}$--namely, the three basis functions $\chi_j(\tensorOne{x}) \tensorOne{n}_f(\tensorOne{x}_j)$, $\chi_j(\tensorOne{x}) \tensorOne{m}_1(\tensorOne{x}_j)$, and $\chi_j(\tensorOne{x}) \tensorOne{m}_2(\tensorOne{x}_j)$ associated to each face $\varphi_j \in \mathcal{F}_f$.
Discrete approximations for the displacement, traction, and pressure can then be expressed as
\begin{align}
  \vec{u}^h_n(\tensorOne{x}) &=
  \sum_{i \in \mathcal{N}_u } \tensorOne{\eta}_i(\tensorOne{x}) u_{i,n} +
  \sum_{i \in \overline{\mathcal{N}}_u } \tensorOne{\eta}_i(\tensorOne{x}) \bar{u}_{i,n},
  &
  \vec{t}^h_n(\tensorOne{x}) &=
  \sum_{j \in \mathcal{N}_t } \tensorOne{\mu}_j(\tensorOne{x}) t_{j,n},
  &
  p^h_n(\tensorOne{x}) &=
  \sum_{k \in \mathcal{N}_{p} } {\chi_k}(\tensorOne{x}) p_{k,n}.
  \label{eq:u_t_q_h}
\end{align}
The unknown nodal displacement components $\{u_{i,n}\}$, face-centered traction components $\{ t_{j,n} \}$, and face-centered pressures $\{ p_{k,n} \}$ at time level $t_n$ are gathered in algebraic vectors $\Vec{u}_n$, $\Vec{t}_n$, and $\Vec{p}_n$.
We emphasize that, at the right-hand side of the expression for $\tensorOne{u}^h$, the first sum represents an approximate displacement solution of the IBVP \eqref{eq:IBVP} satisfying homogeneous prescribed displacement, whereas the second sum is the discrete extension by zero (to the degrees of freedom) of the boundary datum $\bar{\tensorOne{u}}$ over $\partial \Omega_u$.
Hence, $\{ \tensorOne{\eta}_i \}_{i \in \mathcal{N}_u }$ is a basis for $\vecFunSpace{U}^h_0$.

\subsection{Solution strategy}
\label{sec:nonlinear_form_new}

The discrete form of the IBVP \eqref{eq:IBVP} based on the weak form \eqref{eq:weak_form_discr} consists of a nonlinear system of equations and inequalities in the unknowns $\Vec{u}_n$, $\Vec{t}_n$, and $\Vec{p}_n$.
However, if we postulate that active contact regions $\Gamma_{f,n}^{\text{stick}}$ and $\Gamma_{f,n}^{\text{slip}}$ are known at $t_n$, the inequality \eqref{eq:weak_lam_h} can be replaced with a variational equality:
\begin{align}
  \mathcal{R}_t =
	& ( \mu_N, g_{N,n} )_{\Gamma_{f,n}^{\text{stick}} \cup \Gamma_{f,n}^{\text{slip}} } 
  + ( \tensorTwo{\mu}_T, \Delta_n {\vec{g}}_{T} )_{\Gamma_{f,n}^{\text{stick}}} 
  + ( \tensorTwo{\mu}_T, \tensorOne{t}^h_{T,n} - \tensorOne{t}_{T,n}^{*} )_{\Gamma_{f,n}^{\text{slip}}}
  + ( \tensorTwo{\mu}, \tensorTwo{t}^h_n )_{\Gamma_{f,n}^{\text{open}}}
  = 0,
  \label{eq:weak_lam_h_partitioned}
\end{align}
where the four integrals correspond to the weak enforcement of the: (i)  impenetrability condition on $\Gamma_{f,n}^{\text{stick}} \cup \Gamma_{f,n}^{\text{slip}}$, (ii) no slip condition on $\Gamma_{f,n}^{\text{stick}}$, (iii) known tangential traction on $\Gamma_{f,n}^{\text{slip}}$ given by Eq. \eqref{eq:frictional_contact_KKT_2_new}, and (iv) zero traction condition on $\Gamma_{f,n}^{\text{open}}$.
Introducing expressions \eqref{eq:u_t_q_h} into \eqref{eq:weak_form_mom_h}, \eqref{eq:weak_lam_h_partitioned}, and \eqref{eq:weak_form_mass_h} then allows for writing the discrete problem as a nonlinear system of equations that can be solved for the latest solution vectors $\Vec{u}_n$, $\Vec{t}_n$, and $\Vec{p}_n$, 
\begin{align}
  \begin{dcases}
    \Vec{r}_{u} (\Vec{u}_n, \Vec{t}_n, \Vec{p}_n)     = \Vec{0}, \\
    \Vec{r}_{t} (\Vec{u}_n, \Vec{u}_{n-1}, \Vec{t}_n) = \Vec{0}, \\
    \Vec{r}_{p} (\Vec{u}_n, \Vec{u}_{n-1}, \Vec{p}_n) = \Vec{0},
  \end{dcases}
  \label{eq:res_eqs}
\end{align}
with $\Vec{u}_{n-1}$ the known discrete displacement solution from the previous timestep.
To advance one timestep, since the partition of $\Gamma_{f,n}$ into $\Gamma_{f,n}^{\text{stick}}$, $\Gamma_{f,n}^{\text{slip}}$ and $\Gamma_{f,n}^{\text{open}}$ at $t_n$ is of course not known \emph{a priori}, an iterative procedure that includes a Newton method-based solver applied to residual equations \eqref{eq:res_eqs} is needed.

\begin{algorithm}[b]
\caption{Active-set strategy}
\begin{algorithmic}[1]
  \State $\ell$ = 0 \Comment{initialize iterator count}
  \If {$t_n > t_0$} \label{alg:activeset_init_str}
    \State $\Gamma_{f,n}^{\ell} = \Gamma_{f,(n-1)}$ \Comment{set Lagrange multipliers to previous timestep status}
  \Else
    \State $\Gamma_{f,n}^{\ell} = \Gamma_{f,0}^{\text{stick}}$ \Comment{set initial discontinuity status to \textit{stick}}
  \EndIf \label{alg:activeset_init_end}
  \Repeat
  \State $\ell \leftarrow \ell + 1$ \Comment{update active set iteration count}
  \State solve nonlinear problem \eqref{eq:res_eqs} for $\Gamma_{f,n}^{(\ell-1)}$ status
    \label{alg:line:solve}
  \State set \textit{stick}, \textit{slip}, and \textit{open} portions for $\Gamma_{f,n}^{\ell}$  \Comment{update status using last solution} \label{alg:line:update}
  \Until {($\Gamma_{f,n}^{\ell,\text{stick}} = \Gamma_{f,n}^{(\ell-1),\text{stick}}$ and
    $\Gamma_{f,n}^{\ell,\text{slip}} = \Gamma_{f,n}^{(\ell-1),\text{slip}}$ and $\Gamma_{f,n}^{\ell,\text{open}} = \Gamma_{f,n}^{(\ell-1),\text{open}}$)}
\end{algorithmic}
\label{alg:activeset}
\end{algorithm}

\begin{figure}[t]
  \centering
  \includegraphics[width=0.5\linewidth]{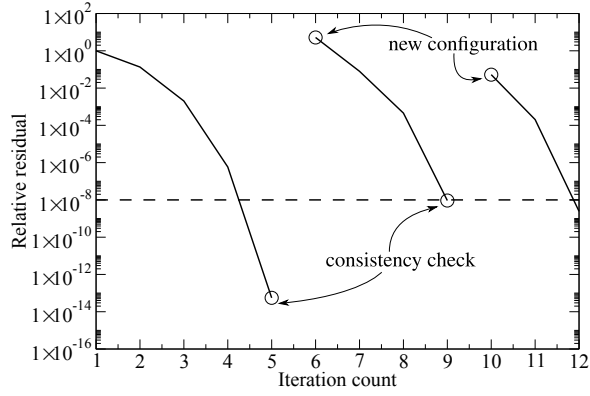}
  \caption{Example of convergence profile using the active set strategy. It is clear the
    jump in the residual after each consistency check.}
  \label{fig:activeSet}
\end{figure}

In this work, we solve the contact problem using an \textit{active-set} strategy, a numerical optimization technique employed in quadratic programming \cite{nocedal2006numerical,antil2018frontiers}.
Algorithm \ref{alg:activeset} summarizes the sequence of steps of an active-set algorithm applied to the contact problem.
From a practical viewpoint, first we assign an initial status (lines \ref{alg:activeset_init_str}--\ref{alg:activeset_init_end}) to all the Lagrange multipliers and solve the discrete nonlinear problem \eqref{eq:res_eqs}.
We highlight that there are just two states, \textit{active} and \textit{inactive}, for both normal and frictional contact.
The initial state is the previous time step solution, whenever available, or the stick state, i.e., $\Gamma_{f,0}^{\text{stick}} = \Gamma_f$, at the beginning of the simulation.
Once the nonlinear problem is solved, the status of all Lagrange multipliers is checked and a new subdivision into stick, slip, and open portions is identified.
If at least one multiplier changes status, we solve a new nonlinear system and check the consistency of the outcome again.
The procedure stops when the starting subdivision of $\Gamma_f$ is consistent with the solution of the nonlinear problem.  Figure~\ref{fig:activeSet} provides an illustration of the resulting nonlinear convergence profile, with three consistency checks required before the final accepted solution.

\begin{rem}
We note that there are no theoretical guarantees for the active-set convergence of
the discrete version of the IBVP \eqref{eq:IBVP}. It can occur that the exact solution
may lie between two discrete solutions, i.e., two similar but different subdivisions of
$\Gamma_f$. In these cases, one solution must be chosen.  In this work, we select the one with
the smaller internal energy.
\end{rem}

In Algorithm \ref{alg:activeset}, Newton's method is used to drive the norm of the combined residual vector below a specified relative tolerance (line \ref{alg:line:solve}).
To better identify each contribution associated with the contact process in the linearized problem at every iteration $\ell$ of the active-set algorithm, we further partition the vector containing the unknown traction degrees of freedom as $\Vec{t}_n^{\ell} = \{\Vec{t}_{S,n}^{\ell}, \Vec{t}_{N,n}^{\ell}, \Vec{t}_{T,n}^{\ell}, \Vec{t}_{O,n}^{\ell} \}$.
Subscripts $S$ and $O$ denote traction degrees of freedom belonging to $\Gamma_{f,n}^{\ell,\text{stick}} $ and $\Gamma_{f,n}^{\ell,\text{open}}$, respectively.
Subscripts $N$ and $T$ identify normal and tangential traction degrees of freedom on $\Gamma_{f,n}^{\ell,\text{slip}}$.
The set of indices of traction basis functions $\mathcal{N}_t$ is also consistently expressed as union of disjoint sets such that $\mathcal{N}_t = \mathcal{N}_t^{\ell,S} \cup \mathcal{N}_t^{\ell,N} \cup \mathcal{N}_t^{\ell,T} \cup \mathcal{N}_t^{\ell,O}$.
The solution of the nonlinear problem \eqref{eq:res_eqs} is then computed as follows.
Given an initial guess for displacement ($\Vec{u}_{n}^{\ell,(0)}$), traction ($\Vec{t}_{\beta,n}^{\ell,(0)}, \beta \in \{S, N, T, O \}$), and pressure ($\Vec{p}_{n}^{\ell,(0)}$) vectors, for $k = 0, 1, \ldots$ until convergence
\begin{align}
  \text{\underline{Solve}} \qquad &
  \begin{bmatrix}
    %
    %
    A_{uu} &
    A_{uS} &
    A_{uN} &
    A_{uT} &
    A_{uO} &
    A_{up} \\
    %
    %
    %
    A_{Su} & 0 & 0 & 0 & 0 & 0 \\
    %
    %
    A_{Nu} & 0 & 0 & 0 & 0 & 0 \\
    %
    %
    A_{Tu} &
    0 &
    A_{TN} &
    A_{TT} &
    0 &
    0 \\
    %
    %
    0 & 0 & 0 & 0 & A_{OO} & 0 \\
    %
    %
    A_{pu} & 0 & 0 & 0 & 0 & A_{pp} \\
  \end{bmatrix}_n^{\ell,(k)}
  \begin{bmatrix}
    \delta \Vec{u} \\
    \delta \Vec{t}_S \\
    \delta \Vec{t}_N \\
    \delta \Vec{t}_T \\
    \delta \Vec{t}_O \\
    \delta \Vec{p}
  \end{bmatrix}
  = -
  \begin{bmatrix}
    \Vec{r}_u \\
    \Vec{r}_{S} \\
    \Vec{r}_{N} \\
    \Vec{r}_{T} \\
    \Vec{r}_{O} \\
    \Vec{r}_p
  \end{bmatrix}_n^{\ell,(k)},
  \label{eq:jac_sys}\\
  \text{\underline{Set}} \qquad &
  \begin{bmatrix}
    \Vec{u} \\
    \Vec{t}_S \\
    \Vec{t}_N \\
    \Vec{t}_T \\
    \Vec{t}_O \\
    \Vec{p}
  \end{bmatrix}_n^{\ell,(k+1)}
  =
  \begin{bmatrix}
    \Vec{u} \\
    \Vec{t}_S \\
    \Vec{t}_N \\
    \Vec{t}_T \\
    \Vec{t}_O \\
    \Vec{p}
  \end{bmatrix}_n^{\ell,(k)}
  +
  \begin{bmatrix}
    \delta \Vec{u} \\
    \delta \Vec{t}_S \\
    \delta \Vec{t}_N \\
    \delta \Vec{t}_T \\
    \delta \Vec{t}_O \\
    \delta \Vec{p}
  \end{bmatrix}.
  \label{eq:newton_update}
\end{align}

\noindent
Detailed expressions for the residual vectors and sub-matrices appearing in the linearized system \eqref{eq:jac_sys} are given in \ref{app:FEFV_vec_mat_blocks}.
Note that matrices $A_{TT}$ and $A_{OO}$ are diagonal with entries equal to the area of the face the Lagrange multiplier degrees of freedom are associated with.
Therefore, $\delta \Vec{t}_T$ and $\delta \Vec{t}_O$ can be eliminated through static condensation leading to the reduced system:
\begin{align}
  \begin{bmatrix}
    A_{uu} + B_{uu} & A_{uS} & A_{uN} + B_{uN} & A_{up} \\
    A_{Su} & 0 & 0 & 0 \\
    A_{Nu} & 0 & 0 & 0 \\
    A_{pu} & 0 & 0 & A_{pp} \\
  \end{bmatrix}_n^{\ell,(k)}
  \begin{bmatrix}
    \delta \Vec{u} \\
    \delta \Vec{t}_S \\
    \delta \Vec{t}_N \\
    \delta \Vec{p} \\
  \end{bmatrix} &= -
  \begin{bmatrix}
    \Vec{r}_n^u \\ 
    \Vec{r}_S \\
    \Vec{r}_N \\
    \Vec{r}_p \\
  \end{bmatrix}_n^{\ell,(k)}
  +
  \begin{bmatrix}
    A_{uT} & A_{uO} \\
    0 & 0 \\
    0 & 0 \\
    0 & 0
  \end{bmatrix}_n^{\ell,(k)}
  \begin{bmatrix}
    A_{TT}^{-1} & 0\\
    0           & A_{OO}^{-1}
  \end{bmatrix}_n^k
  \begin{bmatrix}
    \Vec{r}_{T} \\
    \Vec{r}_{O} \\
  \end{bmatrix}_n^{\ell,(k)},
  \label{eq:jac_sys_reduced}
\end{align}
with
$B_{uu} = - A_{uT} A_{TT}^{-1} A_{Tu}$ and
$B_{uN} = - A_{uT} A_{TT}^{-1} A_{TN}$
\begin{rem}
In the formulation proposed in \cite{Fra18}, the Jacobian system is assembled in its reduced form with matrices $B_{uu}$ and $B_{uN}$ directly assembled as
\begin{subequations}
\begin{align}
  [ B_{uu} ]_{ij} =
  & \left( \llbracket \tensorOne{\eta}_{i,T} \rrbracket, \left.\frac{\partial \vec{t}_T^* }{\partial \Delta_n \tensorTwo{g}_T } \right|_n^{\ell,(k)} \cdot \llbracket \tensorOne{\eta}_{j,T} \rrbracket \right)_{\Gamma_{f,n}^{\ell,\text{slip}}}
  && \forall(i,j) \in \mathcal{N}_u \times \mathcal{N}_u,
  \label{eq:jac_blk_Buu} \\
  [ B_{uN} ]_{ij} =
  & \left( \llbracket \tensorOne{\eta}_{i,T} \rrbracket, \left.\frac{\partial \vec{t}_T^* }{\partial t_N } \right|_n^{\ell,(k)} \mu_{j,N} \right)_{ \Gamma_{f,n}^{\ell,\text{slip}} }
  && \forall(i,j) \in \mathcal{N}_u \times \mathcal{N}_t^N.
  \label{eq:jac_blk_BuN}
\end{align}
\label{eq:jac_blk_reduced}\null
\end{subequations}
Expressions for the partial derivatives in \eqref{eq:jac_blk_reduced} are provided in \ref{app:FEFV_vec_mat_blocks}.
\end{rem}
\begin{rem}
The focus of this work is on the nonlinear algorithm and stabilization only. Thus, the linear solution step
in \eqref{eq:jac_sys} is carried out using a direct solver. The design of a scalable solver for this linear system is clearly nontrivial, however, and is the subject of future research.
\end{rem}

\section{Stabilization}
\label{sec:stabilization}

In this section we explore a family of stabilization techniques to enable the successful use of the spatial discretization presented in Section \ref{sec:variational_form_new}, 
We first elaborate the proposed stabilization procedures considering the simpler stick-contact problem, in the absence of frictional slip or fluid flow.
The extension to more sophisticated physics will then follow seamlessly at the end of the section.

To clearly highlight the source of instability, consider a fracture entirely in stick mode, $\Gamma_f = \Gamma_f^\text{stick}$.  In this case, system \eqref{eq:jac_sys_reduced} reduces to
\begin{equation}
  \begin{bmatrix}
    K   & C \\
    C^T & 0 \\
  \end{bmatrix}
  \begin{bmatrix}
    \delta \Vec{u} \\
    \delta \Vec{t}_S \\
  \end{bmatrix} = -
  \begin{bmatrix}
    \Vec{r}_u \\
    \Vec{r}_S \\
  \end{bmatrix},
  \label{eq:jab_strickonly}
\end{equation}
with $K = A_{uu}$ and $C = A_{uS}$.
For this saddle point system, stability and well-posedeness of the discrete problem requires fulfillment of an inf-sup condition \cite{BreBat90}.
In essence, ``the trace space of the discrete displacement has to be well balanced with the finite-dimensional space for the surface traction''  \cite{wohlmuth2011variationally}.
Unfortunately, the $\vec{Q}_1$-displacement/$\vec{P}_0$-Lagrange multiplier interpolation is not inf-sup stable, leading to unstable approximations \cite[Section 3.1]{wohlmuth2011variationally}.
\begin{figure}[t]
  \hfill
  \begin{subfigure}{0.3\textwidth}
    \centerline{\includegraphics[width=\linewidth]{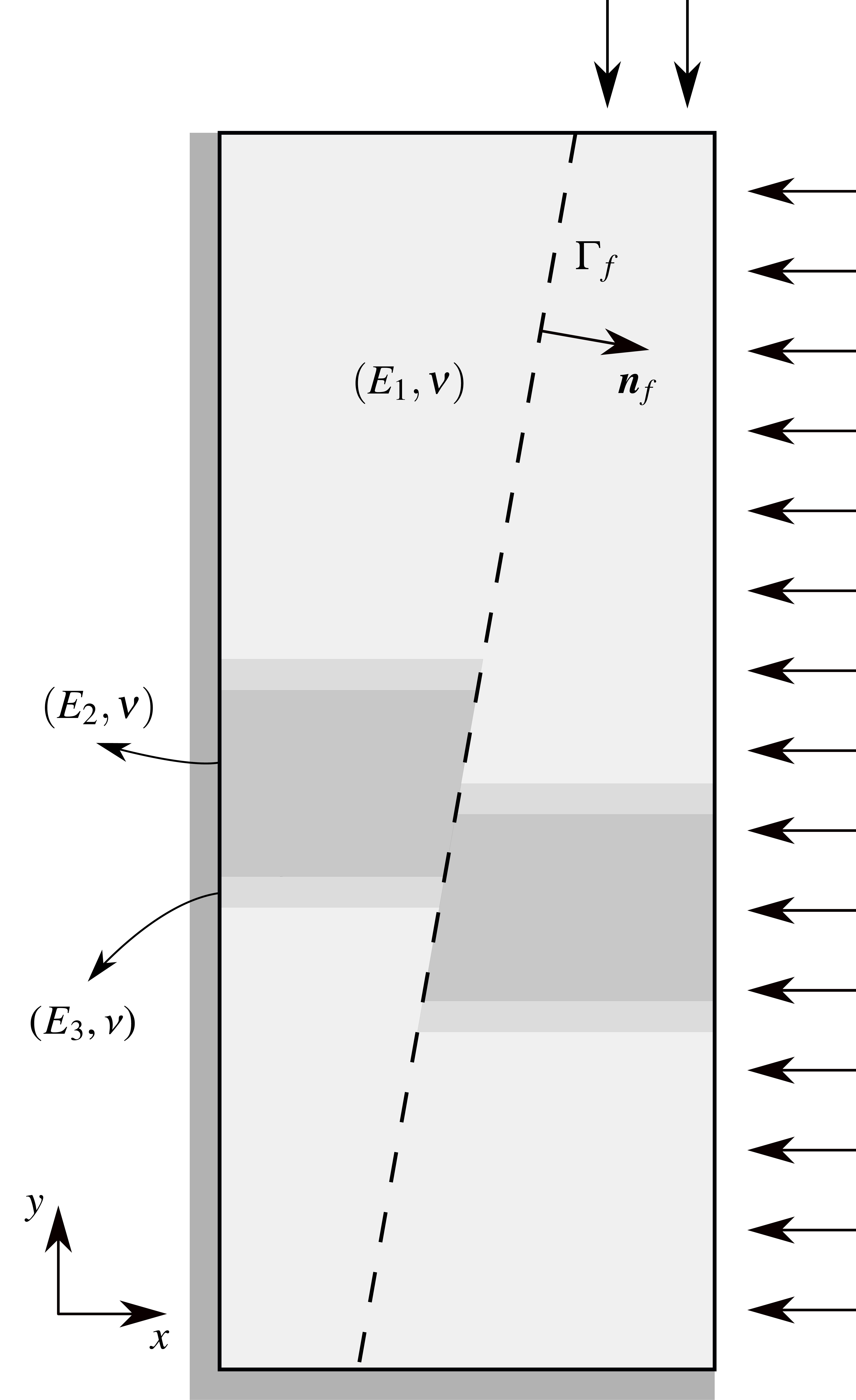}}
    \caption{}
  \end{subfigure}
  \hfill
  \begin{subfigure}{0.3\textwidth}
    \centerline{\includegraphics[width=\linewidth]{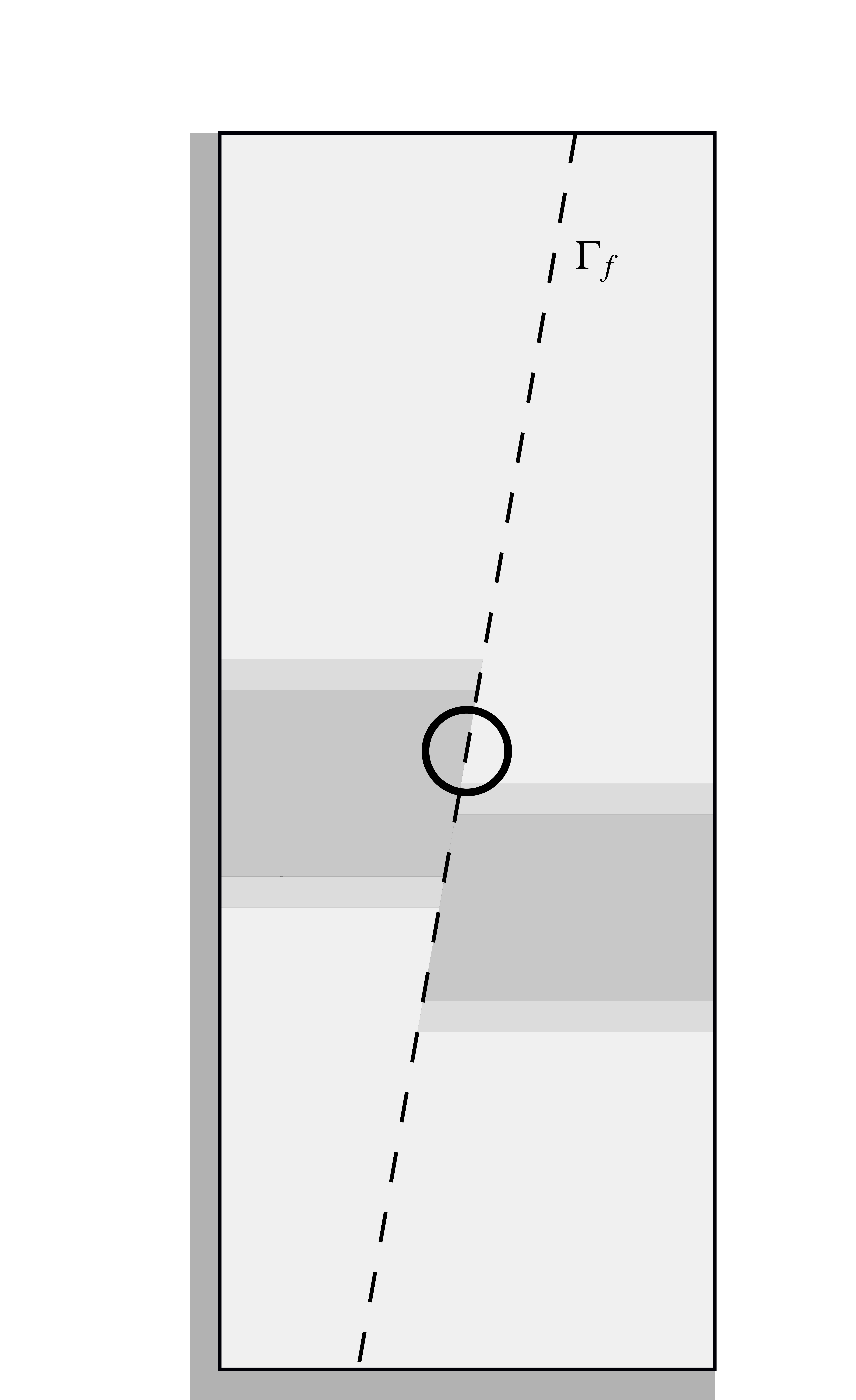}}
    \caption{}
  \end{subfigure}
  \hfill\null
  \caption{Sketch of the 2D test case used to compare stabilizations. The black circle locates the area where a more detailed zoom is taken.}
  \label{fig:stab_2D_sketch}
\end{figure}
\begin{figure}[t]
  \centering
  \includegraphics[width=0.40\linewidth]{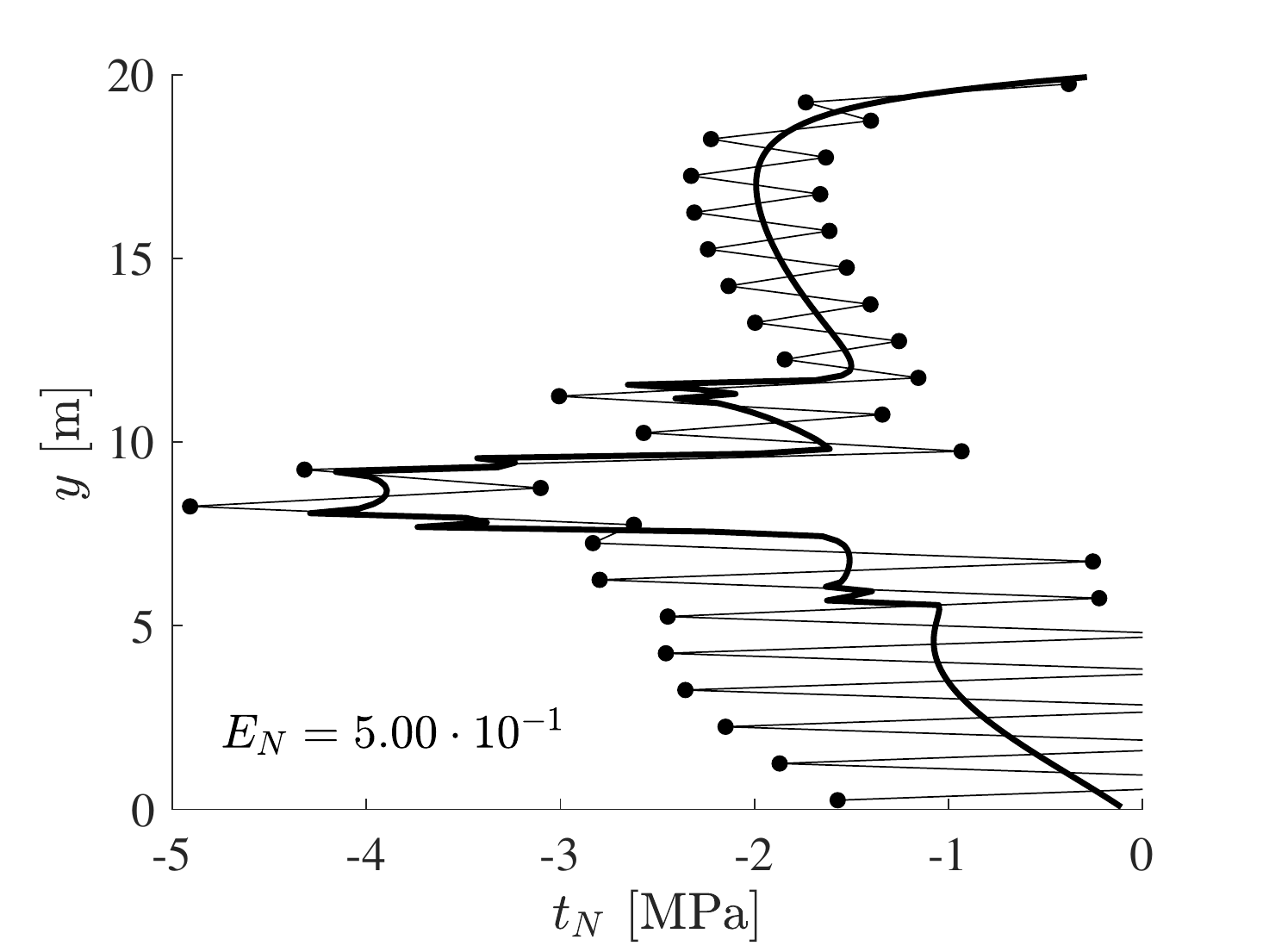}
  \includegraphics[width=0.40\linewidth]{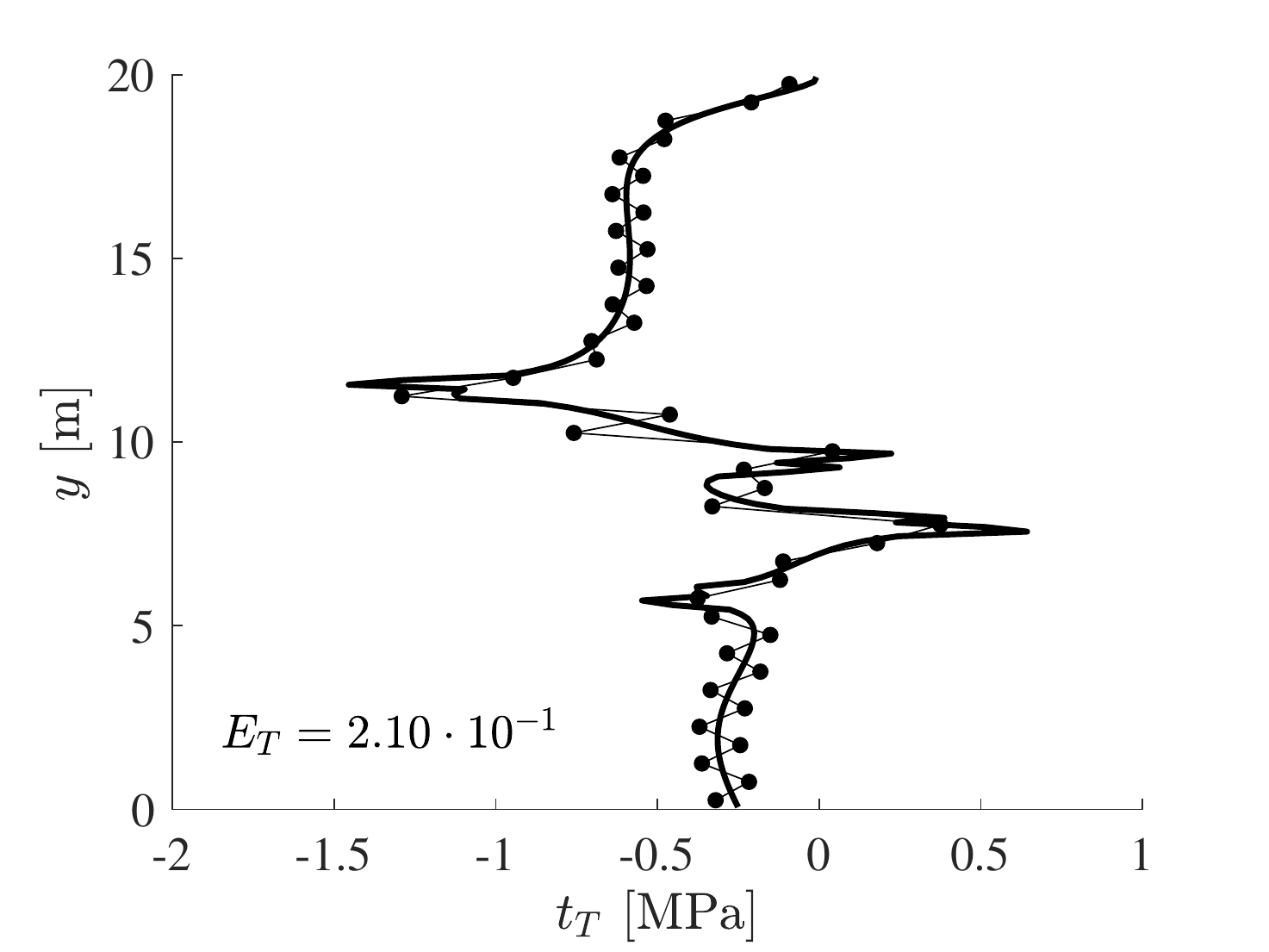}
  \caption{Numerical results for example in Fig. \ref{fig:stab_2D_sketch}, without any stabilization ($16 \times 40$ elements). In these plots, as
    well as in the following ones, the thick line represents the \textit{continuous}
    solution obtained without the fault, while the points are the
    solution using the $\vec{Q}_1$--$\vec{P}_0$ approach.}
  \label{fig:stab_2D_r0}
\end{figure}

As we describe strategies to fix this deficiency, it is useful to apply them to an illustrative example for comparison purposes. We will use the 2D, plain strain problem shown in Fig. \ref{fig:stab_2D_sketch}.
The size of the domain is $8 \times 20 \text{ m}$ and the fracture has a dip of $10^\circ$.
Three material regions are considered, characterized by Young's modulus values $E_1 = 3 \text{ GPa}$, $E_2 = 15 \text{ GPa}$, and $E_3 = (E_1 + E_2)/2 \text{ GPa}$, and a homogeneous Poisson's ratio $\nu = 0.25$.
The horizontal and vertical loads are $2 \text{ MPa} \cdot \text{m}$ and $4 \text{ MPa} \cdot
\text{m}$, respectively.
We compute solutions on a base grid with $16 \times 40$ elements, as well as on an anisotropically refined one with $16 \times 200$ elements.
The resulting traction components orthogonal and parallel to the fracture are identified as $t_N$ and $t_T$,
respectively.
We also compute a reference solution using a highly refined elastic model without a discontinuity (since pure stick conditions are assumed).
In the remainder of the section, we will refer to this numerical solution as the \textit{continuous} one.
For each stabilized solution, we compute two integral relative differences with respect to the continuous one, $E_N$ and $E_T$, for the two components of the traction. 
In Fig. \ref{fig:stab_2D_r0}, results using the $\vec{Q}_1$--$\vec{P}_0$ interpolation without any stabilization are shown. We observe that the traction solution \textit{on average} coincides with the
continuous one, but it exhibits substantial checkerboard oscillations. 

We now consider three possible strategies to fix this issue:
\begin{enumerate}
\item Analytic macroelement stabilization
\item Algebraic macroelement stabilization
\item Algebraic global stabilization
\end{enumerate}
As we will see, the first method makes significant assumptions regarding the grid topology and material heterogeneity, while the latter methods provide greater generality.

\subsection{Analytic macroelement stabilization}
\label{sec:macroelem}

\begin{figure}
  \centering
  \hfill
  \begin{subfigure}{0.4\textwidth}
    \centerline{\includegraphics[width=\linewidth]{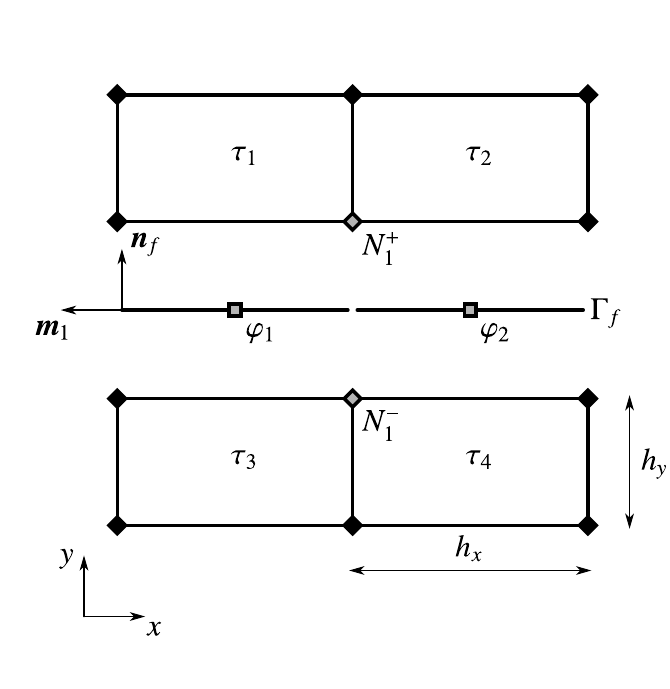}}
    \caption{}
    \label{fig:active_dof_2D}
  \end{subfigure}
  \hfill
  \begin{subfigure}{0.4\textwidth}
    \centerline{\includegraphics[width=\linewidth]{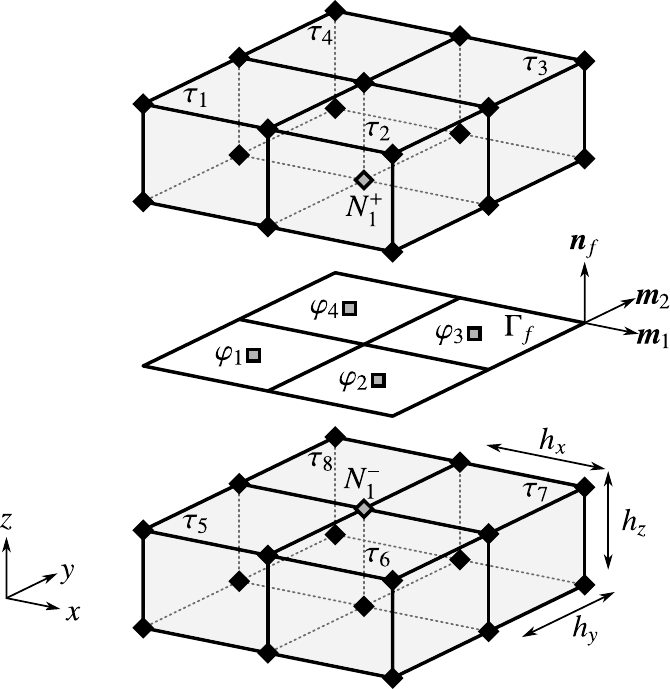}}
    \caption{}
    \label{fig:active_dof_3D}
  \end{subfigure}
  \hfill\null
  \caption{Location of displacement and traction degrees of freedom (DOFs) for a 2D (a) and a 3D (b) macroelement, respectively. Diamonds (\protect\tikz \protect\draw[black,fill=gray!29,line width=.25ex,rotate=45] (0,0) rectangle (0.71ex,0.71ex);) show the location of displacement DOFs, while squares (\protect\tikz \protect\draw[black,fill=gray!29,line width=.25ex] (0,0) rectangle (1ex,1ex);) represent the position of traction DOFs. Black (\protect\tikz \protect\draw[black,fill=black,line width=.25ex,rotate=45] (0,0) rectangle (0.71ex,0.71ex);) and gray (\protect\tikz \protect\draw[black,fill=gray!29,line width=.25ex,rotate=45] (0,0) rectangle (0.71ex,0.71ex);) means fixed and active DOFs, respectively.}
  \label{fig:active_dof}
\end{figure}

We first explore stabilizing the discretization using a macroelement approach \cite{elman2014finite}. 
The method is most easily described using a 2D reference macroelement as shown in Fig. \ref{fig:active_dof_2D}.
This macroelement is formed by two interface elements and four quadrilaterals, two for each side of the fracture.
To create a patch test, the ten displacement degrees of freedom (DOFs) located on the boundary are fixed.
The objective is to derive a stabilization that produces a well-posed problem on an individual patch.
Such a stabilization then implies well-posedness on a grid consisting of stabilized macroelements.

We can achieve this by investigating the kernel modes of the Schur complement for the Lagrange multipliers. For example, for the macroelement of Fig. \ref{fig:active_dof_2D}, there are $4$ displacement DOFs---the $x$- and $y$-components for innermost nodes, above and below the fracture---and $4$ traction DOFs---one normal and one tangential component for each element on $\Gamma_f$.
Assuming the ordering $\Vec{u} = \left[ u_x^{(N_1^+)}, u_y^{(N_1^+)}, u_x^{(N_1^-)}, u_y^{(N_1^-)} \right]$, $\Vec{t} = \left[ t_N^{(\varphi_1)}, t_T^{(\varphi_1)}, t_N^{(\varphi_2)}, t_T^{(\varphi_2)} \right]$ for the unknowns, the explicit expressions of $K$ and $C$ are:
\begin{align}
  K &=
   \frac{E(1-\nu)}{3(1+\nu)(1-2\nu)} \begin{bmatrix} K_2 & 0 \\ 0 & K_2  \end{bmatrix},
  &
  K_{2} &= 
  \begin{bmatrix}
    \frac{\left(1-2\nu\right)}{(1-\nu)}\frac{h_x}{h_y} + \frac{2 h_y}{h_x} & 0 \\
    0 & \frac{(1-2\nu)}{(1-\nu)}\frac{h_y}{h_x} + \frac{2 h_x}{h_y} \\
  \end{bmatrix},
  & 
  C &=
  \frac{h_x}{2} \begin{bmatrix} R_2 & R_2 \\ -R_2 & -R_2  \end{bmatrix},
  \label{eq:K_C_macroelem2D}
\end{align}
with $E$ the Young's modulus, $\nu$ the Poisson's ratio, $h_x$ and $h_y$ the mesh size in the $x$- and $y$-direction, respectively, and $R_2$ a
$2 \times 2$ rotation matrix from the local to global reference system. In the example of Fig.
\ref{fig:active_dof_2D}, $R_2$ reads
\begin{equation}
  R_2 = \begin{bmatrix}
    \vec{n}_f \cdot \vec{e}_x &
    \vec{m}_1 \cdot \vec{e}_x \\
    \vec{n}_f \cdot \vec{e}_y &
    \vec{m}_1 \cdot \vec{e}_y
  \end{bmatrix} = \begin{bmatrix}
    0 & -1 \\
    1 & 0 \\
  \end{bmatrix},
  \label{eq:rot_mat2}
\end{equation}
with $\{ \vec{e}_x, \vec{e}_y \}$ the standard Euclidean basis in $\mathbb{R}^2$.
By definition, the Schur complement $S \in \mathbb{R}^{4\times4}$ is
\begin{equation}
  S = -C^T K^{-1} C =
   - \left(\frac{h_x}{2}\right)^2 \frac{6(1+\nu)(1-2\nu)}{E(1-\nu)}
  \begin{bmatrix}
    K_2^{-1} & K_2^{-1} \\
    K_2^{-1} & K_2^{-1} \\
  \end{bmatrix} 
  .
  \label{eq:schur_2D}
\end{equation}
Performing an eigen-decomposition, its complete set of eigenpairs $\{ \lambda_i, \Vec{v}_i \}_{i=1}^4$ is
\begin{subequations}
\begin{align}
  \{\lambda_1, \Vec{v}_1 \} & = \left\{ 0, \begin{bmatrix*}[r] 1 \\ 0 \\ -1 \\  0 \end{bmatrix*} \right\}, &
  \{\lambda_2, \Vec{v}_2 \} & = \left\{ 0, \begin{bmatrix*}[r] 0 \\ 1 \\  0 \\ -1 \end{bmatrix*} \right\},
  \label{eq:schur_eigenpair_kernel} \\
  \{\lambda_3, \Vec{v}_3 \} & = \left\{ - \frac{(1+\nu)(1-2\nu)}{E(1-\nu)} \frac{ 3 h_x^2 }
{\frac{h_x}{h_y}\frac{\left(1-2\nu\right)}{(1-\nu)} + \frac{2 h_y}{h_x} }, \begin{bmatrix*}[r] 1 \\ 0 \\  1 \\  0 \end{bmatrix*} \right\}, & 
  \{\lambda_4, \Vec{v}_4 \} & = \left\{ - \frac{(1+\nu)(1-2\nu)}{E(1-\nu)} \frac{ 3 h_x^2 }
{\frac{h_y}{h_x}\frac{\left(1-2\nu\right)}{(1-\nu)} + \frac{2 h_x}{h_y} }, \begin{bmatrix*}[r] 0 \\ 1 \\  0 \\  1 \end{bmatrix*} \right\}.
  \label{eq:schur_eigenpair_nonKernel}
\end{align}
\end{subequations}
We observe that $S$ has rank $2$.
Eigenvectors $\Vec{v}_1$ and $\Vec{v}_2$ span the kernel of $S$ and represent the \textit{checkerboard} mode for the traction normal- and tangential-component, respectively, confirming the numerical results shown in Fig. \ref{fig:stab_2D_r0}.
Conversely, the component-wise constant eigenvectors $\Vec{v}_3$ and $\Vec{v}_4$ span the column space of $S$.

The null eigenvectors $\Vec{v}_1$ and $\Vec{v}_2$ are the source of the macroelement instability and need to be removed from the null space of $S$.
To do so, we introduce the symmetric and positive semi-definite matrix $H^* = V V^T$, with $V$ the matrix having columns $\Vec{v}_1$ and $\Vec{v}_2$.
By construction, $\Vec{v}_1$ and $\Vec{v}_2$ are eigenvectors of $H^*$, both corresponding to the eigenvalue 2, which span the range of $H^*$.
Also, $\Vec{v}_3$ and $\Vec{v}_4$  are now a basis for the kernel of $H^*$.
Let $H = \alpha H^*$ be a scaled matrix, with $\alpha$ a scalar stabilization constant having units of squared length per pressure, i.e. the same units as the eigenvalues of $S$. The system to be solved is modified as
\begin{equation}
  \begin{bmatrix}
    K   & C \\
    C^T & - H \\
  \end{bmatrix}
  \begin{bmatrix}
    \delta \Vec{u} \\
    \delta \Vec{t}_S \\
  \end{bmatrix} = -
  \begin{bmatrix}
    \Vec{r}_u \\
    \Vec{r}_S - H_S \Vec{t}_S \\
  \end{bmatrix},
  \label{eq:jab_strickonly_stab}
\end{equation}
where a stabilizing contribution now replaces the zero block of the original matrix.  The resulting modified Schur complement is then $\tS = -C^T K^{-1} C - \alpha H^*$.
The eigenpairs of $\tS$ are the same as those of $S$ with the only difference that $\Vec{v}_1$ and $\Vec{v}_2$ are now associated with the same non zero eigenvalue $2\alpha$.
The scaling constant $\alpha$ is chosen in such a way that, in a homogeneous case, the eigenvalues of $\tS$ are bounded between $\min(\lambda_3,\lambda_4)$ and $\max(\lambda_3,\lambda_4)$, so that the matrix conditioning does not depend on the stabilization constant.
For a regular Cartesian grid with $h_x = h_y = h$, from Eq. \eqref{eq:schur_eigenpair_nonKernel} we can write
\begin{equation}
  \lambda = \lambda_3 = \lambda_4 = -\frac{\left(1+\nu\right)\left(1-2\nu\right)}{E}
\frac{3 h^2}{\left(3-4\nu\right)},
\end{equation}
thus, the ``optimal'' scaling factor simply reads:
\begin{equation}
  \alpha^* = \frac{\left(1+\nu\right)\left(1-2\nu\right)}{E} \frac{3}{2}\frac{h^2}
{\left(3-4\nu\right)}.
  \label{eq:stab_2D_alpha}
\end{equation}
In case of geometric anisotropy, $h$ can be computed as the average length of the interface elements composing the macroelement.

We emphasize that the stabilization matrix $H$ represents an example of a minimal stabilization operator \cite{ElmSilWat02} that does not pollute the physical eigenpairs.  It requires, however, explicit knowledge of the eigenvectors associated with non-zero eigenvalues. 
For the macroelement of Fig. \ref{fig:active_dof_2D}, $H$ can also be interpreted as a macroelement stabilization matrix for the $\alpha$-weighted inter-element traction (component-wise) jump since the following relationship holds true
\begin{equation}
  \Vec{t}^T H \Vec{t} =
  \alpha \Vec{t}^T H^* \Vec{t} =
  \alpha \left[ \left( t_N^{(\varphi_1)} - t_N^{(\varphi_2)} \right)^2 + \left( t_T^{(\varphi_1)} - t_T^{(\varphi_2)} \right)^2  \right].
\end{equation}

\begin{figure}
  \centering
  \includegraphics[width=0.40\linewidth]{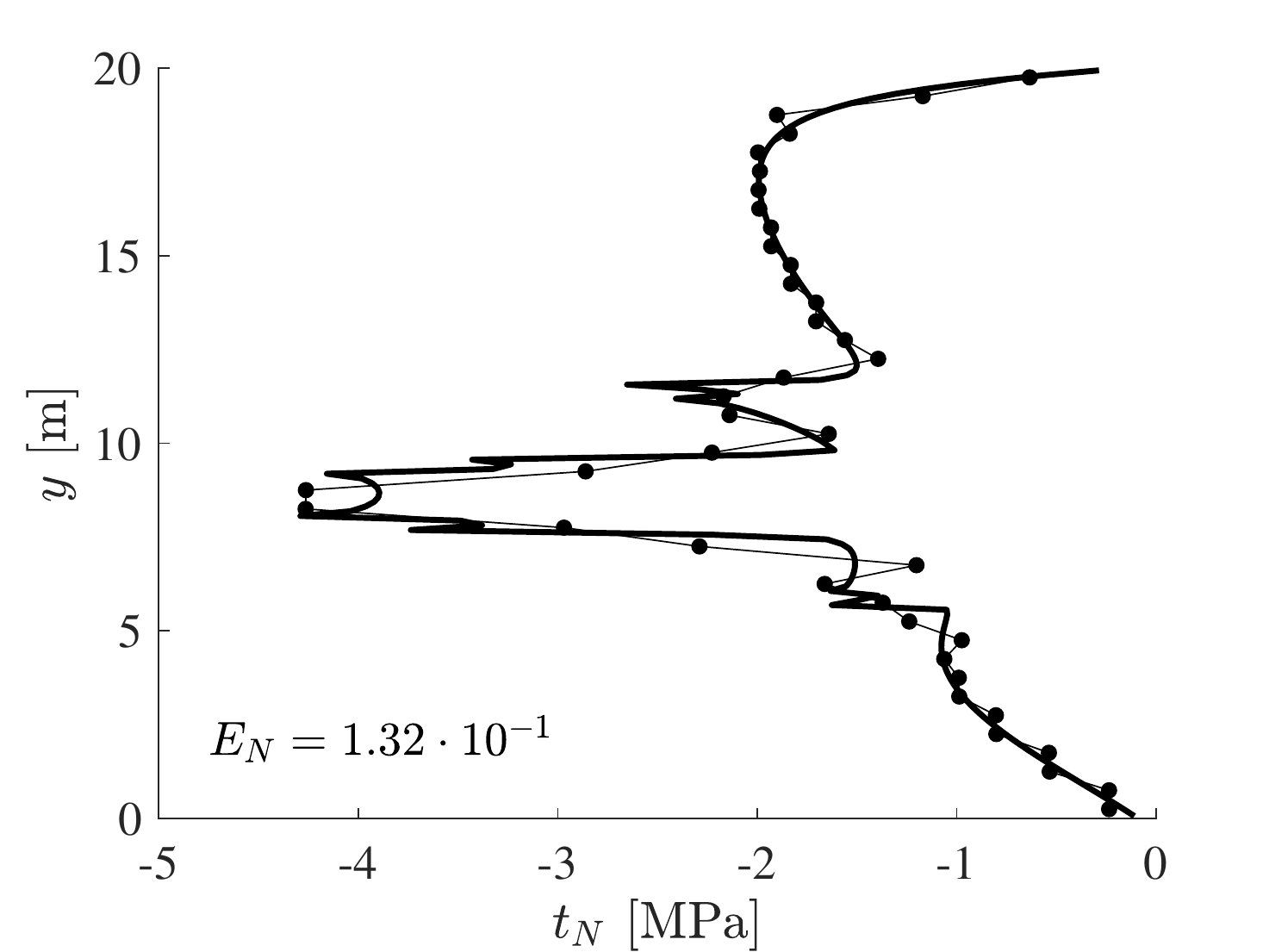}
  \includegraphics[width=0.40\linewidth]{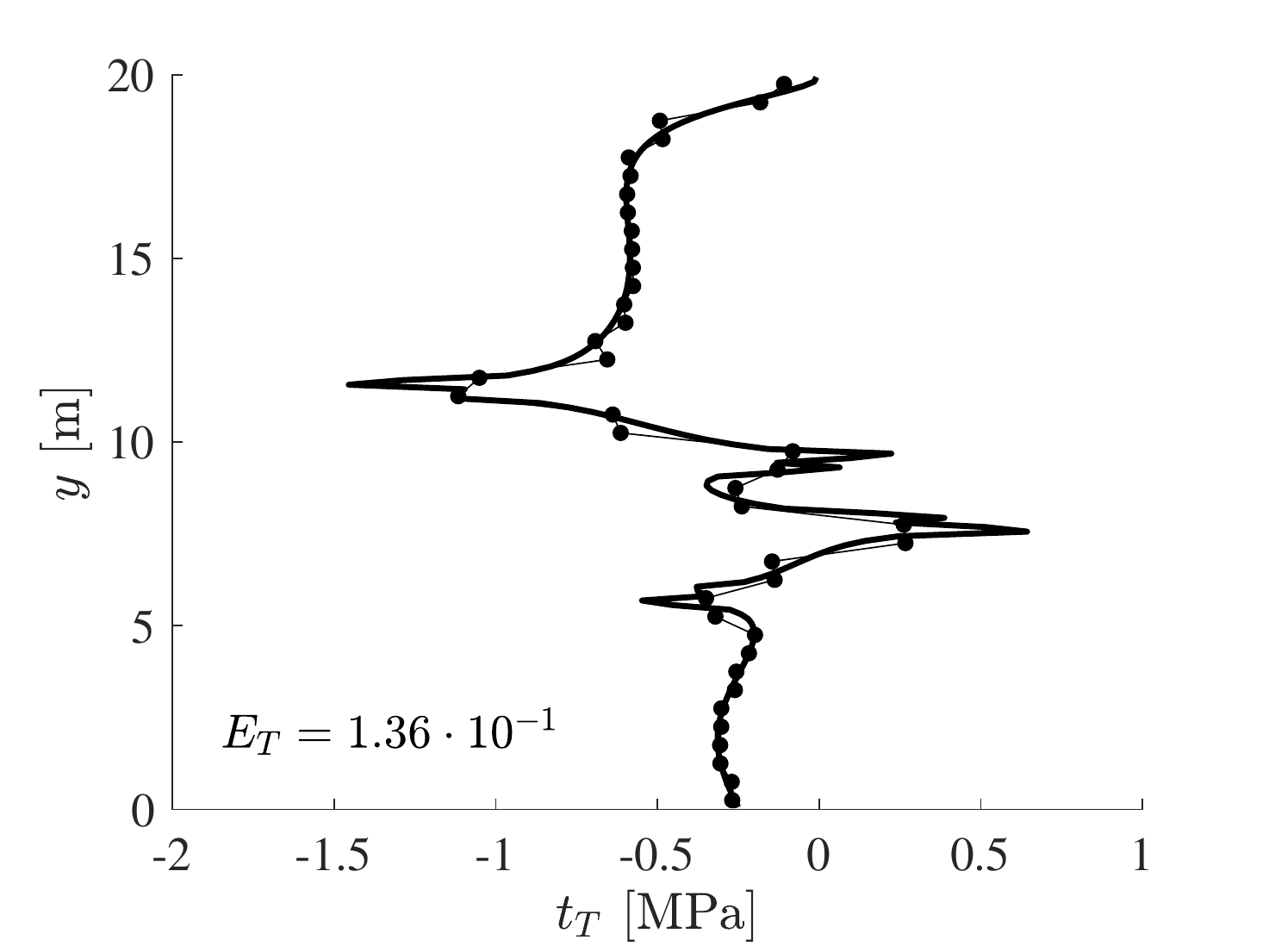}
  \includegraphics[width=0.40\linewidth]{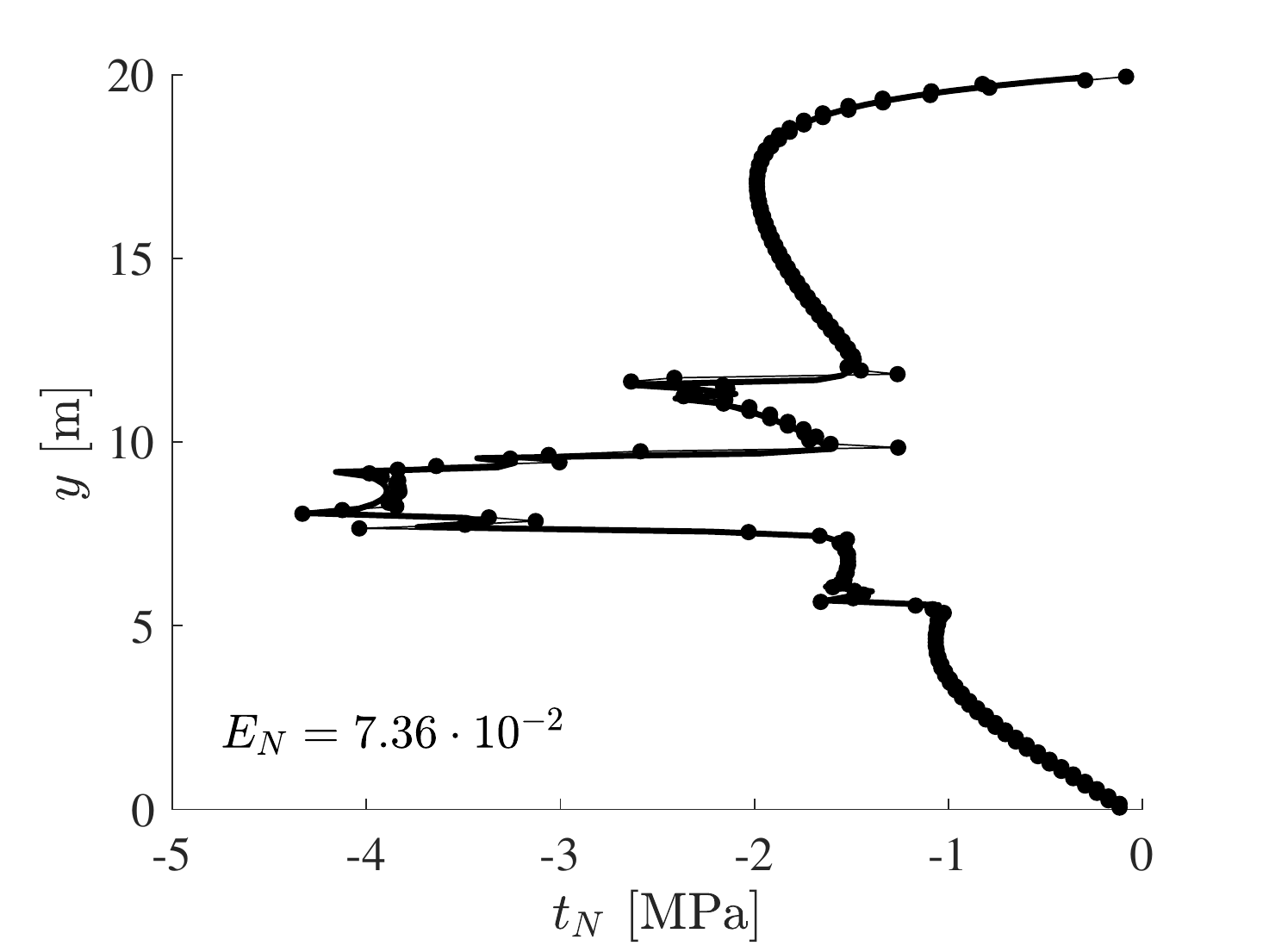}
  \includegraphics[width=0.40\linewidth]{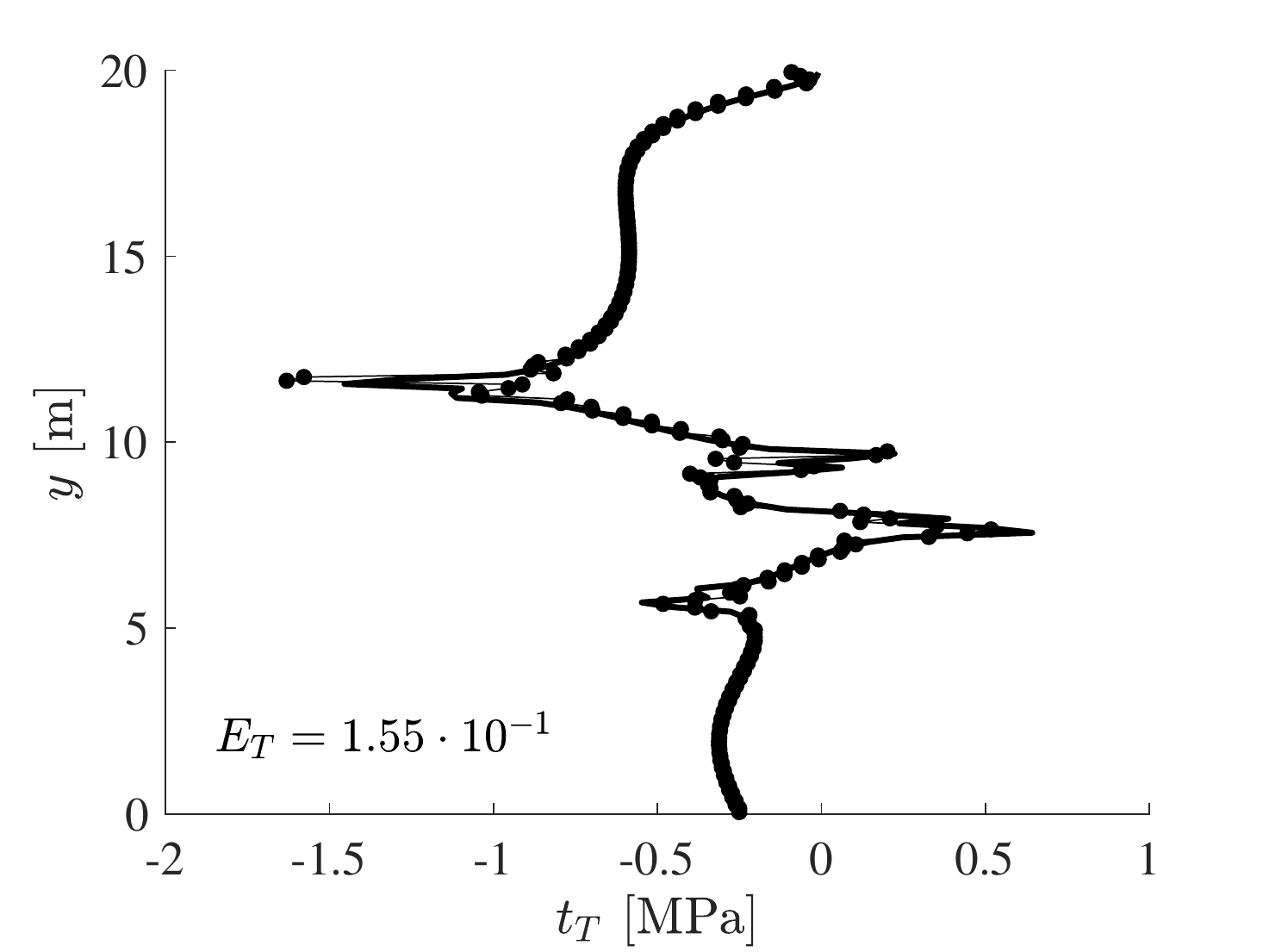}
  \caption{Numerical results for example in Fig. \ref{fig:stab_2D_sketch} using the analytic macroelement stabilization. First row: base grid, last row:
    vertically refined grid.}
  \label{fig:stab_2D_r1}
\end{figure}

Applying this stabilization to the 2D test case (Fig. \ref{fig:stab_2D_sketch}) yields a substantial reduction of the oscillations observed in the unstabilized case as shown in Fig. \ref{fig:stab_2D_r1}, with the errors $E_N = 1.32 \cdot 10^{-1}$, $E_T = 1.36 \cdot 10^{-1}$ and $E_N = 7.36 \cdot 10^{-2}$, $E_T = 1.55 \cdot 10^{-1}$ for the base and the refined grids, respectively.

\begin{rem}
The 2D reference macroelement (Fig. \ref{fig:active_dof_2D}) consists of four equal rectangular elements---hence, the stiffness matrix $K$, see Eq. \eqref{eq:K_C_macroelem2D}, is diagonal---with the normal $\vec{n}_f$ aligned with the $y$-axis.
Consequently, normal- and tangential-component of the traction are decoupled as the eigenvectors defined in Eqs. \eqref{eq:schur_eigenpair_kernel}-\eqref{eq:schur_eigenpair_nonKernel} reveal.
This is not the case in a more general geometry where the two traction components are typically related.
\end{rem}

\begin{figure}
  \centering
  \includegraphics[width=0.20\linewidth]{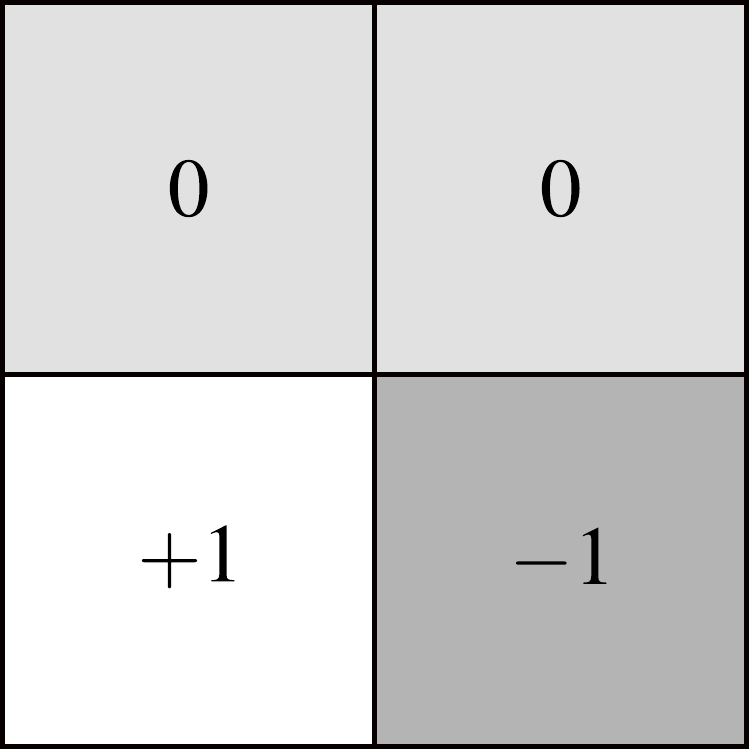}\qquad
  \includegraphics[width=0.20\linewidth]{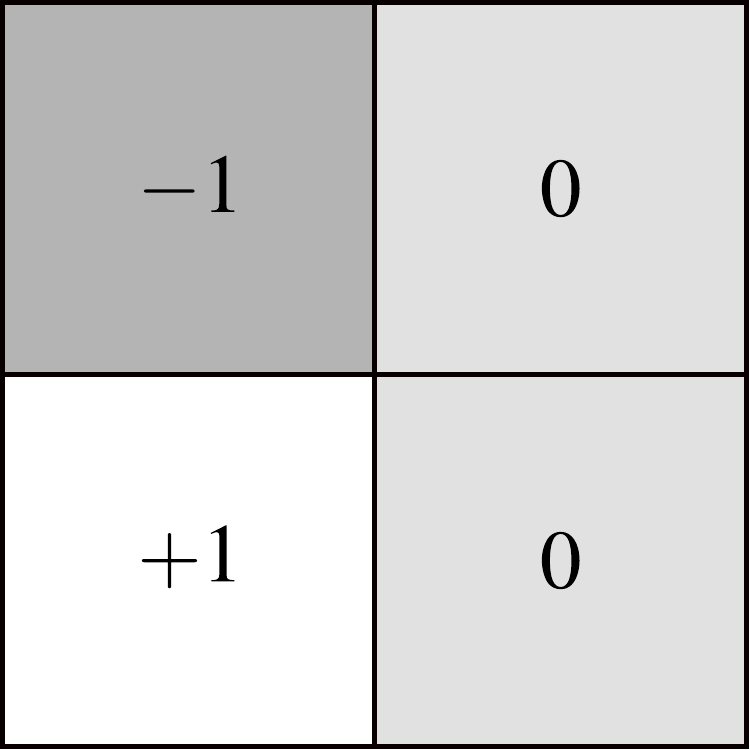}\qquad
  \includegraphics[width=0.20\linewidth]{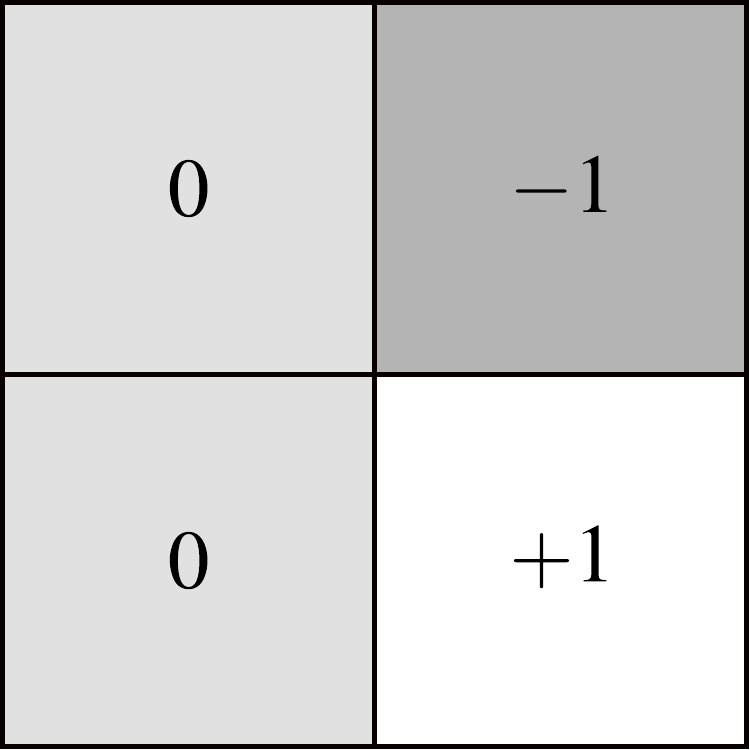}
  \caption{Three different modes that need to the stabilized for the 3D case. Each
    of them applies to the three components (one normal and two tangential) of the traction vector,
    thus the kernel size is 9. Note that the remaining mode is a linear combination
    of these three, thus it is part of the same space, the Schur complement kernel.}
  \label{fig:stab3D_modes}
\end{figure}

Extension to 3D is based on the reference macroelement shown in Fig. \ref{fig:active_dof_3D}, which consists of four quadrilateral interface elements and eight hexahedral elements.
The derivation follows the same steps discussed above for the 2D case.
We omit the computations and simply report the main results.
The Schur complement $S \in \mathbb{R}^{12\times12}$ has rank three with column space spanned by three component-wise constant eigenvectors.
There are nine spurious traction modes that need to be stabilized.
A convenient basis for the kernel of $S$ is shown in Fig. \ref{fig:stab3D_modes} and corresponds to checkerboard-like modes for each component of the traction vector with respect to three internal edge of the interface element patch.
Note that the choice of the three edges is arbitrary.
The optimal $\alpha$ value to obtain a stabilization contribution that falls in the already present lower/upper spectral Schur complement bounds is:

\begin{equation}
  \alpha^* = \frac{\left(1+\nu\right)\left(1-2\nu\right)}{E} \frac{3 h^3}
{\left(16-24\nu\right)}.
\end{equation}

\noindent
As for the 2D case, this value is exact in case of $h_x = h_x = h_z = h$ and constant physical properties.
In general, $h$ can be computed as the cubic root of the average volume of the elements surrounding the fracture and sharing a face.

\subsection{Algebraic macroelement stabilization}
\label{sec:macroelem_noA}

The scaling factor $\alpha$ in the macroelement approach above is a scalar value collecting both mechanical and geometric information.
Its definition can be quite difficult, especially for 3D problems with distorted grids and material heterogeneity.
In this section, we propose an algebraic alternative for computing the stabilization matrix at the macroelement level that circumvents the need for introducing the factor $\alpha$.
The development is based on the following observation: $\alpha$ is fundamentally needed to scale matrix $H^*$ introduced in Sec. \ref{sec:macroelem} so that the spectrum of the stabilized Schur complement $\tS$ is bounded between the smallest and largest nonzero eigenvalues of $S$.
Thus, whenever we provide a different but admissible scaling, $\alpha$ can be avoided.

\begin{figure}
  \hfill
  \begin{subfigure}{0.45\textwidth}
    \centerline{\includegraphics[width=\linewidth]{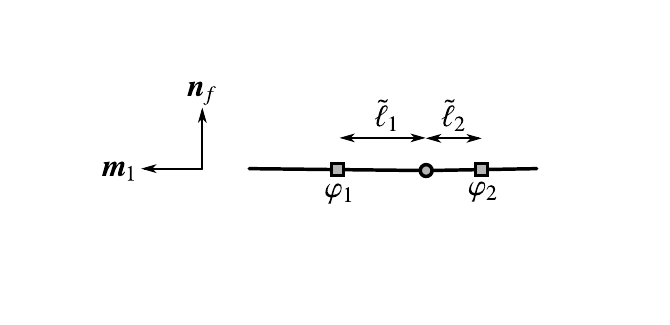}}
    \caption{}
    \label{fig:unstructured_patch_2D}
  \end{subfigure}
  \hfill
  \begin{subfigure}{0.45\textwidth}
    \centerline{\includegraphics[width=\linewidth]{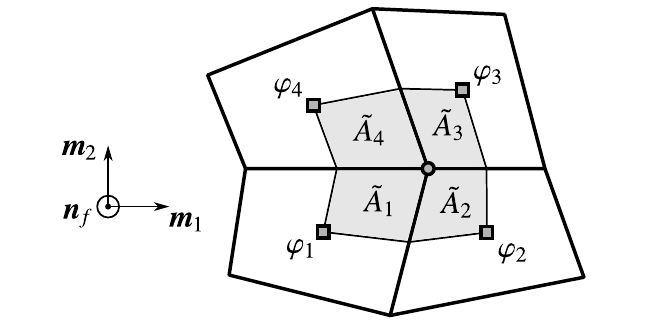}}
    \caption{}
    \label{fig:unstructured_patch_3D}
  \end{subfigure}
  \hfill\null
  \caption{General interface element patches for a 2D (a) and a 3D (b) macroelement, respectively. Squares (\protect\tikz \protect\draw[black,fill=gray!29,line width=.25ex] (0,0) rectangle (1ex,1ex);) indicate the barycenter of the interface elements. Quantities $\tilde{\ell}_i$, $i \in \{1,2 \}$ and $\tilde{A}_j$, $j \in \{1,\ldots,4 \}$ denote the length and area, respectively, associated to the common vertex (\protect\tikz \protect\draw[black,fill=gray!29,line width=.25ex] (0,0) circle (0.5ex);).}
  \label{fig:unstructured_patch}
\end{figure}

Let us consider a general patch of interface elements in a non regular macroelement (Fig. \ref{fig:unstructured_patch}).  The only topological assumption we make is that the fracture elements within an individual macroelement are co-linear in 2D or co-planar in 3D.
For the two dimensional case of Fig. \ref{fig:unstructured_patch_2D}, matrix $C$ given in Eq. \eqref{eq:K_C_macroelem2D} reads
\begin{equation}
  C =
  \begin{bmatrix*}[r]
     \tilde{\ell}_1 R_2 &  \tilde{\ell}_2 R_2 \\
    -\tilde{\ell}_1 R_2 & -\tilde{\ell}_2 R_2 \\
  \end{bmatrix*},
  \label{eq:C_macro_2d_unstructured}
\end{equation}
where $\tilde{\ell}_1$ and $\tilde{\ell}_2$ are the interface element lengths associated to the vertex in common between elements $\varphi_1$ and $\varphi_2$, i.e. the integral of the standard hat function associated to that vertex over $\varphi_1$ and $\varphi_2$, respectively.
The rows of the following matrix $\tC$, formed from $C$ swapping block columns and changing sign to the first block column,
\begin{equation}
  \tC = \begin{bmatrix*}[r]
    -\tilde{\ell}_2 R_2 &  \tilde{\ell}_1 R_2 \\
     \tilde{\ell}_2 R_2 & -\tilde{\ell}_1 R_2 \\
  \end{bmatrix*},
\end{equation}
are orthogonal by construction to $C$, i.e., $C  \tC^T = 0$.
Thus, columns of $\tC^T$ belong to the kernel 
of the Schur complement $S = -C^T K^{-1} C$.
Also, $\tC$ has rank 2, the number of modes that are known to require the stabilization.
Hence, $\tC^T \tC$ is a stabilizing contribution to the Schur complement.
It has to be scaled, but from the observation that $\tC$ has the same entries of $C$, except for the order (they are swapped on a Lagrange multiplier base) and the sign, it is natural to use the inverse of the stiffness matrix in \eqref{eq:K_C_macroelem2D}, i.e. $K^{-1}$, to scale it.
Indeed, a scaling that incorporates element size and material properties information is needed and $K$ collects both of them.
Numerical tests show that the inverse of its diagonal, denoted as $D_{\widetilde{K}} \in \mathbb{R}^{4 \times 4}$, is enough.
Therefore, a macroelement stabilization matrix can be define as $H = \tC^T D_{\widetilde{K}}^{-1} \tC$.
Using the same notation of Sec. \ref{sec:macroelem}, such a stabilization can also be expressed in the following compact form:
\begin{align}
  H &= V D V^T,
  &
  V &=
  \begin{bmatrix*}[r]
    -\tilde{\ell}_2 R_2^T \\ 
     \tilde{\ell}_1 R_2^T
  \end{bmatrix*},
  &
  D &=
  D_{\widetilde{K},N_1}^{-1} + D_{\widetilde{K},N_2}^{-1},
  \label{eq:Hs_macro_2d_unstructured}
\end{align}
with $D \in \mathbb{R}^{2 \times 2}$ a diagonal matrix given by the sum of the inverse of the diagonal portion of the stiffness matrix associated with the two nodes (on both the discontinuity surfaces) shared by the two interface elements forming the macroelement (i.e., $N_1^+$ and $N_1^-$ in Fig. \ref{fig:active_dof_2D}). 
Note that whenever applied to a regular grid, with $h_x = h_y$, and homogeneous material properties, this stabilization reduces to the $\alpha$-based method of the previous section.

\begin{figure}
  \centering
  \includegraphics[width=0.30\linewidth]{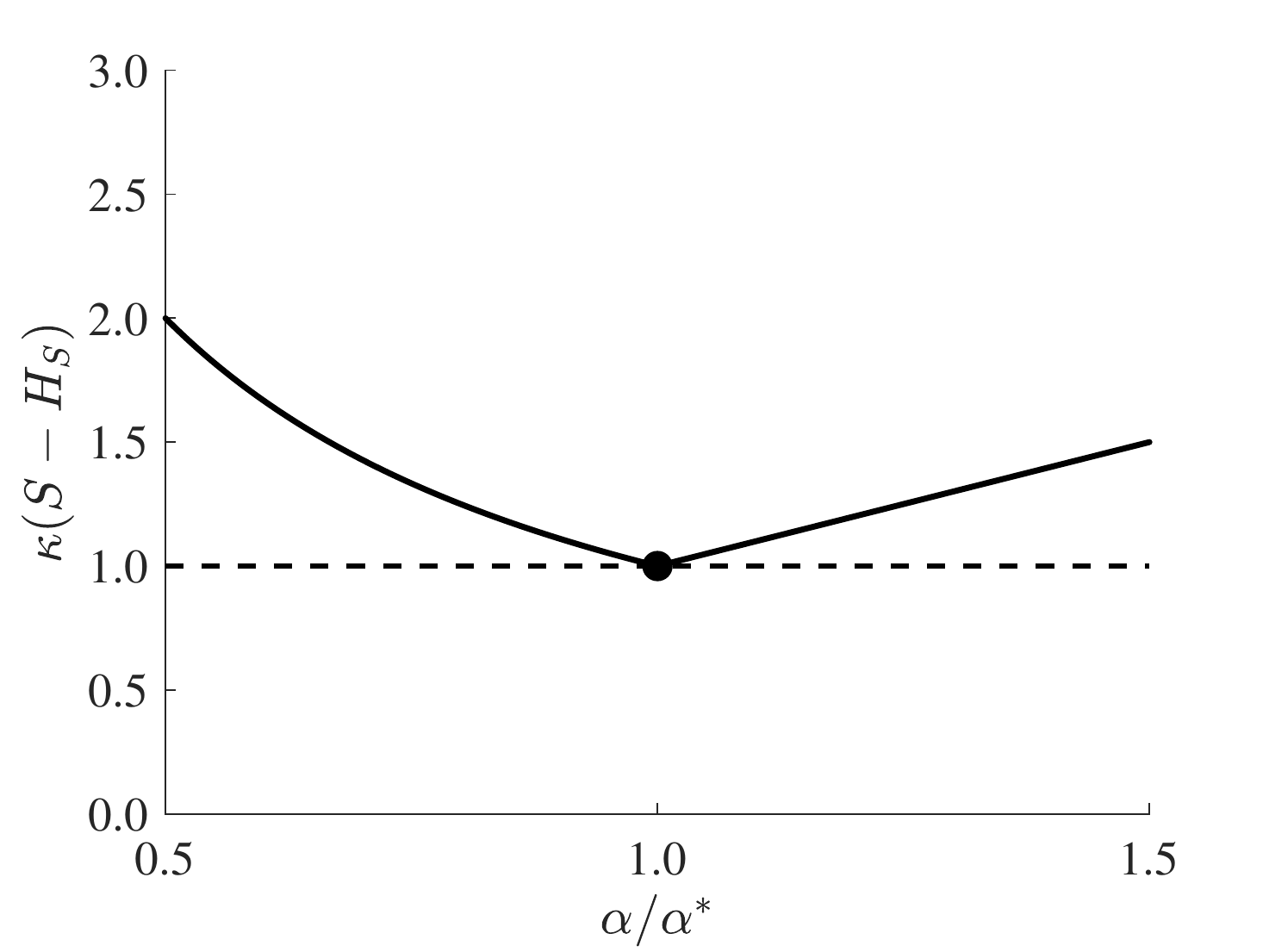}
  \includegraphics[width=0.30\linewidth]{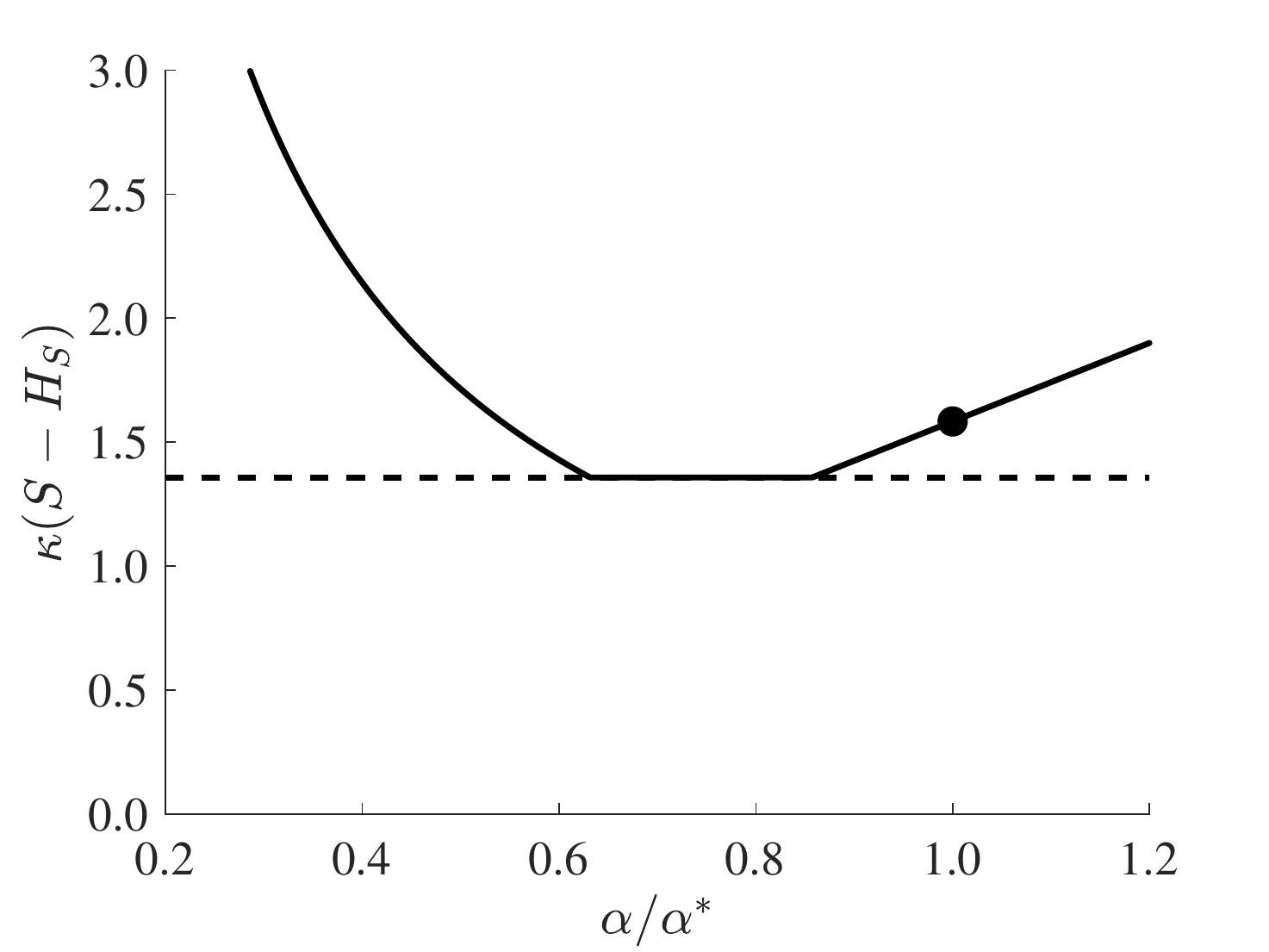}
  \includegraphics[width=0.30\linewidth]{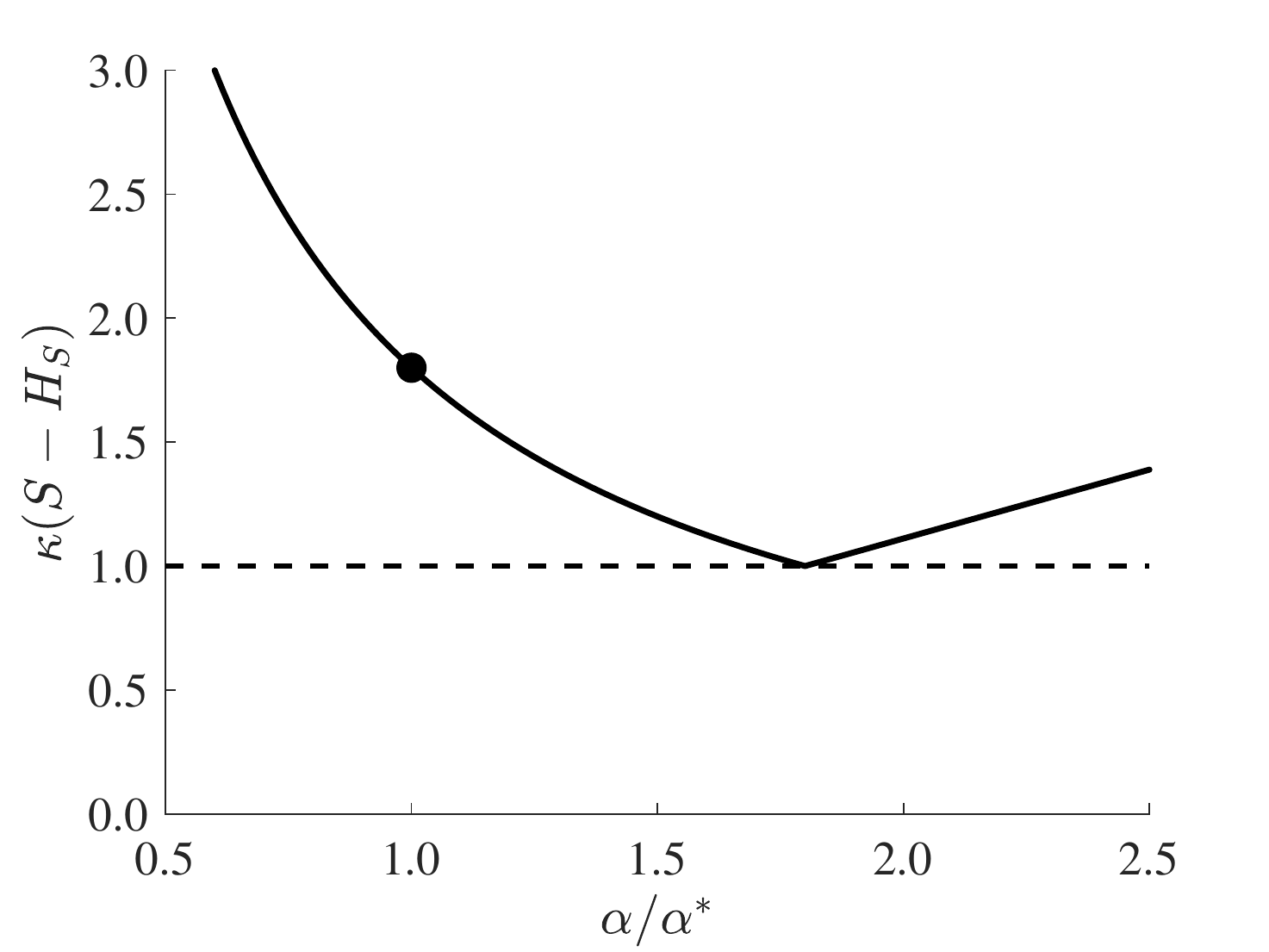}
  \caption{Comparison between the analytic (continuous line) and the algebraic
    (dashed line) macroelement stabilizations.
    On the $y$ axis is the condition number $\kappa(\tS) = \lambda_{\max}(\tS)/
    \lambda_{\min}(\tS)$ of the stabilized Schur complement. Panels from left to the right show
    (i) regular grid and homogeneous material properties, (ii) stretched grid with
    homogeneous material properties, and (iii) regular grid with heterogeneity in the
    material properties.}
  \label{fig:stab_2D_cond}
\end{figure}

To compare the $\alpha$-based and algebraic macroelement stabilizations, in Fig. \ref{fig:stab_2D_cond} we report the behavior of the condition number
$\kappa(\tS) = \lambda_{\max}(\tS)/\lambda_{\min}(\tS)$ of the stabilized Schur
complement for different macroelement configurations. The $x$-axis is normalized with respect to the
``optimal'' value according to Eq. \eqref{eq:stab_2D_alpha}, the dot is the conditioning
arising from the $\alpha$-based stabilization and the dashed line is the conditioning of
the algebraic stabilization. From the left to the right, we have three cases:
\begin{itemize}
  \item Regular grid, with homogeneous material properties. It can be observed that the
    two techniques offer the same optimal result.
  \item Stretched grid (Fig. \ref{fig:stab_2D_sketch}) with $h_{y,{\tau_1}} = 5
    h_{y,{\tau_3}}$.  The material is still homogeneous. In this case, the stabilization
    based on $\alpha$ is not able to get the minimal conditioning, while the algebraic one
    is.
  \item Regular grid (Fig. \ref{fig:active_dof_2D}), with heterogeneity in the Young's modulus: $E_{\tau_1} = E_{\tau_2} = 5E_{\tau_3} = 5 E_{\tau_4}$. Again, the stabilization that uses a scalar value is not able to
    provide the minimal conditioning, while the algebraic one is.
\end{itemize}
Observing the plots in Fig. \ref{fig:stab_2D_cond}, it can be noticed that the
conditioning has a unique minimum when the grid is regular, but there is a set of minima
in the middle case, with a stretched grid. The former occurrence implies that the two
eigenvalues $\lambda_3$ and $\lambda_4$ are the same, while in the latter that they are
different. This finding is consistent with the analytical definition of the Schur
complement eigenvalues (see Eqs. \eqref{eq:schur_eigenpair_nonKernel}), where the aspect
ratio $h_x/h_y$ is a parameter. Using a more refined definition of $\alpha^*$, considering
also the geometric anisotropy, would improve the conditioning number of the $\alpha$-based
$\tS$ for this specific case.  In general, however, as assumptions regarding mesh geometry and material properties are relaxed, an algebraic approach is increasingly appealing in its simplicity.

\begin{figure}
  \centering
  \includegraphics[width=0.40\linewidth]{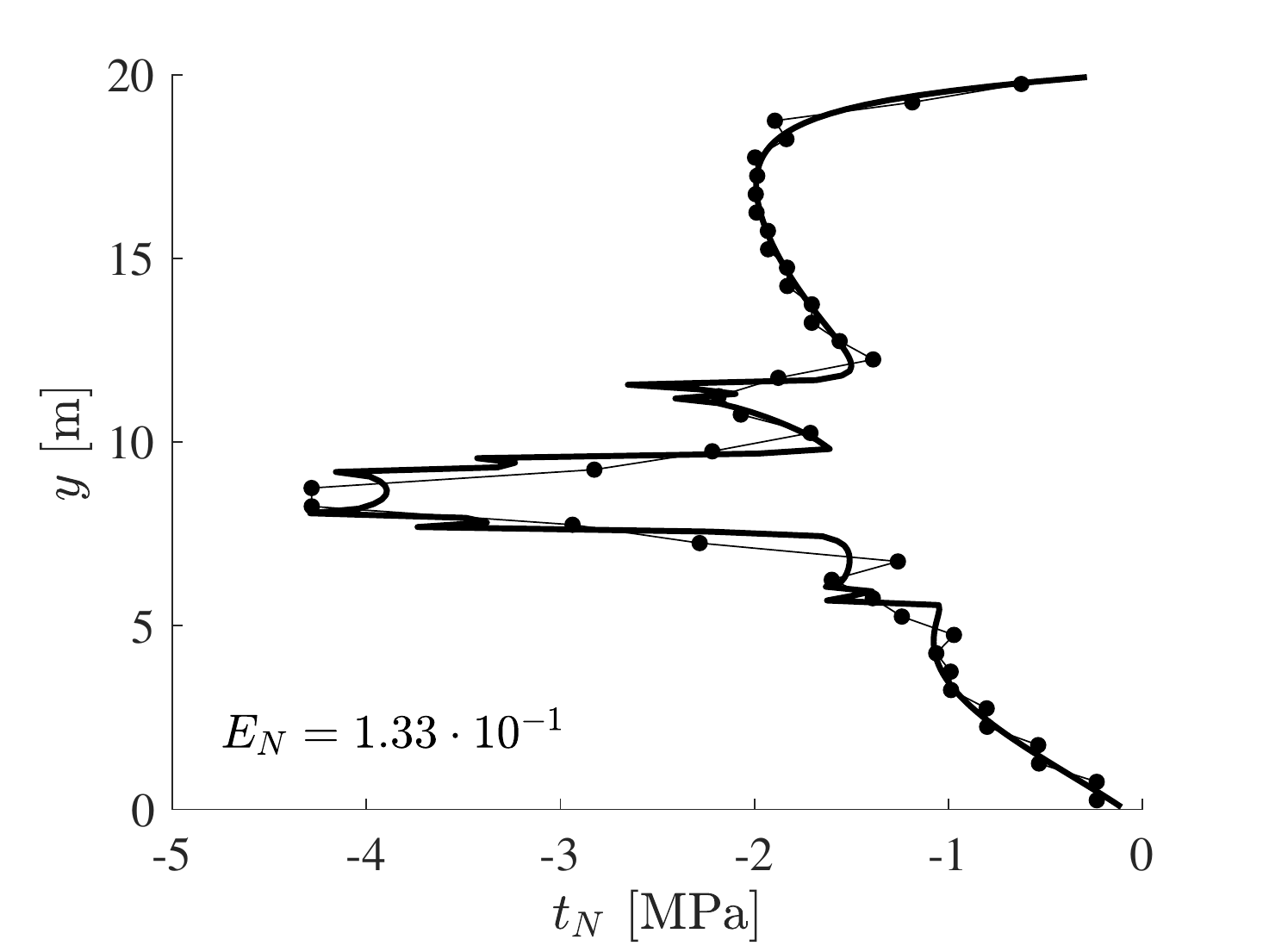}
  \includegraphics[width=0.40\linewidth]{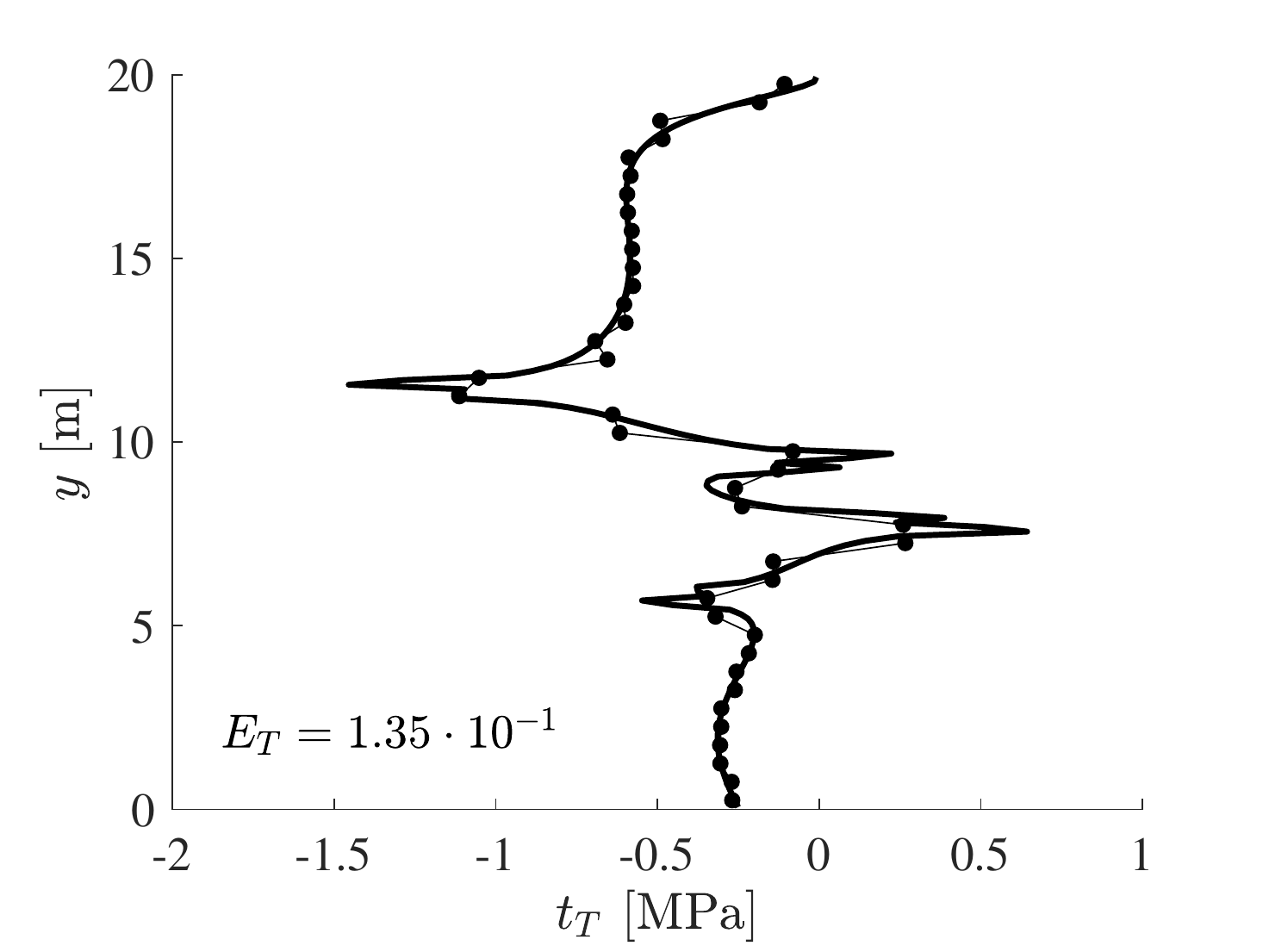}
  \includegraphics[width=0.40\linewidth]{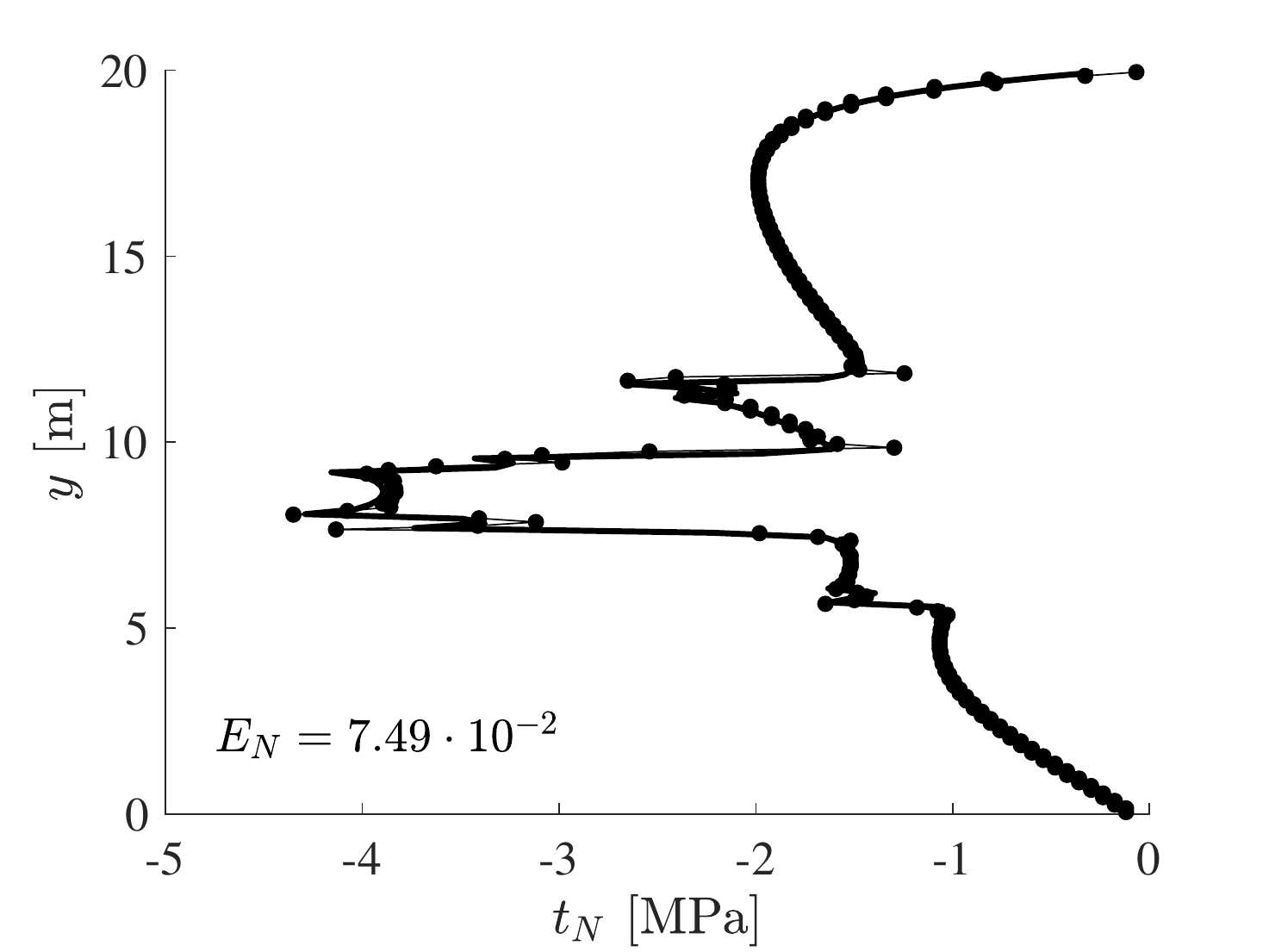}
  \includegraphics[width=0.40\linewidth]{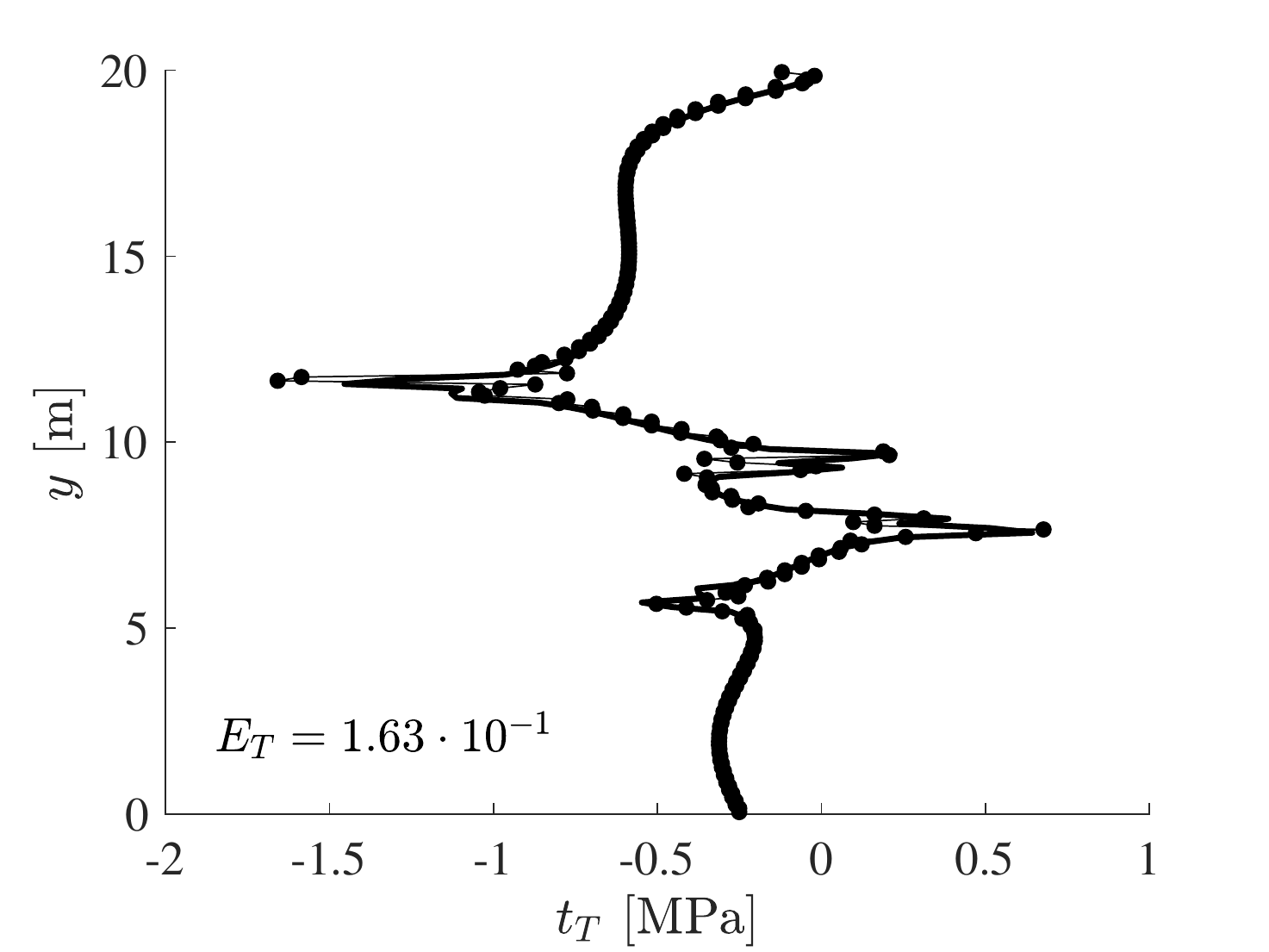}
  \caption{Numerical results for example in Fig. \ref{fig:stab_2D_sketch} using the algebraic macroelement stabilization. First row: base grid, last row:
    vertically refined grid.}
  \label{fig:stab_2D_r2}
\end{figure}

To further test this approach, we focus again on the model problem depicted in Fig. \ref{fig:stab_2D_sketch}.
Numerical results are provided in Fig.
\ref{fig:stab_2D_r2}, with the errors $E_N = 1.33 \cdot 10^{-1}$, $E_T = 1.35 \cdot
10^{-1}$ and $E_N = 7.49 \cdot 10^{-2}$, $E_T = 1.63 \cdot 10^{-1}$ for the base and the
refined grids, respectively. Comparing this technique with the $\alpha$-based approach, we observe quite similar performance.

We conclude this section by discussing the extension to 3D problems.
Matrix $C$ for a 3D macroelement reads
\begin{equation}
  C =
  \begin{bmatrix*}[r]
     \tilde{A}_1 R_3 &  \tilde{A}_2 R_3 &  \tilde{A}_3 R_3 &  \tilde{A}_4 R_3 \\
    -\tilde{A}_1 R_3 & -\tilde{A}_2 R_3 & -\tilde{A}_3 R_3 & -\tilde{A}_4 R_3 
  \end{bmatrix*},
  \label{eq:C_macro_3d_unstructured}
\end{equation}
where the areas $\tilde{A}_i$, $i = \{ 1, \ldots, 4 \}$, are the counterpart of lengths $\tilde{\ell}_1$ and $\tilde{\ell}_2$ (Fig. \ref{fig:unstructured_patch_3D}), and $R_3$ is a 3D rotation matrix.
Using the same arguments as in the 2D case, matrix  $H$ can be obtained by stabilizing the traction jumps across internal edges of the macroelement as $H = \sum_{i = 1}^{4} V_i D V_i^T$, where $D = D_{\widetilde{K},N_1}^{-1} + D_{\widetilde{K},N_2}^{-1} \in \mathbb{R}^{3 \times 3}$ and
\begin{align}
  V_1 &=
  \begin{bmatrix*}
    -\tilde{A}_2 R_3^T \\ 
     \tilde{A}_1 R_3^T \\ 
                     0 \\ 
                     0 
  \end{bmatrix*},
  &
  V_2 &=
  \begin{bmatrix*}
                     0 \\ 
    -\tilde{A}_3 R_3^T \\ 
     \tilde{A}_2 R_3^T \\ 
                     0  
  \end{bmatrix*},
  &
  V_3 &=
  \begin{bmatrix*}
                     0 \\ 
                     0 \\ 
    -\tilde{A}_4 R_3^T \\ 
     \tilde{A}_3 R_3^T 
  \end{bmatrix*},
  &
  V_4 &=
  \begin{bmatrix*}
     \tilde{A}_4 R_3^T \\ 
                     0 \\ 
                     0 \\ 
    -\tilde{A}_1 R_3^T 
  \end{bmatrix*}.
  \label{eq:Hs_macro_3d_unstructured}
\end{align}
\begin{rem}
No explicit assembly process in needed in the construction of $V$ and $D$ both in 2D and 3D. 
A local gathering from the global stiffness and coupling matrix is enough to extract the required blocks.
\end{rem}

\subsection{Algebraic global stabilization}
\label{sec:global_stab}
In this section, we extend the algebraic approach of Sec. \ref{sec:macroelem_noA} to overcome the key limitation of any macroelement-based approach, the topological restriction that the mesh be partitioned into macroelements in the first place.
As previously stated, the key point is to stabilize the jump between two interface elements. Rather than doing this for macroelement-internal edge alone, we now simply stabilize \emph{all} possible jumps between each pair of interface elements.  

\begin{figure}[t]
  \centering
  \includegraphics[width=0.20\linewidth]{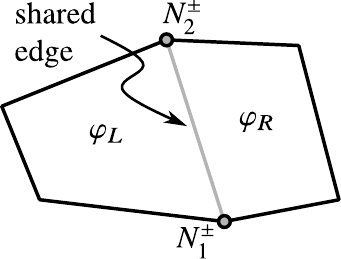}
  \caption{Sketch of two elements sharing an edge. Common nodes on both surfaces of the
    fracture are highlighted.}
  \label{fig:sketch_gloS}
\end{figure}

To do so, it is sufficient to compute the matrix $\tC_\text{loc}$ for all degrees-of-freedom shared between two elements discretizing the fracture, gather the local matrix $K_\text{loc}$, and then assemble individual $H_\text{loc}$ contributions into a global stabilization matrix $H$.  Algorithm \ref{alg:globalStab} describes the approach for the 2D case---see Figure~\ref{fig:sketch_gloS}---using a Matlab-style notation.  Note that the proposed algorithm assumes, implicitly, that fractures are planar. Also recall that multiple DOFs are associated to each geometric object (nodes and faces) and thus the gather/scatter indexing refers to vector sets of component indices.

We note that the resulting stabilization matrix has the same stencil as a two-point flux approximation,
having been assembled using points (edges) shared between pairs of elements in 2D (3D). For a quadrilateral (hexahedral)
mesh, it is a 3-point (5-point) block Laplacian in 2D (3D). By
\emph{block} Laplacian, we mean that each matrix super-entry is composed of a $2 \times 2$ or $3 \times 3$
dense block, according to the space dimension, with one row/column for each component of the
traction vector.

This Laplacian stencil can be easily built using the pre-assembled global $C^T$ matrix, because its pattern provides the DOF-connectivity between faces and nodes for the fracture discretization. In 3D, we must take care
to avoid the assembly of elements sharing just one node, instead of an edge.  We can do so by computing the
product $C_1^T C_1$, where $C_1$ is a matrix with the sparsity pattern of $C$ but all
non-zero entries set to one.  Each row/column entry of the product matrix is then the number of shared nodes between the respective elements.  We may then readily filter out any element pairs that share just one node.

\begin{algorithm}[t]
\caption{Algebraic global stabilization}
\begin{algorithmic}[1]
  \ForAll {shared edges}
    \State find connected elements $\varphi_L$, $\varphi_R$ and nodes $N_1^+$, $N_1^-$, $N_2^+$,
      $N_2^-$ \Comment{See Fig. \ref{fig:sketch_gloS}}
    \State gather local $C_{\text{loc}} =
      C([\mathcal{N}_{u,1^+},\mathcal{N}_{u,1^-},\mathcal{N}_{u,2^+},\mathcal{N}_{u,2^-}],[\mathcal{N}_{t,L},\mathcal{N}_{r,R}])$ 
    \State form $\tC$ by swapping columns of $C_{\text{loc}}$ and changing sign
    \State gather local $K_{\text{loc}} =
      K([\mathcal{N}_{u,1^+},\mathcal{N}_{u,1^-},\mathcal{N}_{u,2^+},\mathcal{N}_{u,2^-}],
        [\mathcal{N}_{u,1^+},\mathcal{N}_{u,1^-},\mathcal{N}_{u,2^+},\mathcal{N}_{u,2^-}])$ 
    \State form $D_{\widetilde{K}} = \text{diag}(K_{\text{loc}})$  
    \State compute stabilization contribution $H_{\text{loc}} = \tC^T
      D_{\widetilde{K}}^{-1} \tC$ 
    \State add contribution to global stabilization matrix
      $H([\mathcal{N}_{t,L},\mathcal{N}_{r,R}],[\mathcal{N}_{t,L},\mathcal{N}_{r,R}] = H_{\text{loc}}$
  \EndFor
\end{algorithmic}
\label{alg:globalStab}
\end{algorithm} 

\begin{figure}
  \centering
  \includegraphics[width=0.40\linewidth]{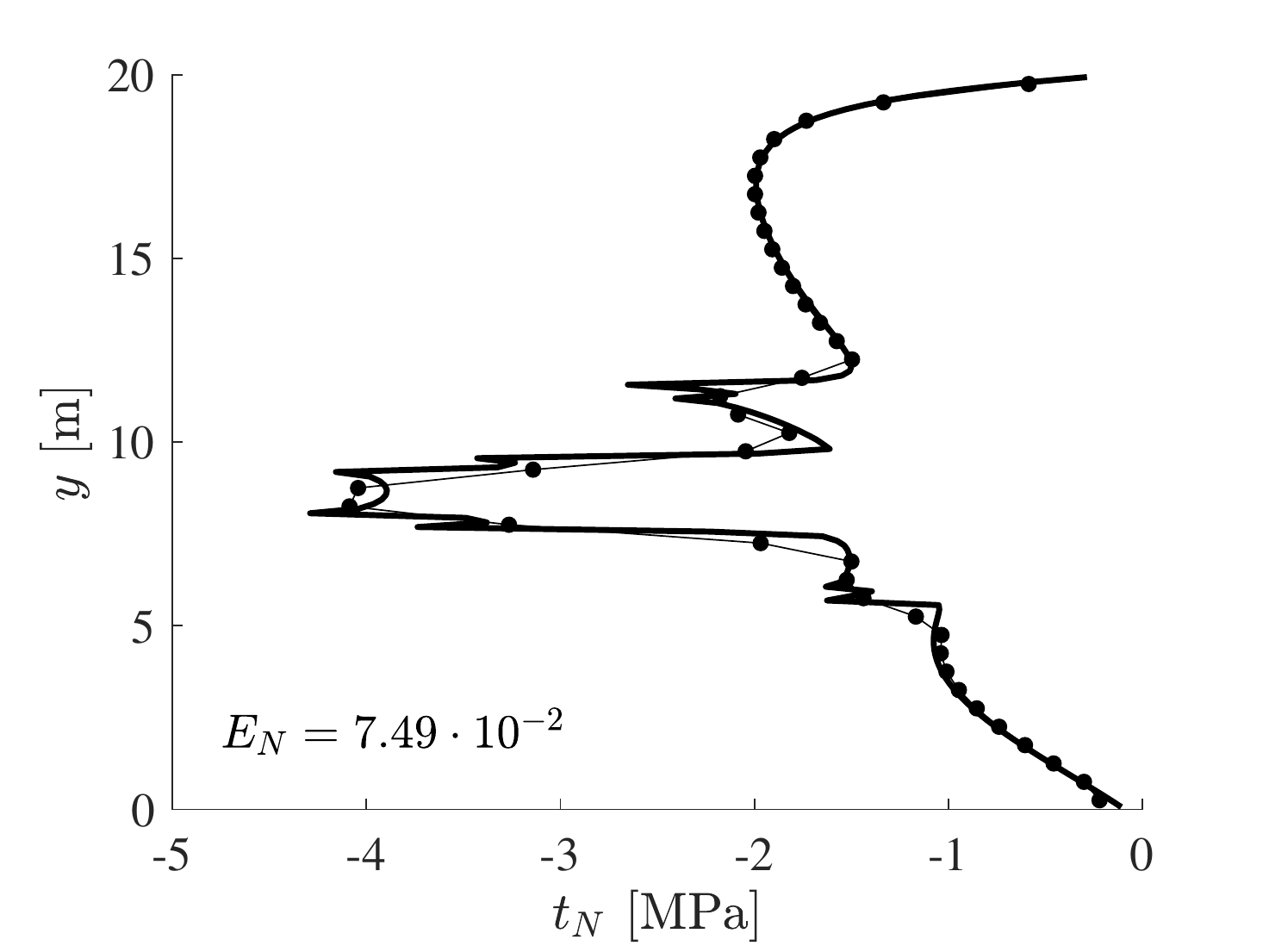}
  \includegraphics[width=0.40\linewidth]{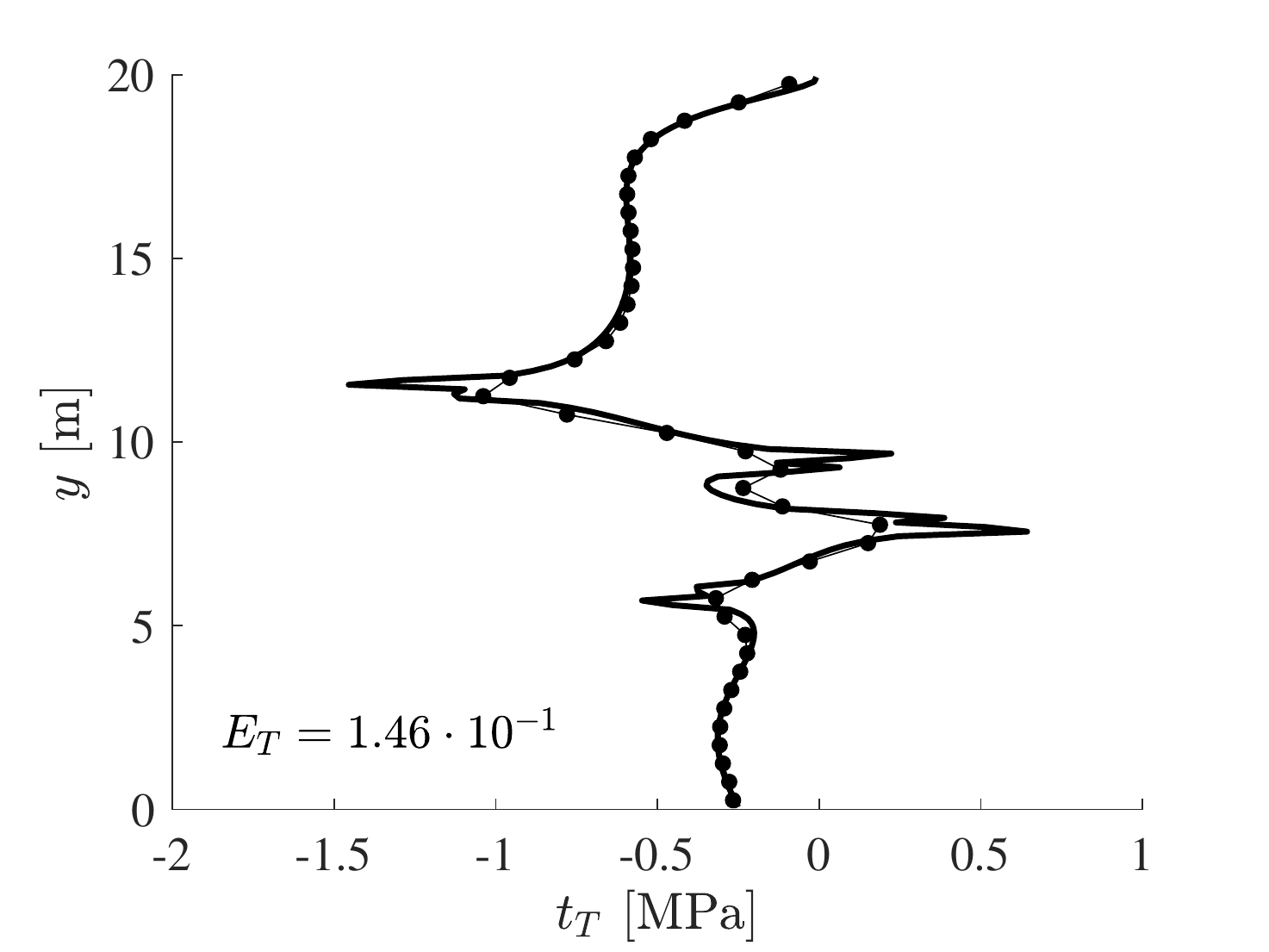}
  \includegraphics[width=0.40\linewidth]{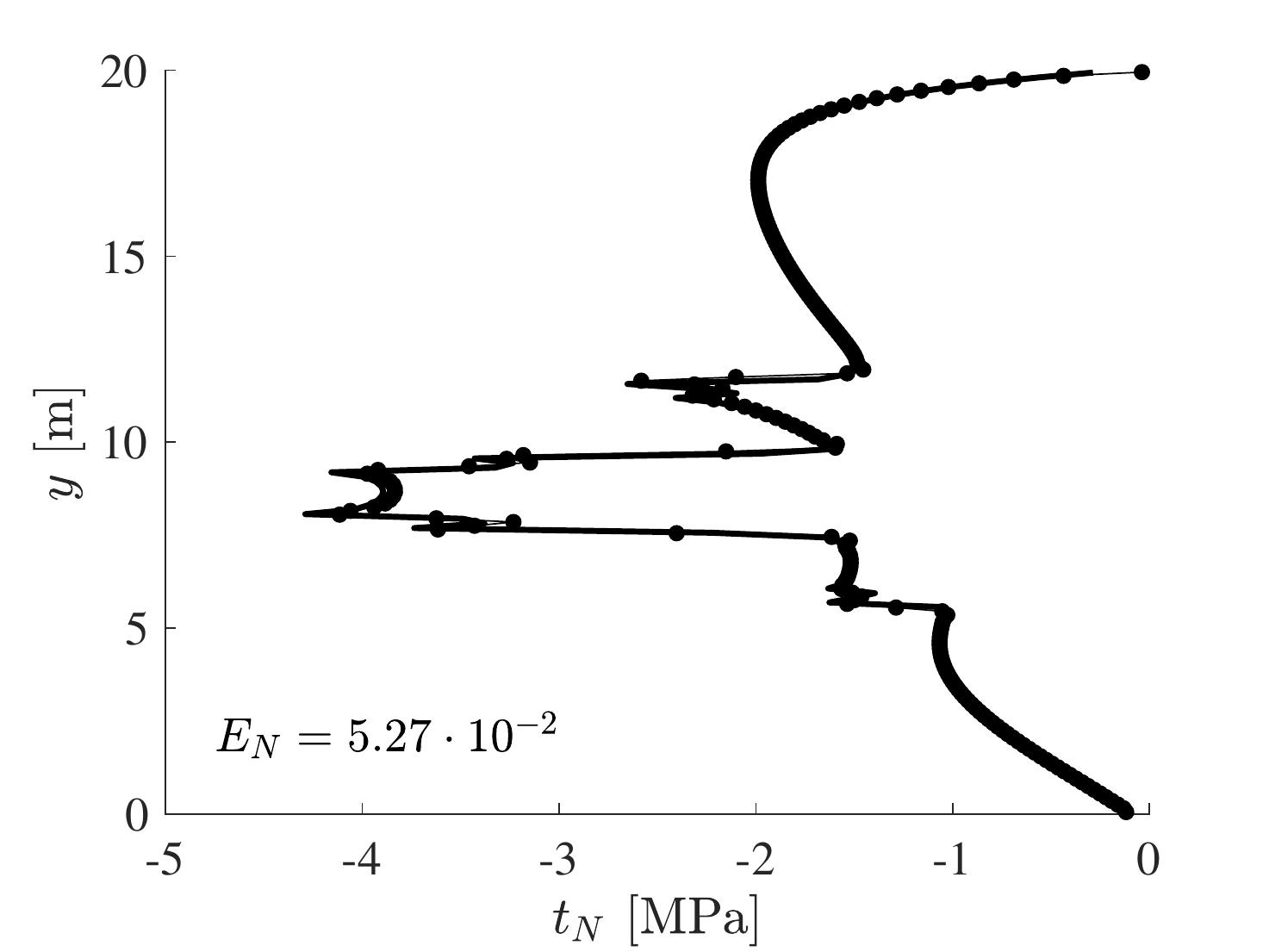}
  \includegraphics[width=0.40\linewidth]{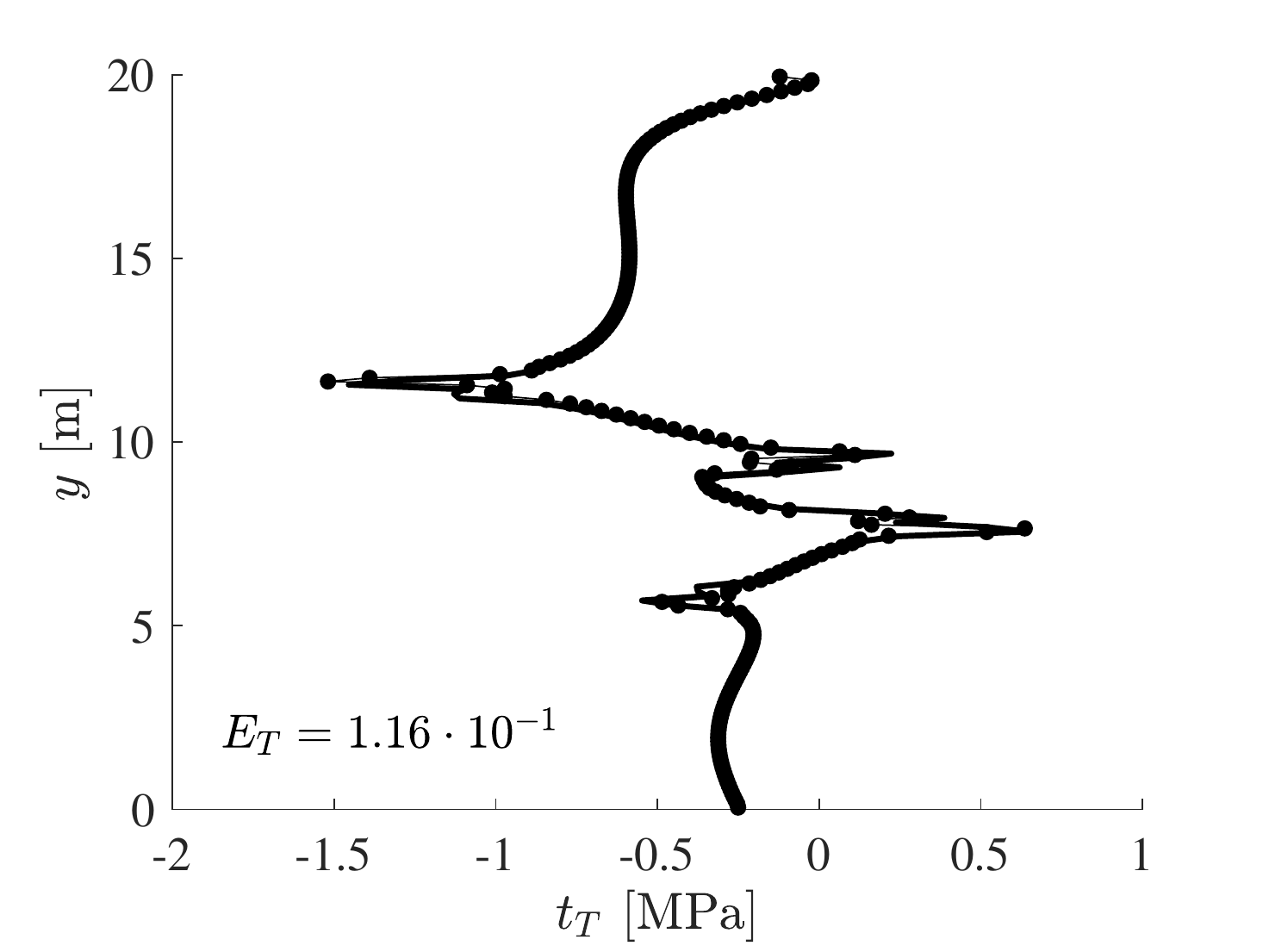}
  \caption{Numerical results for example in Fig. \ref{fig:stab_2D_sketch} using the algebraic global stabilization. First row: base grid, last row:
    vertically refined grid.}
  \label{fig:stab_2D_r3}
\end{figure}

As with other approaches, Fig. \ref{fig:stab_2D_r3} shows the performance on this method on our comparison problem.  The resulting errors are $E_N = 7.49 \cdot 10^{-2}$, $E_T = 1.46 \cdot
10^{-1}$ and $E_N = 5.27 \cdot 10^{-2}$, $E_T = 1.16 \cdot 10^{-1}$ for the base and the
refined grids, respectively. Comparing this technique with the other two approaches, this method produces smaller errors.
\begin{figure}
  \centering
  \includegraphics[width=0.33\linewidth]{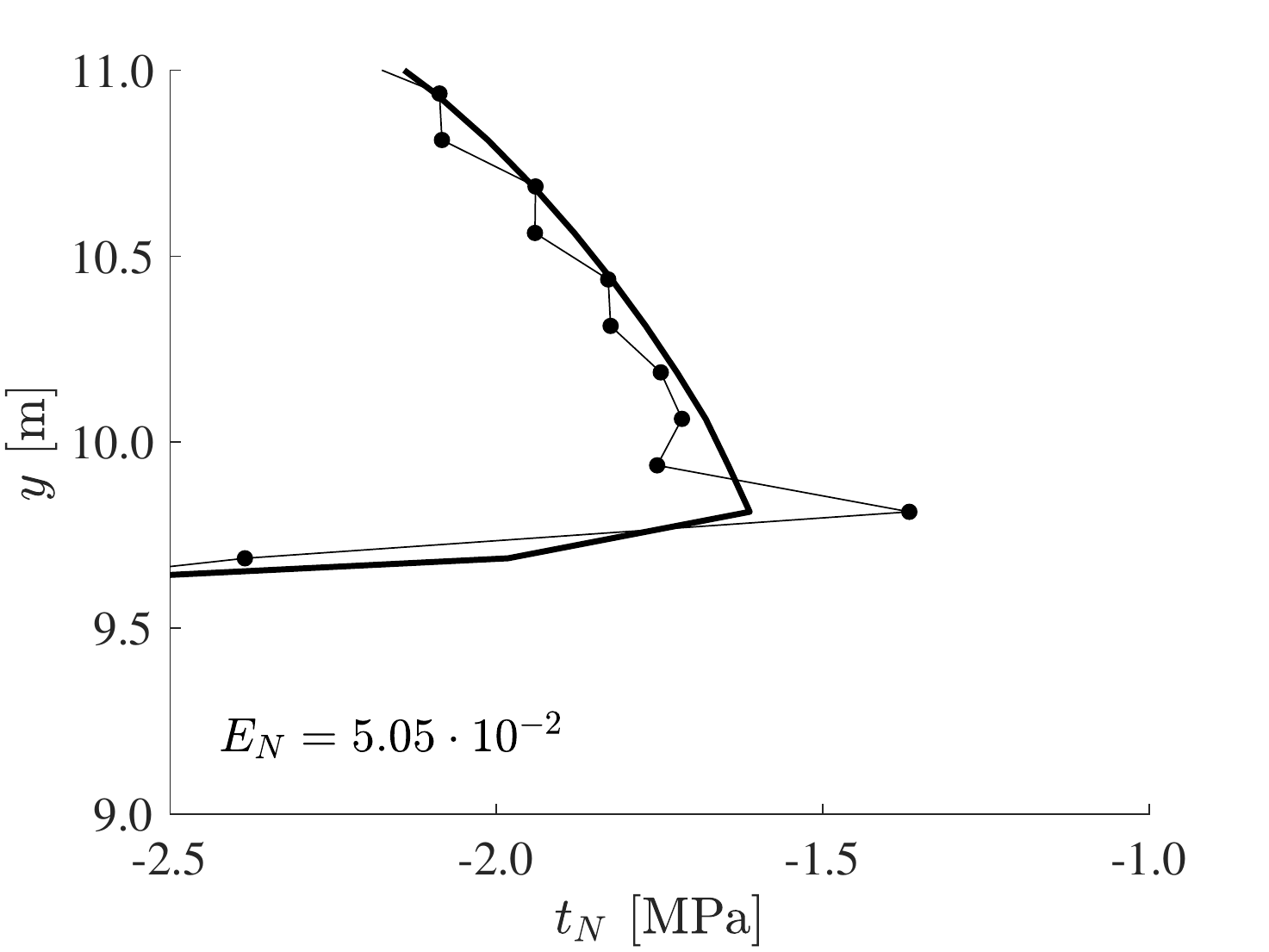}
  \includegraphics[width=0.33\linewidth]{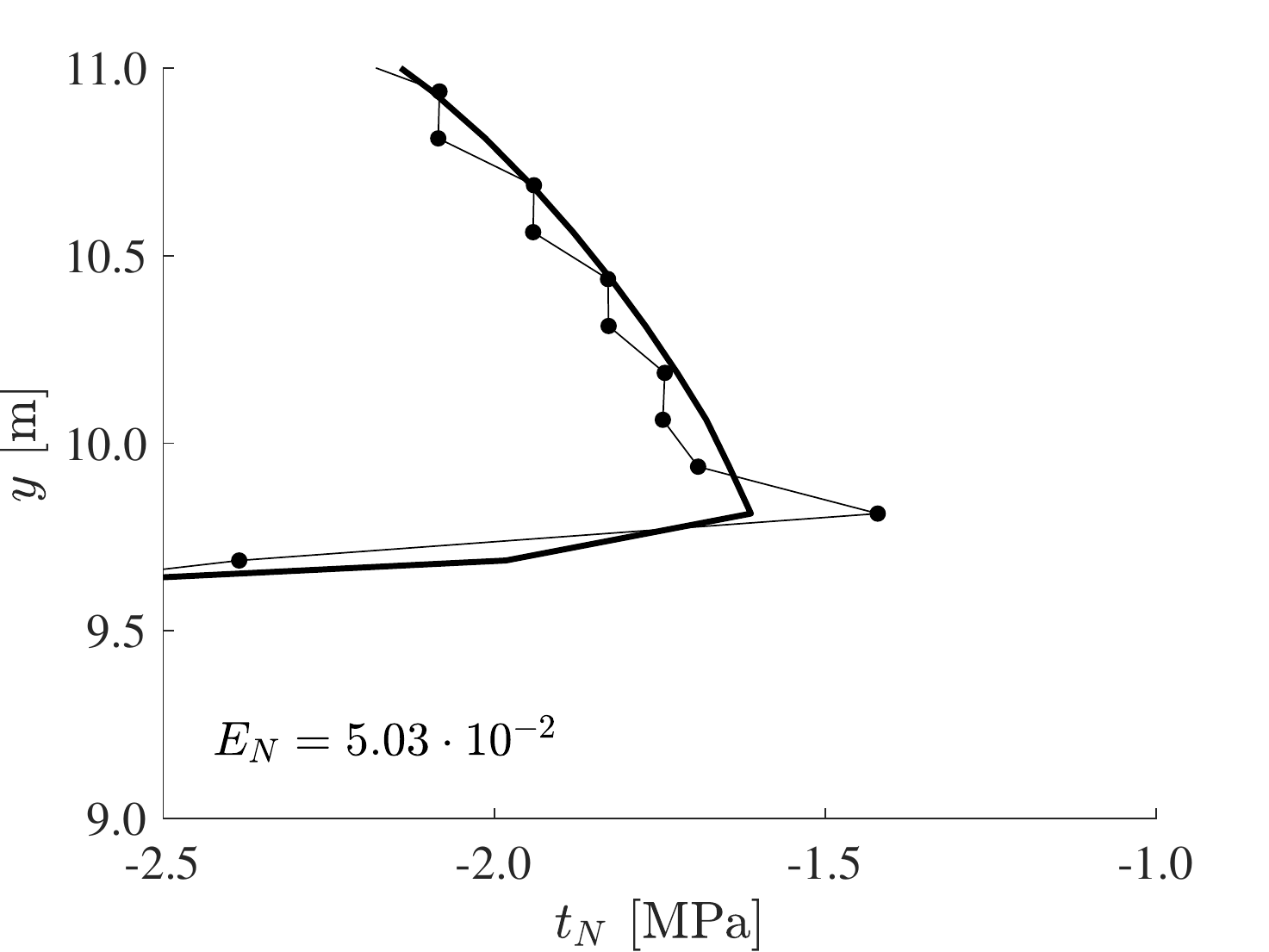}
  \includegraphics[width=0.33\linewidth]{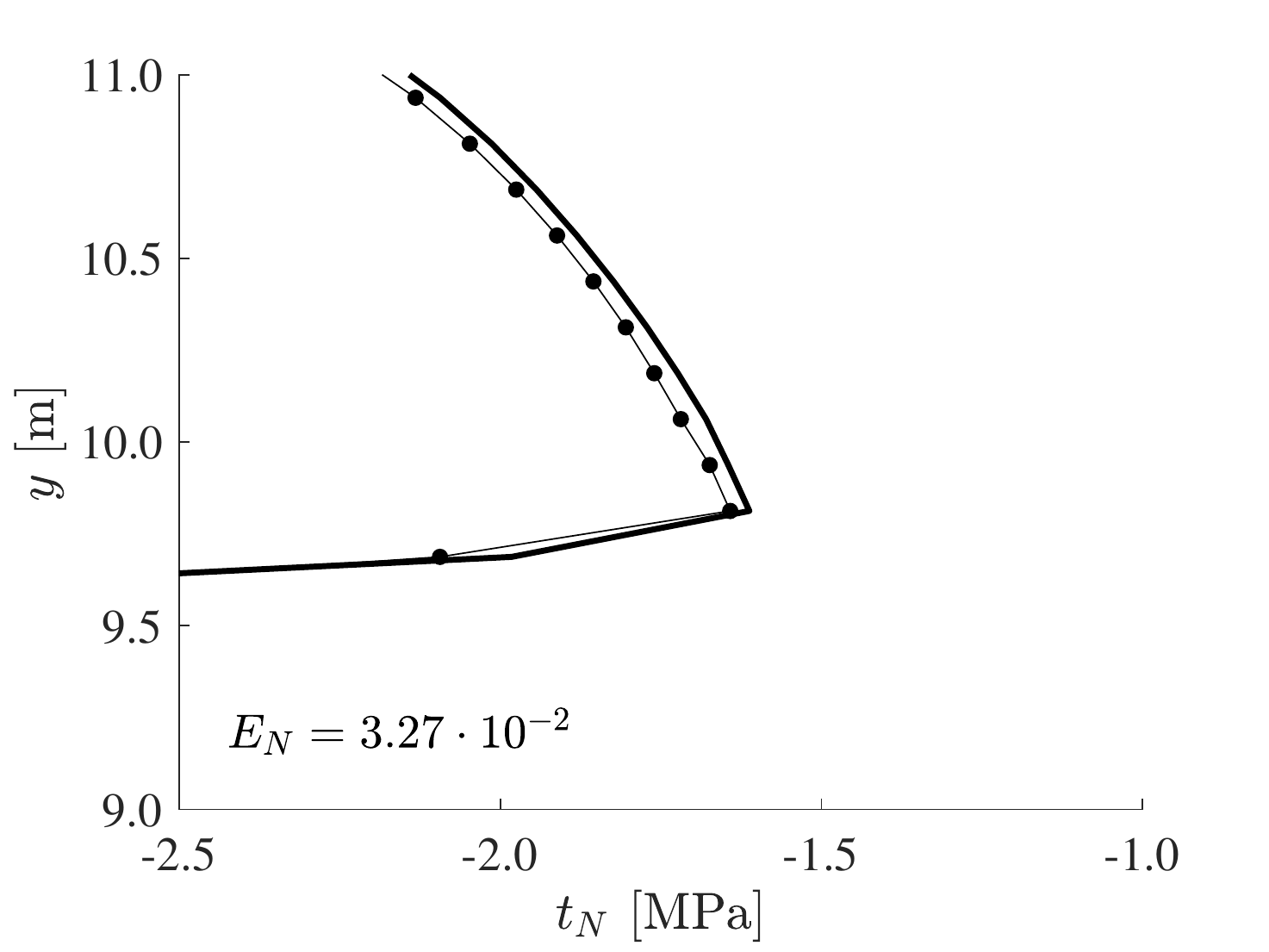}
  \caption{Comparison of different stabilization techniques for a refined regular grid.
    From the left to the right: analytic macroelement, algebraic macroelement, and algebraic global
    stabilization approach.}
  \label{fig:stab2D_comp}
\end{figure}
To better highlight the different behavior of the three methods, we solve again the model problem of Fig. \ref{fig:stab_2D_sketch} using a more refined grid with $64 \times 160$ elements. We report results in terms of normal component $t_N$
only, being representative of overall performance. In Fig. \ref{fig:stab_2D_sketch}, the black circle
identifies the portion of the fracture represented in Fig. \ref{fig:stab2D_comp}. The
global relative error $E_N$, shown therein, is computed as before.

Given its straightforward computation and good performance, we adopt the algebraic global stabilization as our default approach.

\begin{rem}\label{rem:sign_schur}
We observe that, for all three approaches, the stabilization matrix $H$ is
symmetric and positive semi-definite (SPSD). We recall that, being the Schur complement of a
negative definite matrix, it is added with the minus sign, i.e., $\tS = S - H$.
\end{rem}

\begin{rem}
The idea of using a global stabilization based on jump penalization is not completely new, having been suggested for other discretization and interpolation schemes
\cite{BurHan10,puso2012embedded,puso2015embedded}.
\end{rem}

\subsection{Inclusion of friction and fluid flow}

The methods described for the stick problem form the building blocks for stabilizing the full physical model, allowing for frictional sliding and fluid flow.  In particular, all that is required is to modify the original system~(\ref{eq:jac_sys}) to its stabilized version,
\begin{align}
  \begin{bmatrix}
    %
    %
    A_{uu} &
    A_{uS} &
    A_{uN} &
    A_{uT} &
    A_{uO} &
    A_{up} \\
    %
    %
    %
    A_{Su} & -H_{SS} & -H_{SN} & 0 & 0 & 0 \\
    %
    %
    A_{Nu} & -H_{NS} & -H_{NN} & 0 & 0 & 0 \\
    %
    %
    A_{Tu} &
    0 &
    A_{TN} &
    A_{TT} &
    0 &
    0 \\
    %
    %
    0 & 0 & 0 & 0 & A_{OO} & 0 \\
    %
    %
    A_{pu} & 0 & 0 & 0 & 0 & A_{pp}+\frac{1}{\Delta t} H_{pp} \\
  \end{bmatrix}_n^{\ell,(k)}
  \begin{bmatrix}
    \delta \Vec{u} \\
    \delta \Vec{t}_S \\
    \delta \Vec{t}_N \\
    \delta \Vec{t}_T \\
    \delta \Vec{t}_O \\
    \delta \Vec{p}
  \end{bmatrix}
  = -
  \begin{bmatrix}
    \Vec{r}_u \\
    \Vec{r}_{S} - H_{SS}\Vec{t}_S- H_{SN}\Vec{t}_N \\
    \Vec{r}_{N} - H_{NS}\Vec{t}_S- H_{NN}\Vec{t}_N\\
    \Vec{r}_{T} \\
    \Vec{r}_{O} \\
    \Vec{r}_p + \frac{1}{\Delta t} H_{pp} \Vec{p}
  \end{bmatrix}_n^{\ell,(k)} \,,
\end{align}
where five stabilizing sub-matrices $H_{ij}$ have been added to particular components of
the multiphysics problem. In particular, the traction components $\Vec{t}_S$ and
$\Vec{t}_N$ and the pressure field $\Vec{p}$ must all receive stabilizing contributions.
Stability is primarily provided by the matrices $H_{SS}$, $H_{NN}$, and $H_{PP}$, which
are assembled on macroelements using the same procedures as described in the previous
section, but now applied to different traction and pressure components as warranted. The
off-diagonal terms $H_{SN}$ and $H_{NS}$ are necessary to capture contributions from
stabilized edges between elements in a ``mixed'' state---that is, where an
element in a stick state is adjacent to one in frictional sliding state. We observe that
$\vec{t}_T$---i.e., the friction component of the traction on
$\Gamma_f^{\text{slip}}$---is a function of the relative displacement rate $\Delta
\vec{g}_T$ and therefore does not need to be stabilized. Similarly, for
$\Gamma_t^{\text{open}}$, there is no need for stabilization, as the traction $\Vec{t}_O$
is \emph{a priori} known to be zero. As before, because the matrices $A_{TT}$ and $A_{OO}$
are diagonal, a Schur-complement reduction can be performed to eliminate $\Vec{t}_T$ and
$\Vec{t}_O$. The reduced system will closely resemble Eq.~(\ref{eq:jac_sys_reduced}), but
with the stabilizing blocks above included in the appropriate places.

\section{Analytical benchmarks for contact mechanics without fluid flow}
\label{sec:anal_cont_mech}
In this section, we validate the formulation, the discretization, and the stabilization
strategies for  \textit{pure} contact mechanics problems, without fluid flow.

\subsection{Analytical benchmark for shear behavior: Constant solution}
\label{sec:borja_2D}

\begin{figure}
  \begin{subfigure}{0.3\textwidth}
    \centerline{\includegraphics[height=5cm]{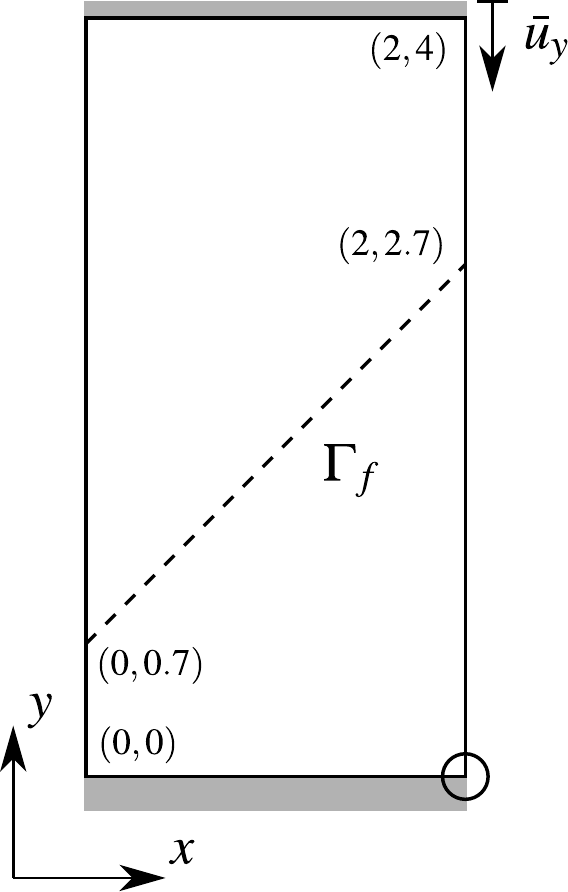}}
    \caption{}  
    \label{fig:borja_sketch}
  \end{subfigure}
  \hfill
  \begin{subfigure}{0.3\textwidth}
    \centerline{\includegraphics[height=5cm]{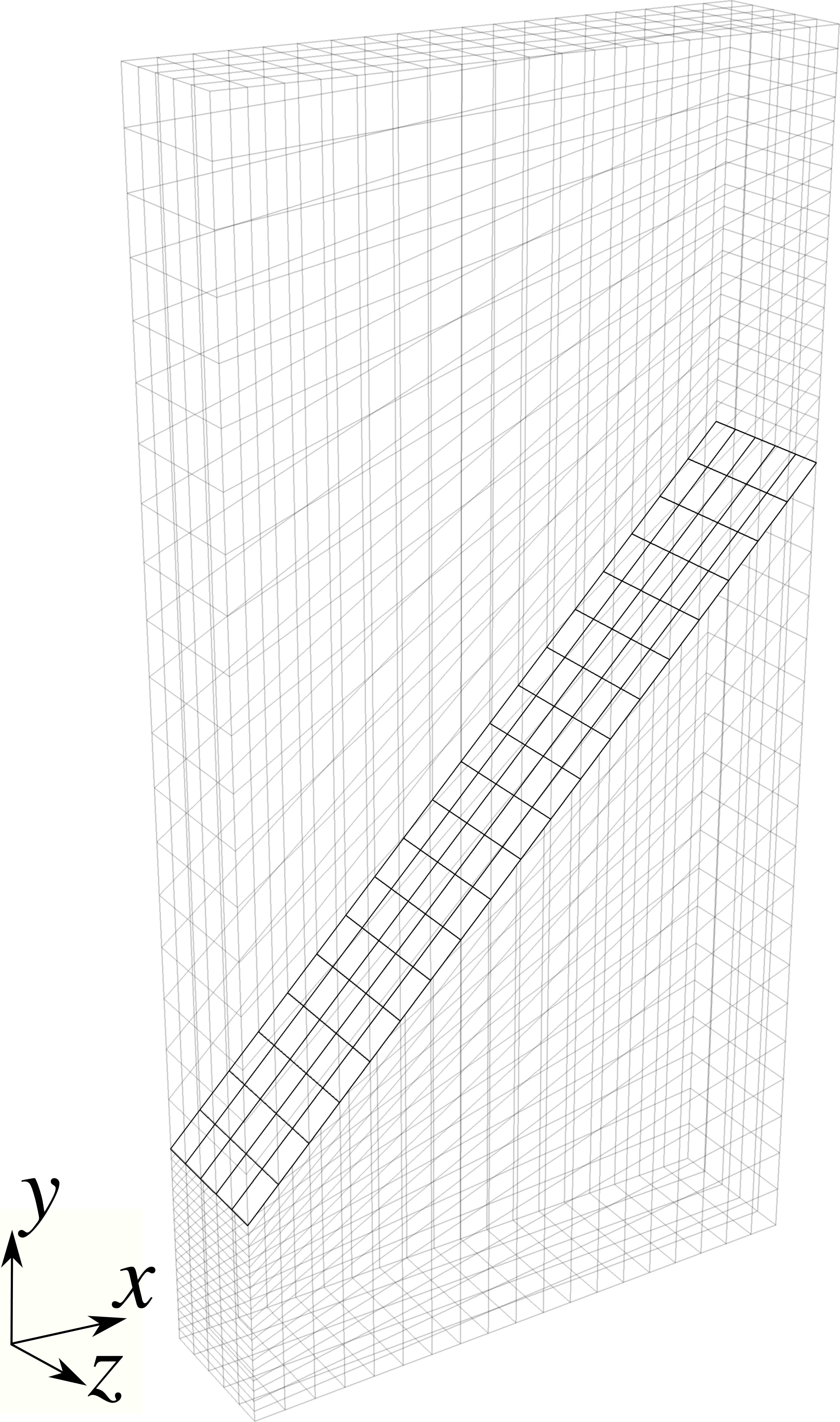}}
    \caption{}  
    \label{fig:borja_mesh}
  \end{subfigure}
  \hfill
  \begin{subfigure}{0.3\textwidth}
    \centerline{\includegraphics[height=5cm]{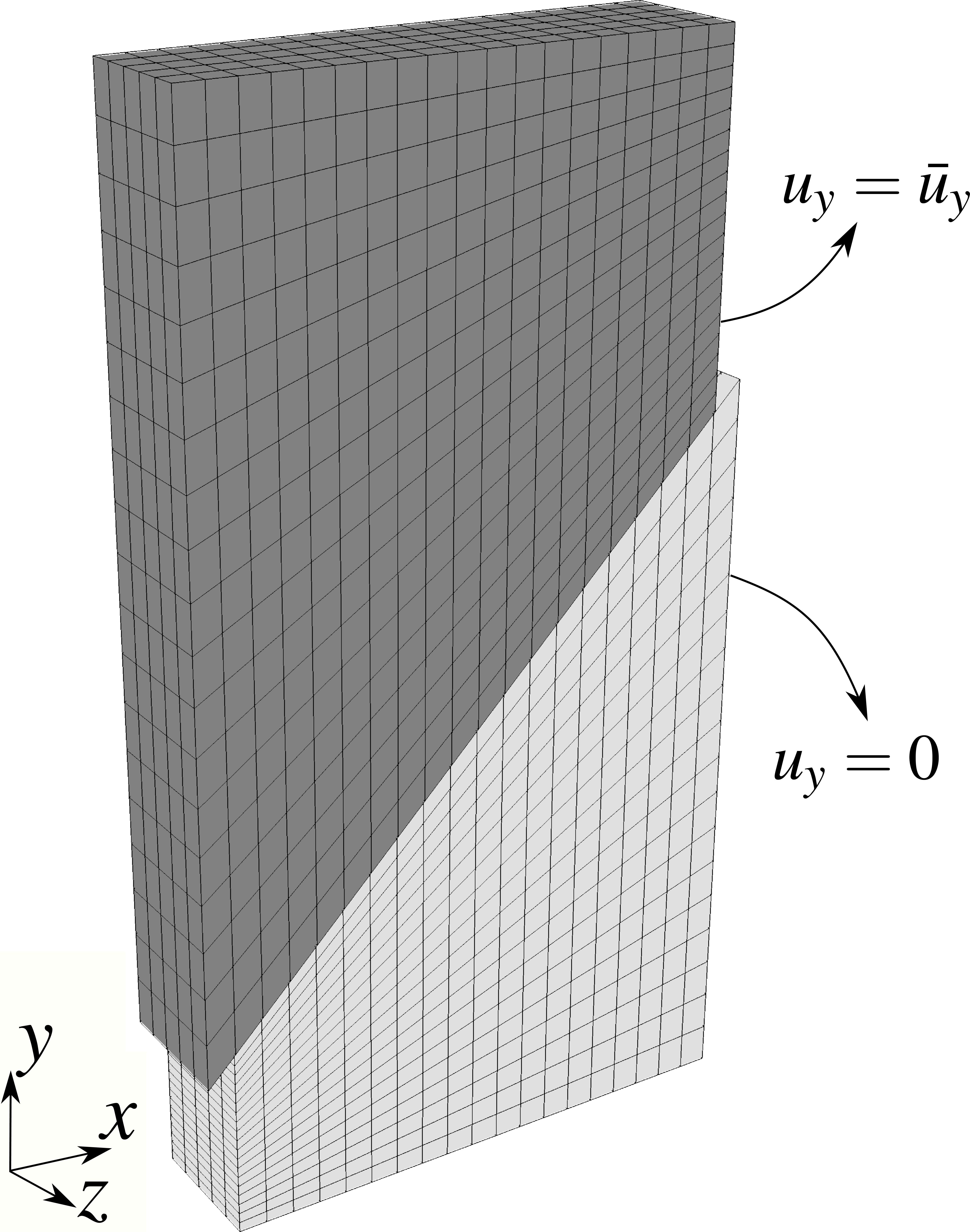}}
    \caption{}  
    \label{fig:borja_sol_contour}
  \end{subfigure}   
  \caption{Analytical benchmark for shear behavior \cite{borja2008assumed}: (a) sketch of the domain with coordinates expressed in meters; (b) computational 3D mesh; (c) contour plot for the simulated vertical displacement (magnified).}
\end{figure}

This analytical benchmark was originally proposed in \cite{borja2008assumed}. The
representation of the 2D model domain is reported in Fig. \ref{fig:borja_sketch}, where
the lower boundary is $y$-constrained, the circled corner is constrained also in the $x$
direction and on the upper boundary a uniform displacement is imposed ($\bar{u}_y = 0.10
\text{ m}$). The material is homogeneous with elastic parameters $E = 5000 \text{ MPa}$
and $\nu = 0.25$. Coulomb's frictional parameters are $c = 0$ and $\theta = 5.71^\circ$,
such that the friction coefficient is $0.1$. The solution is a constant sliding on the
fracture of value $g_T = ||\tensorTwo{g}_T||_2 = 0.1 \sqrt{2} \approx 0.1414 \text{ m}$.

We highlight that the simulation is carried out on a 3D mesh (Fig. \ref{fig:borja_mesh}), because in a 2D setting the
direction of the shear vector is known; thus, the Coulomb frictional contact condition is
not needed and the shear behavior is linear.
To test the nonlinearity introduced by the dependency of the direction on the relative
displacement rate, when the failure condition is matched and the fracture slides, we need
to work in a 3D setting.
To match a 2D
analytical solution, we respect the plane strain assumption and build the domain extruding
the surface shown in Fig. \ref{fig:borja_sketch} by $0.5 \text{ m}$ and $z$-constraining
the two surfaces parallel to the $x$-$y$ plane.

\begin{figure}
  \hfill
  \begin{subfigure}{0.4\textwidth}
    \centerline{\includegraphics[width=\linewidth]{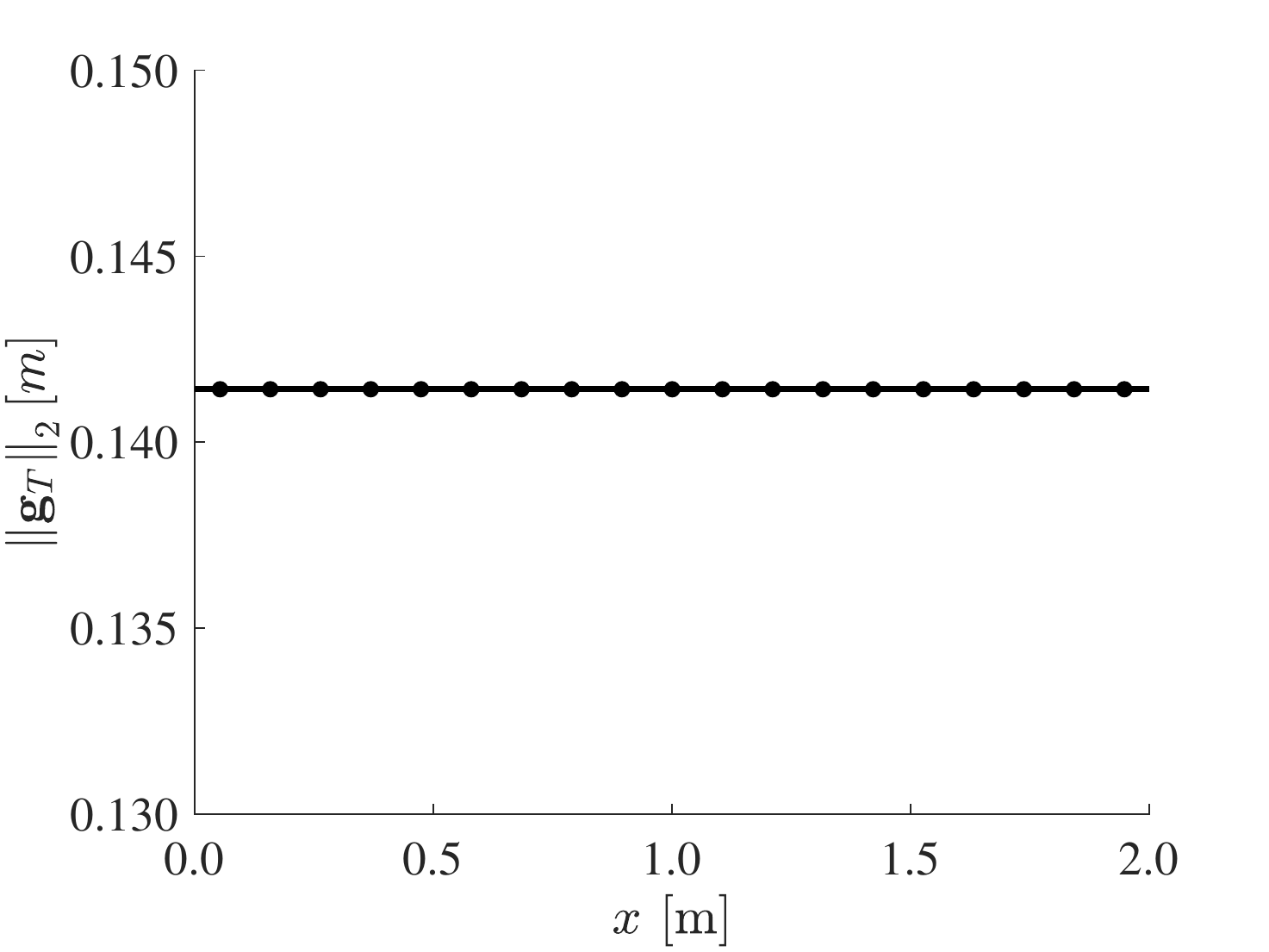}}
    \caption{}  
    \label{fig:borja_sol}
  \end{subfigure}
  \hfill
  \begin{subfigure}{0.4\textwidth}
    \centerline{\includegraphics[width=\linewidth]{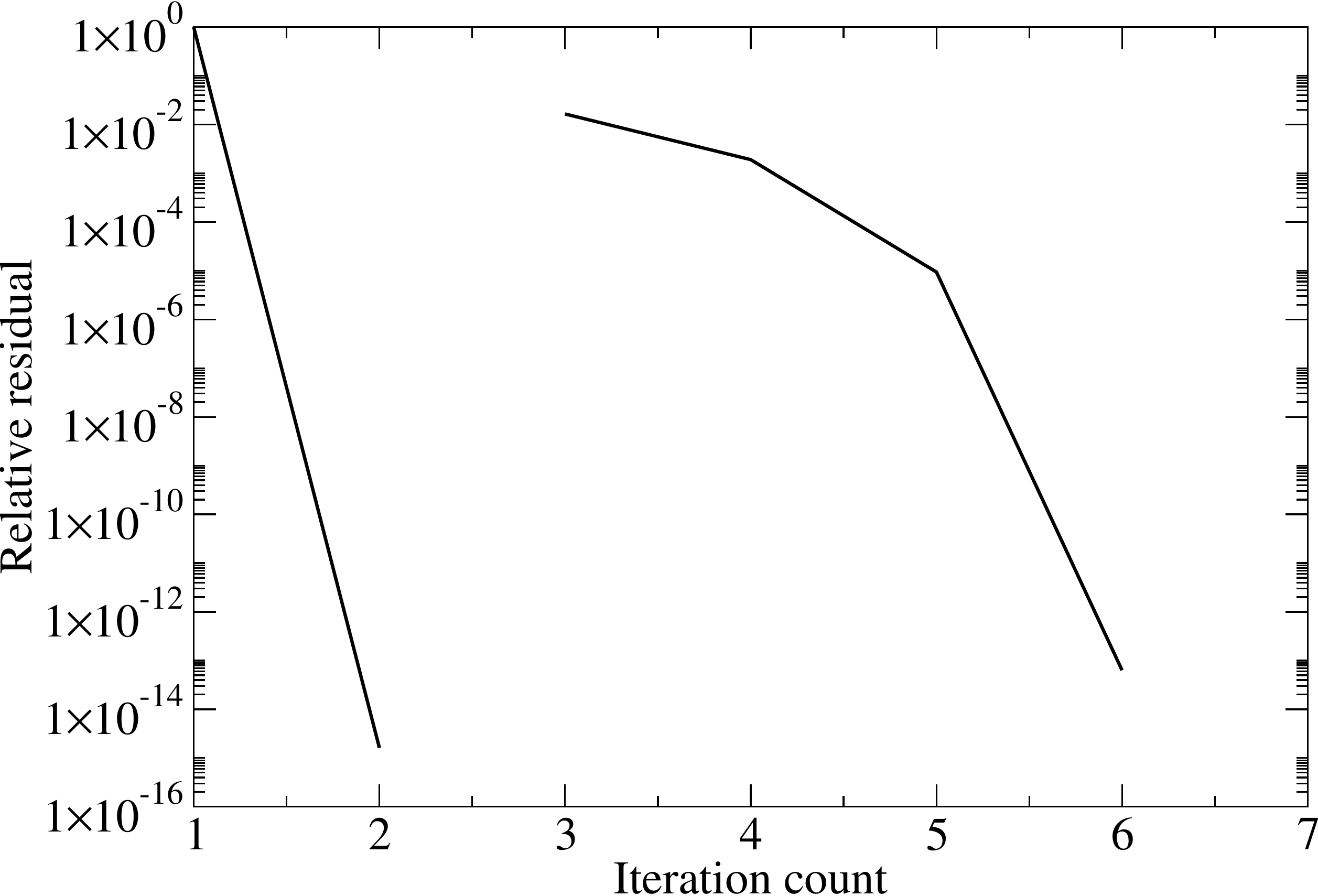}}
    \caption{}  
    \label{fig:borja_convprof}
  \end{subfigure}
  \hfill\null
  \caption{Analytical benchmark for shear behavior: (a) slip on the fracture plane; (b) nonlinear convergence profile.}
\end{figure}

From Figs. \ref{fig:borja_sol_contour}-\ref{fig:borja_sol}, it is clear that the model is able to match the analytic
linear solution even with a small number of elements (4680 nodes, 3610 hexahedra and 95
quadrilaterals for the fracture).
We report also the convergence 
profile of the nonlinear solution algorithm \ref{alg:activeset} in Fig.
\ref{fig:borja_convprof}.
In this specific case, after the elastic step,
i.e., the first solution with $\Gamma_f = \Gamma_f^ \text{stick}$, the final configuration
of $\Gamma_f$ is obtained and no other outer loop iterations are required.

\subsection{Single crack under compression}
\label{sec:phan_2D}

\begin{figure}
  \hfill
  \begin{subfigure}[b]{0.35\textwidth}
    \centerline{\includegraphics[width=\linewidth]{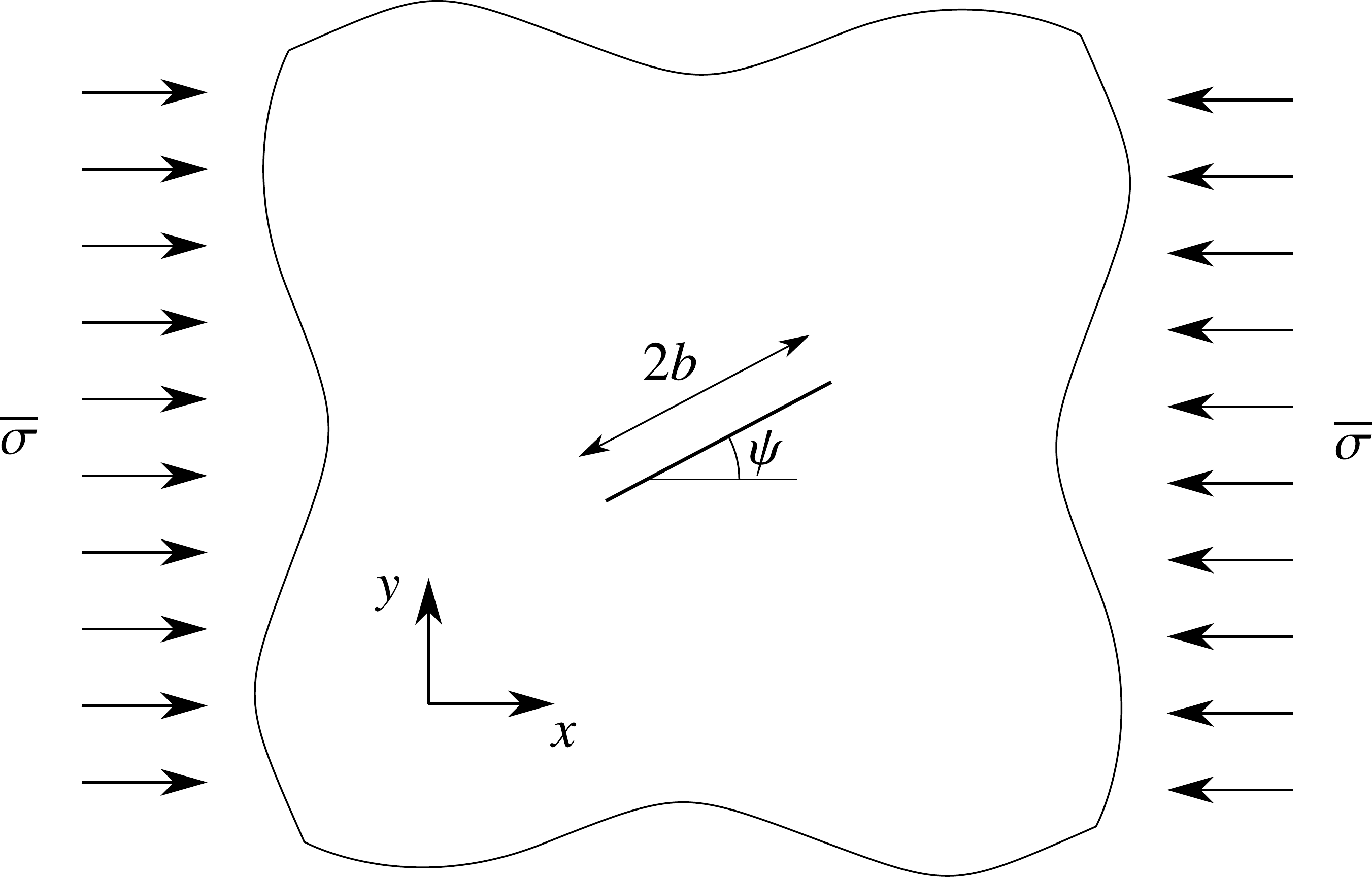}}
    \caption{}  
    \label{fig:phan_geom}
  \end{subfigure}
  \hfill
    \begin{subfigure}[b]{0.22\textwidth}
    \centerline{\includegraphics[width=\linewidth]{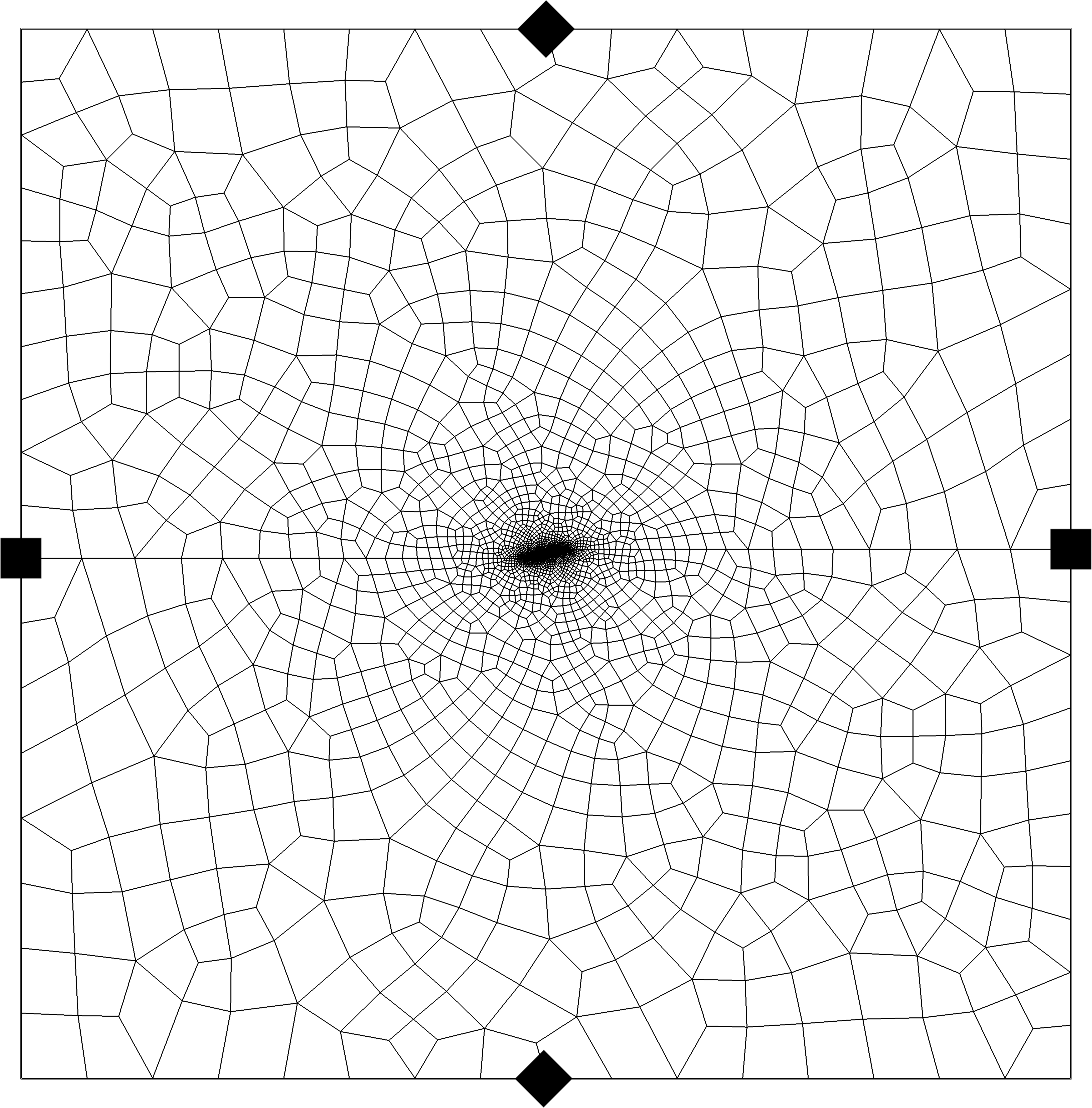}}
    \caption{}  
    \label{fig:phan_grid}
  \end{subfigure}
  \hfill
  \begin{subfigure}[b]{0.35\textwidth}
    \centerline{\includegraphics[width=\linewidth]{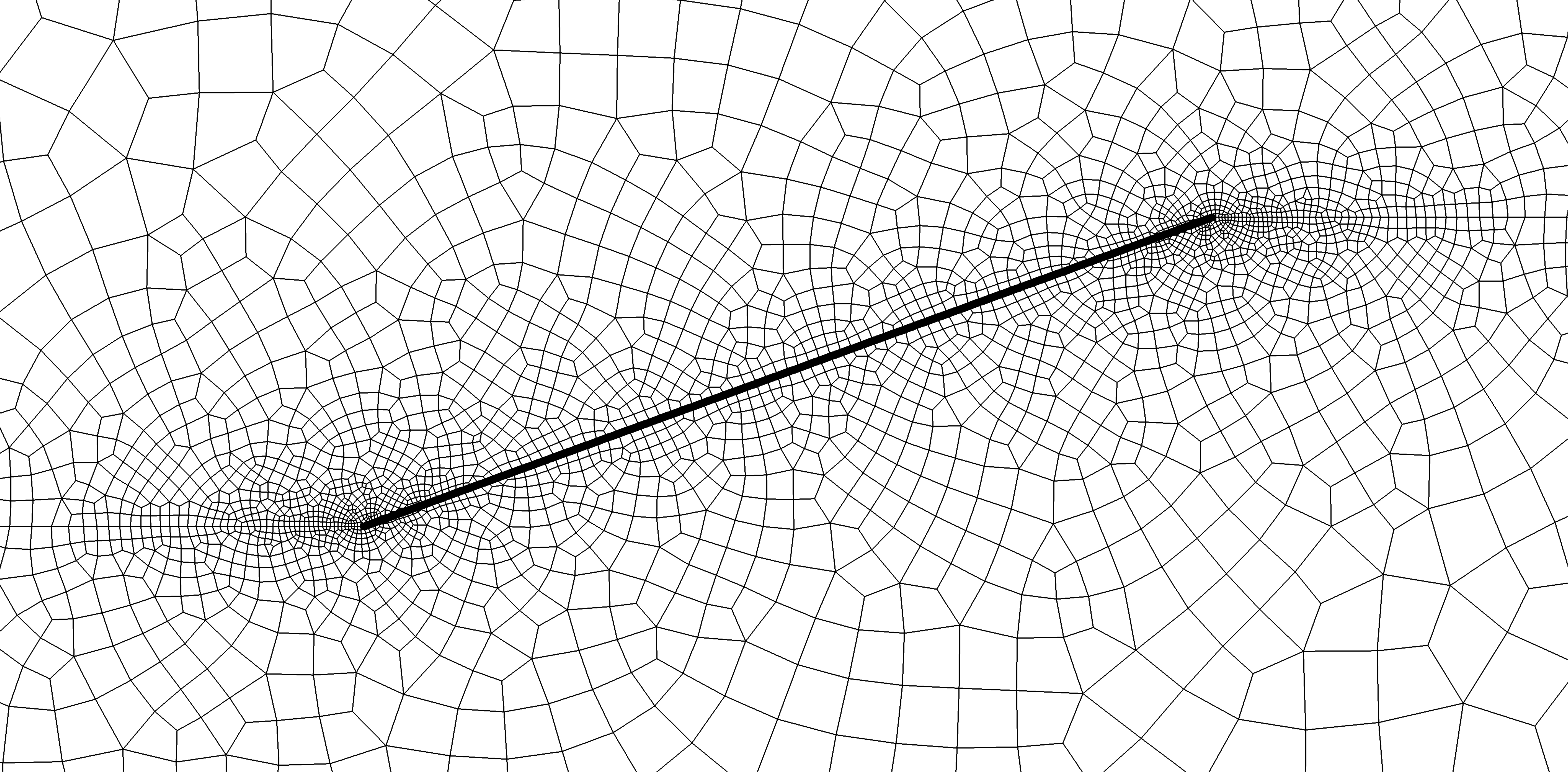}}
    \caption{}  
    \label{fig:phan_grid_zoom}
  \end{subfigure}
  \hfill\null
  \caption{Single crack under compression in infinite domain: (a) sketch of setup; (b) computational mesh; (c) zoom close to the crack. Diamonds (\protect\tikz \protect\draw[black,fill=black,line width=.25ex,rotate=45] (0,0) rectangle (0.71ex,0.71ex);) and squares (\protect\tikz \protect\draw[black,fill=black,line width=.25ex] (0,0) rectangle (1ex,1ex);) indicate $u_x = 0$ and $u_y = 0$ Dirichlet boundary conditions, respectively.}
  \label{fig:phan}
\end{figure}

The second example is a single crack in a 2D infinite domain under a constant uniaxial 
compression. The benchmark geometry is described in detail by \cite{phan2003symmetric} and
reproduced in Fig. \ref{fig:phan_geom}. Being $E$ and $\nu$ the linear elastic parameters,
$\overline{\sigma}$ the compressive stress, $\psi$ the fracture inclination, $2b$ its length and
$\theta$ the friction angle for Coulomb's criterion (with zero cohesion), the analytical
solution provides the normal traction on the fracture and the sliding on it:
\begin{equation}
  t_N = -\overline{\sigma} \sin^2\psi \quad \text{and} \quad
  g_T = \left\|\vec{g}_T\right\|_2 = \frac{4\left(1-\nu^2\right)}{E}\left(\overline{\sigma}
\sin\psi \left( \cos\psi - \sin\psi\tan\theta\right)\right) \sqrt{b^2 -
\left(b-\xi\right)^2},
  \label{eq:phan_sol}
\end{equation}
respectively, where $0 \le \xi \le 2b$ is a local coordinate on the fracture. A plane
strain status is assumed. For the simulation, we used the following values: $E = 25000 \text{ MPa}$ and $\nu = 0.25$, the friction angle is $\theta =
30^\circ$, the fracture is tilted by $\psi = 20^\circ$ and extends for $2b = 2
\text{ m}$, and $|\overline{\sigma} : (\vec{e}_x \otimes \vec{e}_x) | = 100 \text{ MPa}$.

The model is discretized with different resolutions in the $x$-$y$ plane, from 5K to 68K quadrilaterals and from
108 to 432 interface elements, for the coarsest and the finest grid, respectively. The final 3D
domain is obtained by extrusion.
The less refined computational domain is shown in Fig.
\ref{fig:phan_grid}, where a zoom on the mesh around the fracture is also provided. The
boundary conditions for $u_x$ and $u_y$ are set in order to respect the symmetry of the
expected solution (see Fig. \ref{fig:phan_grid}). The two faces parallel to the $x$-$y$
plane are constrained in the $z$ direction, because of the plane strain assumption.

\begin{figure}
  \hfill
  \begin{subfigure}[b]{0.35\textwidth}
    \centerline{\includegraphics[width=\linewidth]{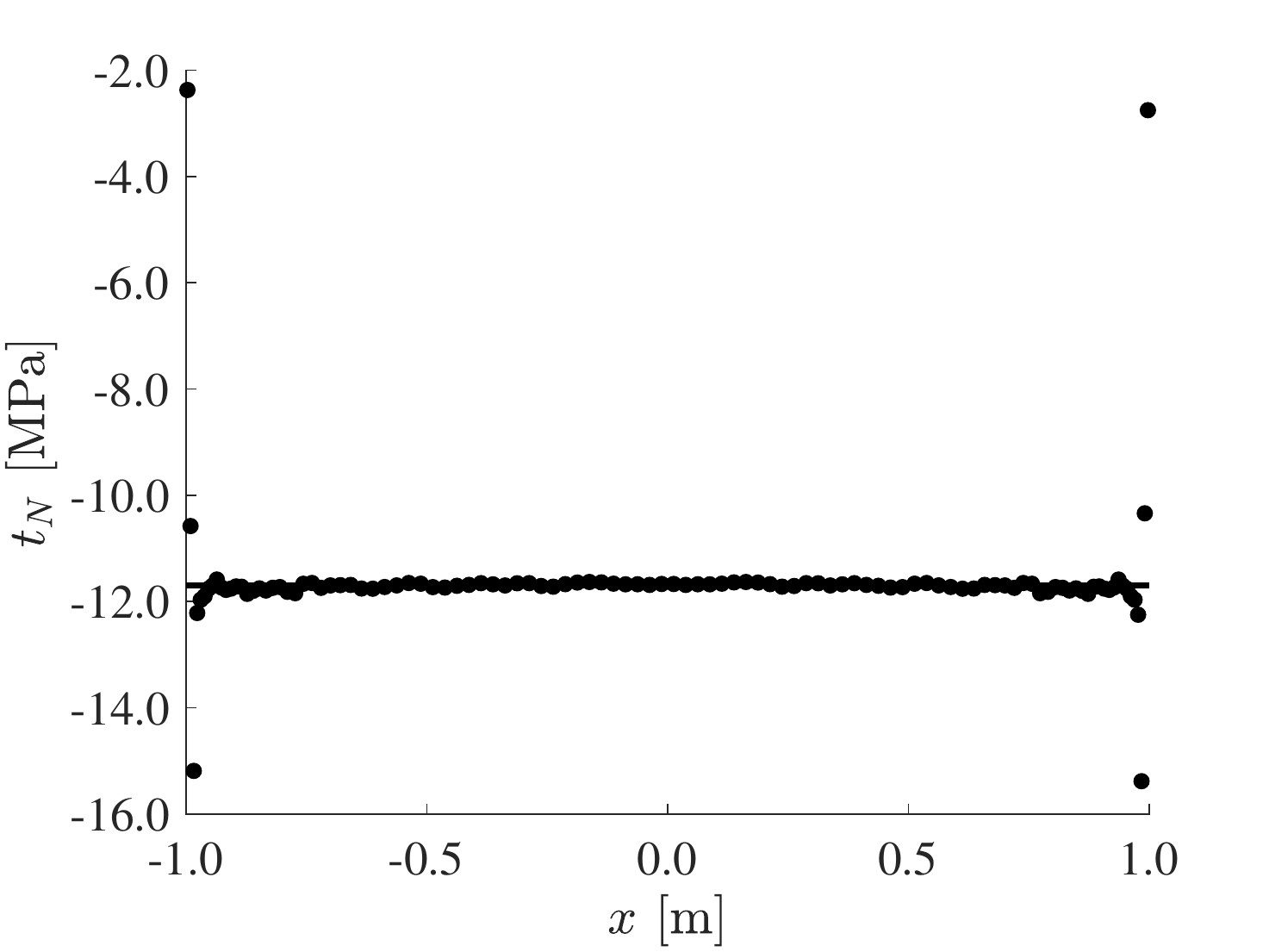}}
    \caption{}  
    \label{fig:phan_sol_tN}
  \end{subfigure}
  \hfill
  \begin{subfigure}[b]{0.35\textwidth}
    \centerline{\includegraphics[width=\linewidth]{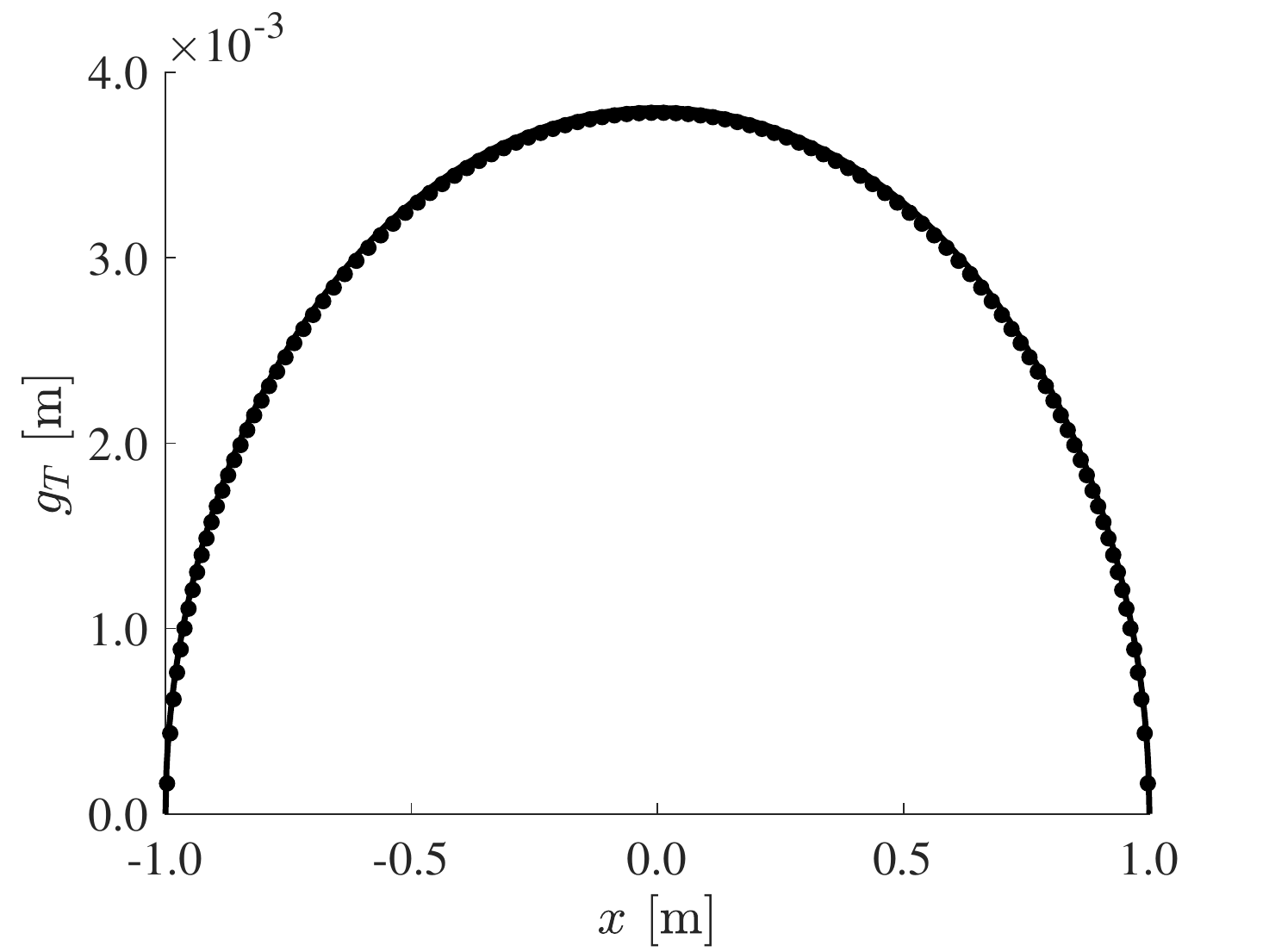}}
    \caption{}  
    \label{fig:phan_sol_gT}
  \end{subfigure}
  \hfill\null
  \caption{Single crack under compression in an infinite domain: comparison between analytical and numerical solutions.}
  \label{fig:phan_sol}
\end{figure}

Comparisons between numerical and analytical solution are provided for the frictional traction (Fig. \ref{fig:phan_sol_tN}) and the relative displacement (Fig. \ref{fig:phan_sol_gT}).
It can be noticed that the computed displacement is in good agreement with the expected one everywhere, while the traction is quite different very close to the fracture tip, where some oscillations are observed.
Refining the mesh, these oscillations are still present, but always closer to the fracture tip.
From a physical viewpoint, these can be explained as a
\textit{locking} phenomenon \cite{borja2008assumed}.
Given the imposed external load, the fracture tries to slip, but the elements close to the tip are not allow for because of the tip itself.
To accommodate the expected behavior, the fracture has to open near the tip, even if there is a compressive normal traction on the crack.
With a less refined mesh, the interface element with a fixed edge on the tip opens, yielding a positive normal traction.
In agreement with these considerations, the rightmost and leftmost traction values are much lower than the average, tending to zero, i.e., the near-tip elements tend to open.

\begin{figure}
  \hfill
  \begin{subfigure}[b]{0.35\textwidth}
    \centerline{\includegraphics[width=\linewidth]{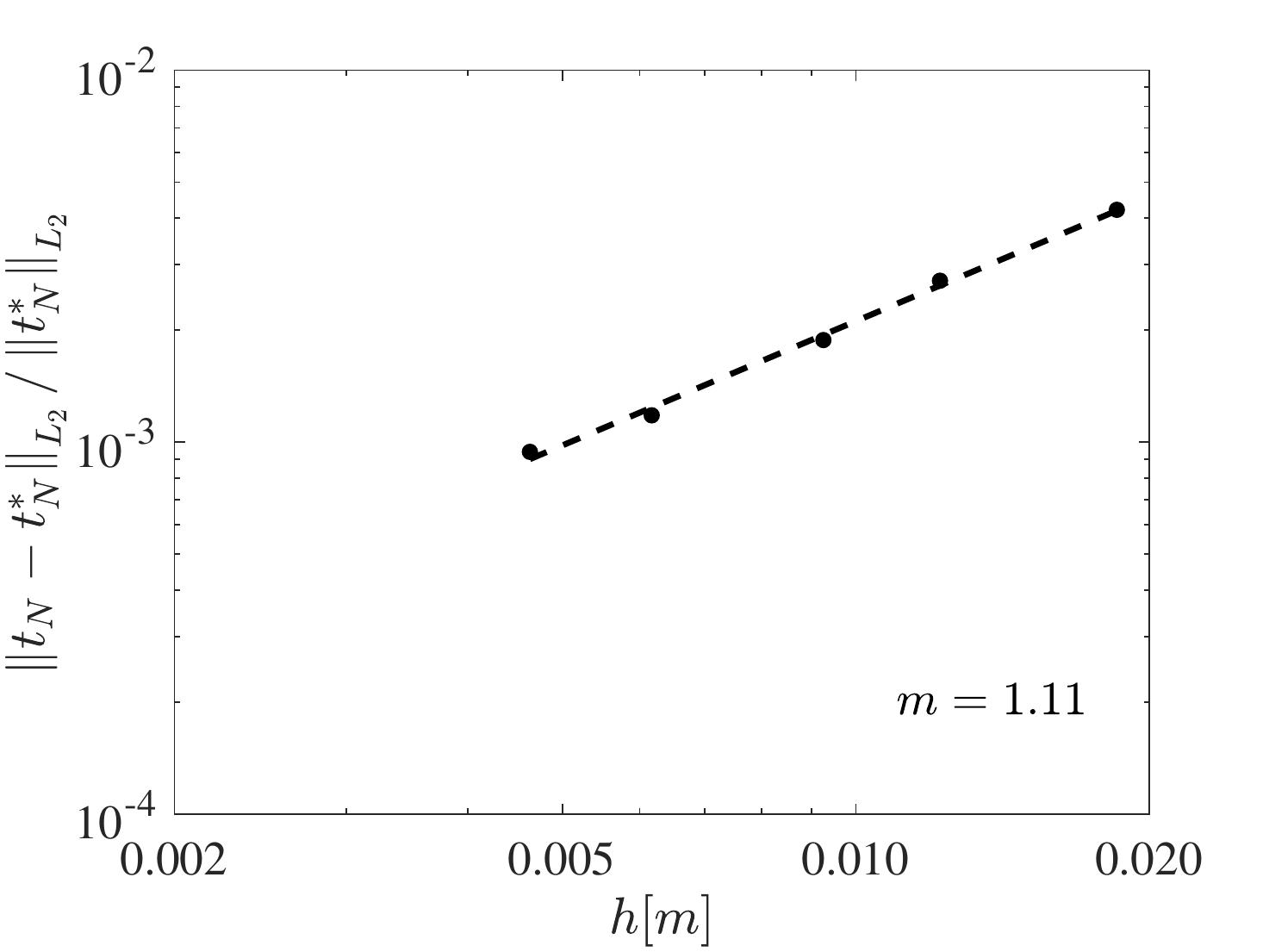}}
    \caption{}  
    \label{fig:convRate_phan_sig}
  \end{subfigure}
  \hfill
  \begin{subfigure}[b]{0.35\textwidth}
    \centerline{\includegraphics[width=\linewidth]{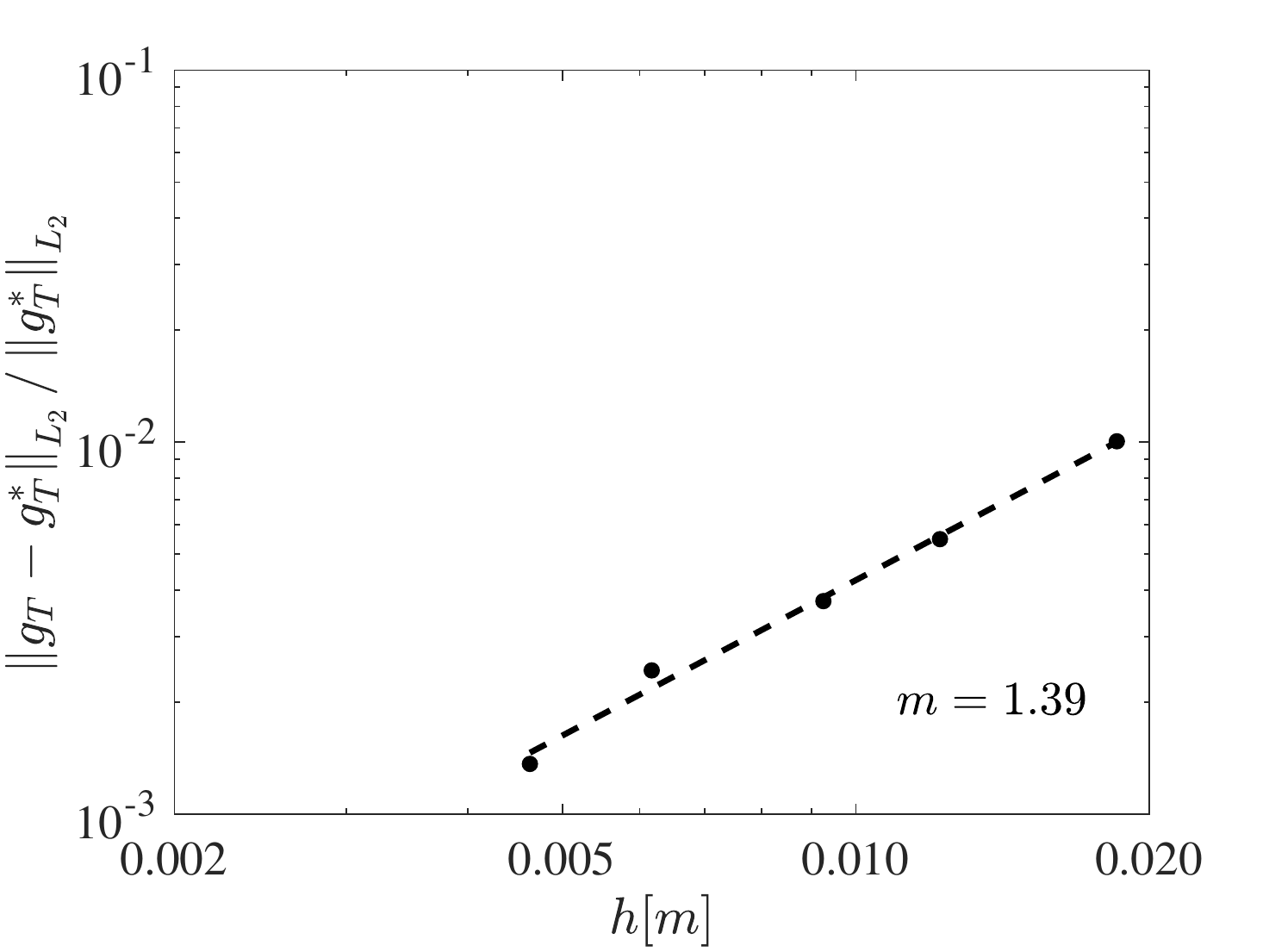}}
    \caption{}  
    \label{fig:convRate_phan_duT}
  \end{subfigure}
  \hfill\null  
  \caption{Single crack under compression in an infinite domain. Convergence rates for: (a) normal traction component, and (b) tangential component of the relative displacement across the crack. The order of convergence $m$ is computed through a least-square linear interpolation.}
  \label{fig:convRate_phan}
\end{figure}

Using the previous example, where an analytical solution for displacement and traction field on the fracture is available, we study the convergence rate, i.e., the error dependence on the mesh size.
First of all, we observe that the oscillations on $t_N$ (see Fig. \ref{fig:phan_sol_tN}) will prevent any traction error norm from convergence, thus, as done in \cite{manzini2008finite}, to compute a meaningful norm we neglect the extreme portions of the domain, i.e., the traction related norms are computed on the central $90 \%$ of the fracture trace.
In Fig. \ref{fig:convRate_phan}, we show the convergence of the two error norms, on the normal Lagrange multiplier (Fig. \ref{fig:convRate_phan_sig}) and on the sliding component of the displacement (Fig. \ref{fig:convRate_phan_duT}).
While the first one is a properly defined norm on the domain $\Gamma_f$, where $\vec{t}$ is defined, the latter one is just the norm of a continuous field projected on a surface, indeed $\vec{u}$ is defined on $\Omega$ but the norm is computed on $\Gamma_f$ only.
The order of convergence $m$ is slightly larger than $1$ for the Lagrange multiplier norm, while it is close to $1.4$ for the projection of the displacement.
There is no well-defined value for the convergence rate of the mixed finite elements space used in this work.
The rate $1.11$ obtained here is in agreement with the theoretical results given in \cite[Section 4]{wohlmuth2011variationally}.

\subsection{Zipper crack problem}
\label{sec:griffith}

\begin{figure}
  \hfill
  \begin{subfigure}[b]{0.45\textwidth}
    \centerline{\includegraphics[width=\linewidth]{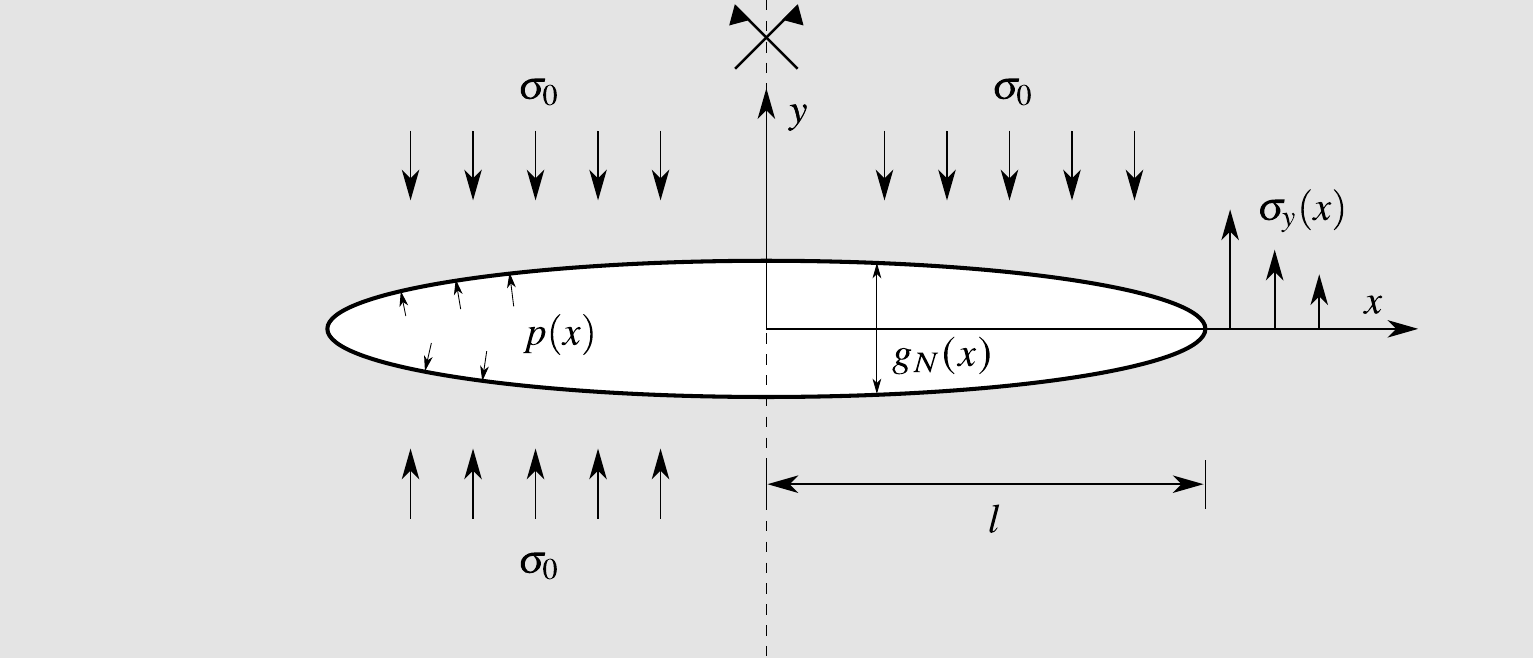}}
    \caption{}  
    \label{fig:griffith_sketch}
  \end{subfigure}
  \hfill
  \begin{subfigure}[b]{0.35\textwidth}
    \centerline{\includegraphics[width=\linewidth]{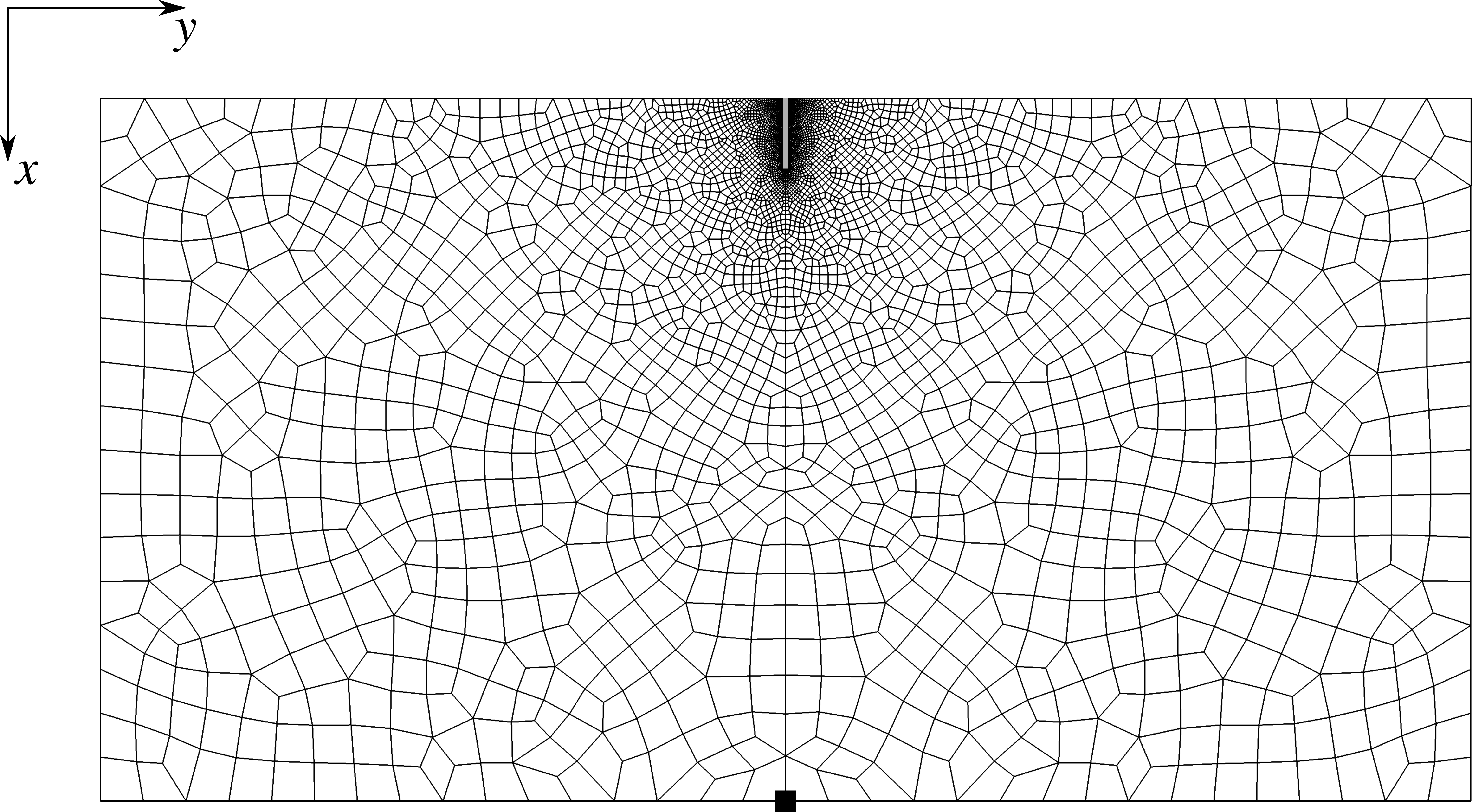}}
    \caption{}  
    \label{fig:griffith_mesh}
  \end{subfigure}
  \hfill\null
  \caption{Zipper crack problem: (a) sketch of the setup; (b) computational domain with the fracture highlighted as a gray segment.}
\end{figure}

To complete the validation of the model, we use the line crack problem as described in \cite{valko1995hydraulic}.
The zipper crack case is also known as Griffith problem \cite{griffith1920phenomena} and consists of an infinite plane ($x$-$y$ plane) with a single linear crack of length $2l$ in the $x$-direction.
Inside the fracture, there is a fluid with a given pressure $p(x)$.
The assumption is of plane strain.
Fig. \ref{fig:griffith_sketch} represents a sketch of the setup. The analytical solution provide the opening of the fracture for every location $x$ and the stress in the continuous medium on the $x$ direction ($y=0$), for $x > l$.
The far-field stress orthogonal to the fracture is $\sigma_0$.
We impose a constant pressure $p_0$ only on one part of the fracture, in such a way that the tip closes smoothly.
The pressurized length ends at $x_0$, with:
\begin{equation}
  x_0 = l \sin\frac{\pi \sigma_0}{p_0}.
\end{equation}
Introducing $q_1 = \sqrt{l^2-x_0^2}$, $q_2 = \sqrt{l^2-x^2}$ and $q_3 = \sqrt{\left|x_0^2-x^2\right|}$, the analytical solutions for opening and stress in $y$ direction are \cite{valko1995hydraulic}:
\begin{align}
  g_N(x) &= \frac{2\left(1-\nu^2\right)}{\pi E} p_0 \left(4 x_0 \log\frac{q_1+q_2}{q_3} +
x\log\frac{l^2x_0^2 - 2q_1q_2x_0x + l^2x^2 - 2x_0^2x^2}{l^2x_0^2 + 2q_1q_2x_0x + l^2x^2 -
2x_0^2x^2}\right),
    && 0 \le x \le l,\\
  \sigma_y(x,0) &= -p_0 + \frac{2}{\pi}\left(p_0 \arctan \frac{x q_1}{x_0 q_3}\right),
    && x \ge l.
\end{align}
In Fig. \ref{fig:griffith_mesh}, we show the computational domain used to reproduce this solution.
The 3D domain, whose size is $150 \times 300 \times 0.3 \text{ m}$, is discretized with 17772 nodes, 13050 hexahedra and 228 interface finite elements.
The objective is to simulate a fracture length of $l = 10 \text{ m}$, but we discretize a longer fracture, $l_1 = 15 \text{ m}$, to verify if the remaining length $l_1 - l = 5 \text{ m}$ remains closed.
We need a large domain to attenuate the boundary effect as the analytical solution is obtained for an infinite domain.
It is $10 \times 20$ times the actual fracture length.
Referring to Fig. \ref{fig:griffith_mesh}, the only constrained boundary is the one intersecting the fracture, where symmetric conditions are imposed.
The only fixed point is the one highlighted in the figure.
The 3D domain is obtained extruding a 2D domain, and the faces parallel to the $x$-$y$ plane are $z$-constrained, to fulfill the conditions of the assumed plane strain state.
Regarding material properties, we have $E = 25 \text{ GPa}$ and $\nu = 0.25$.
Fluid pressure and far-field stress are $p_0 = 15 \text{ MPa}$ and $\sigma_0 = 10 \text{ MPa}$.
With the chosen set of parameters, we have $\lim_{x \rightarrow l} \sigma_y(x,0) = -15 + \frac{30}{\pi}\arctan\frac{2}{\sqrt{3}} \approx -6.816 \text{ MPa}$, thus there is a discontinuity in the stress field at the fracture tip.

\begin{figure}
  \hfill
  \begin{subfigure}[b]{0.32\textwidth}
    \centerline{\includegraphics[width=\linewidth]{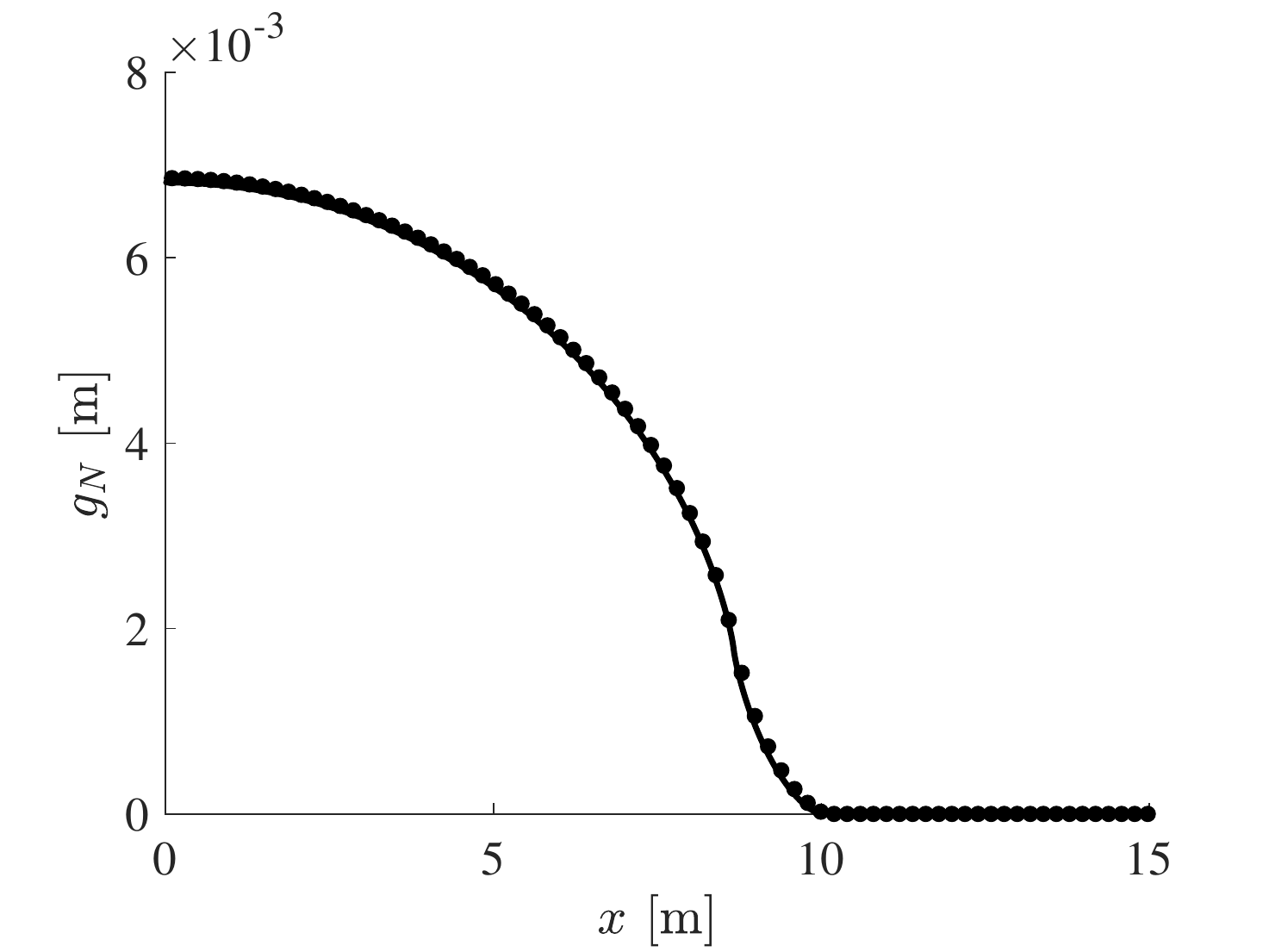}}
    \caption{}  
    \label{fig:griffith_sol_duN}  
  \end{subfigure}
  \hfill
  \begin{subfigure}[b]{0.32\textwidth}
    \centerline{\includegraphics[width=\linewidth]{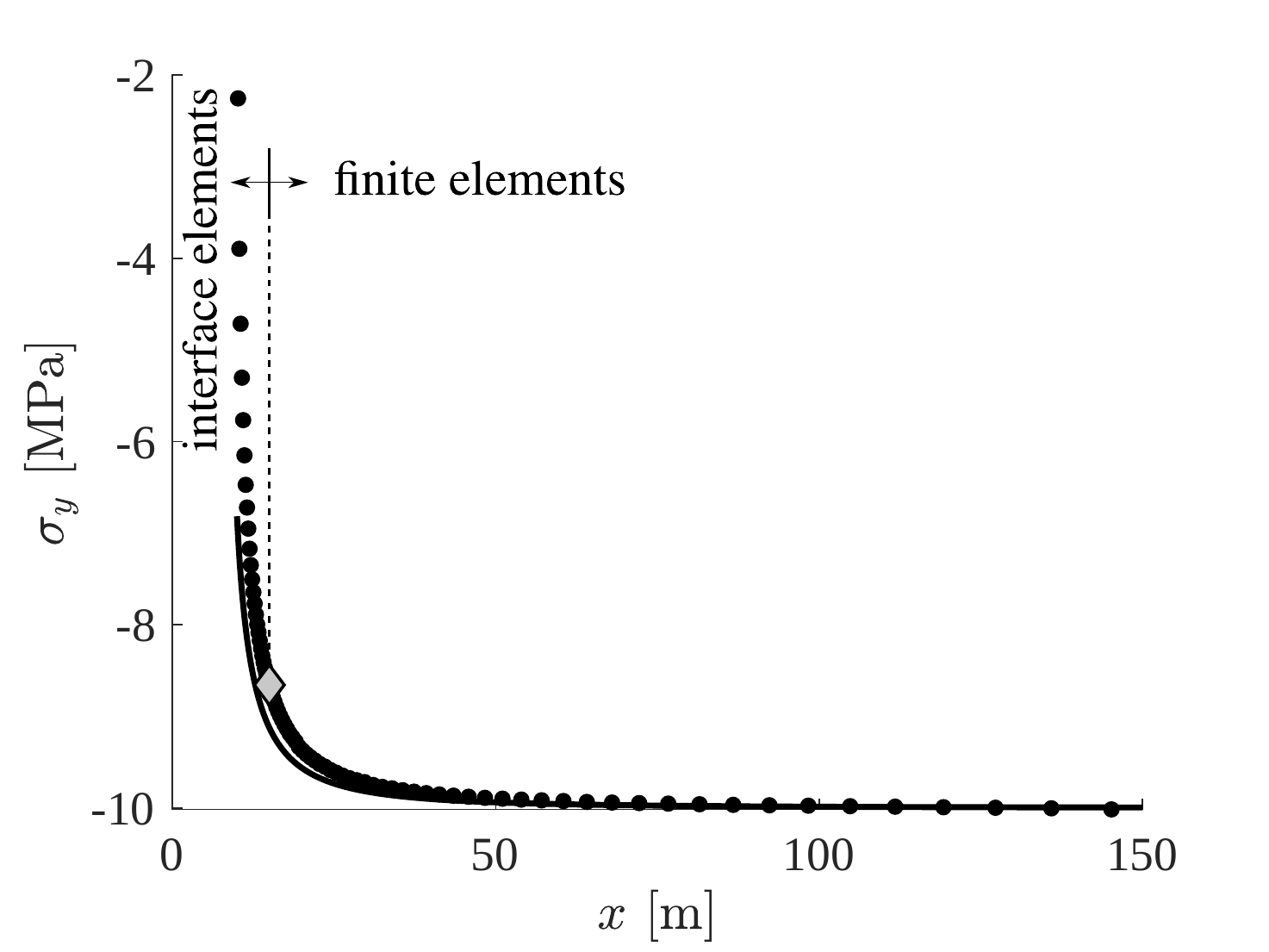}}
    \caption{}  
    \label{fig:griffith_sol_sigmayy}  
  \end{subfigure}
  \hfill
  \begin{subfigure}[b]{0.32\textwidth}
    \centerline{\includegraphics[width=\linewidth]{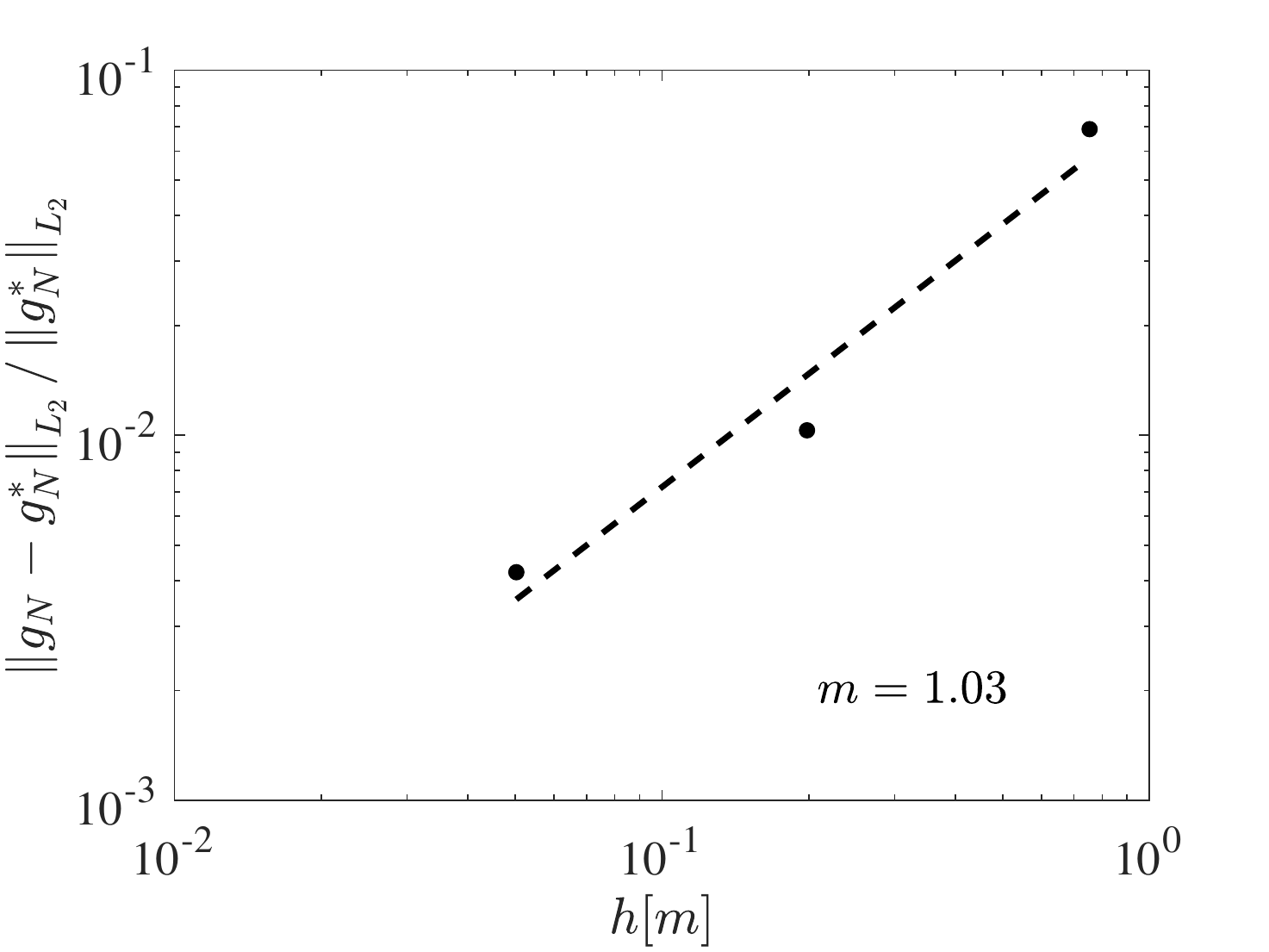}}
    \caption{}
    \label{fig:griffith_conv}  
  \end{subfigure}
  \hfill\null
  \caption{Zipper crack problem: (a) comparison between numerical and analytical solution for fracture opening; (b) comparison between numerical and analytical solution for $\sigma_y$-stress int the continuous domain; (c) convergence rate for fracture opening with the the order of convergence $m$ computed through a least-square linear interpolation. In (b), the gray diamond (\protect\tikz \protect\draw[black,fill=gray!22,line width=.25ex,rotate=45] (0,0) rectangle (0.71ex,0.71ex);) indicate where the discretization with interface elements ends and classical finite elements are used.}
  \label{fig:griffith_sol}
\end{figure}

In Figs. \ref{fig:griffith_sol_duN}-\ref{fig:griffith_sol_sigmayy}, we report the comparison between the outcomes of the numerical model and the analytical solutions.
It can be notice, for $x > l$, the fracture is closed and the transition is smooth, as predicted by the theoretical solution.
The stress behavior is very similar to the analytical one, except very close to the tip, where the numerical model shows a smaller jump between the fracture and the continuous material.
Nevertheless, we have a smooth transition between the solutions computed with two different discretizations---on a portion of the closed length interface elements are used, while on the other classical finite elements are employed.
Overall, we can see a good agreement, with an integral relative error of $1.0 \%$ for the displacements and $2.5 \%$ for the stress. For the error evaluation related to the stress, we neglected the near-tip portion of the domain, i.e., we considered only $x > 5/4 l = 12.5 \text{ m}$.

As a concluding remark for this analytical benchmark, Fig. \ref{fig:griffith_conv} shows the convergence rate for the $L_2$-norm of the error on the fracture aperture $g_N(x)$ only, as the stress field is unbounded close to the tip, so unsuited for such a test.
As expected, the order of convergence is around $1$ because of the discontinuity in the stress field, indeed, the asymptotic behavior is almost always lost whenever there are singularities in the solution \cite{zienkiewicz2000finite}.
The coarser mesh has $1K$ finite elements and $60$ interface finite elements, while the finest has $194K$ and $900$ finite and interface finite elements, respectively.

\section{Analytical benchmarks for contact mechanics with fluid flow}
\label{sec:anal_flow}
To verify the numerical model for the full IBVP \eqref{eq:IBVP}, we use two classic analytical solutions: (i) the Kristianovic-Geerstma-deKlerk (KGD) problem  \cite{khristianovic1955formation,geertsma1969rapid,valko1995hydraulic,adachi2002self,detournay2004propagation,rahman2010review,
Set_etal17} and (ii) the penny-shaped crack problem \cite{selvadurai1980penny,
valko1995hydraulic,savitski2002propagation,detournay2004propagation,Set_etal17,abe1976growth}.
The first test case is mainly a 2D problem, while the latter is a real 3D case.

\subsection{KGD problem}
\label{sec:KGD}
\begin{figure}
  \hfill
  \begin{subfigure}[b]{0.45\textwidth}
    \centerline{\includegraphics[width=\linewidth]{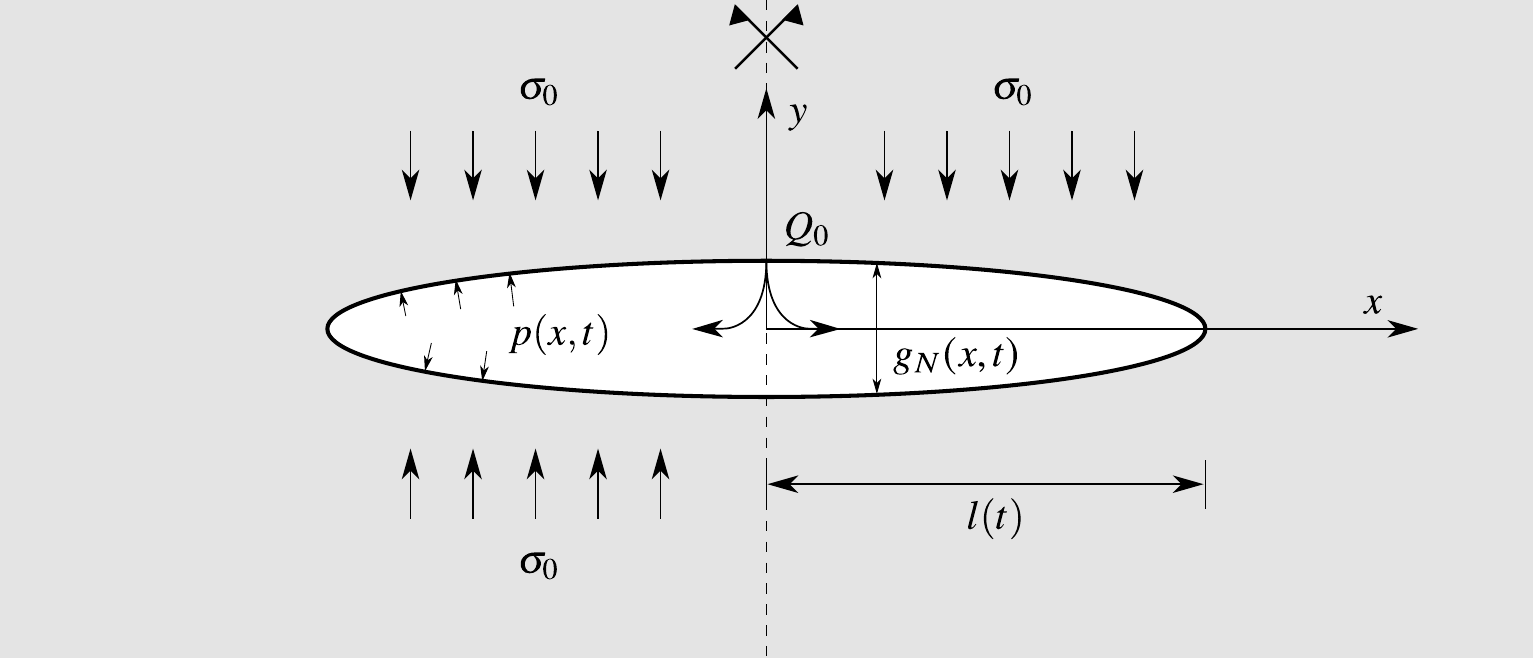}}
    \caption{}  
    \label{fig:KGD_sketch}
  \end{subfigure}
  \hfill
  \begin{subfigure}[b]{0.35\textwidth}
    \centerline{\includegraphics[width=\linewidth]{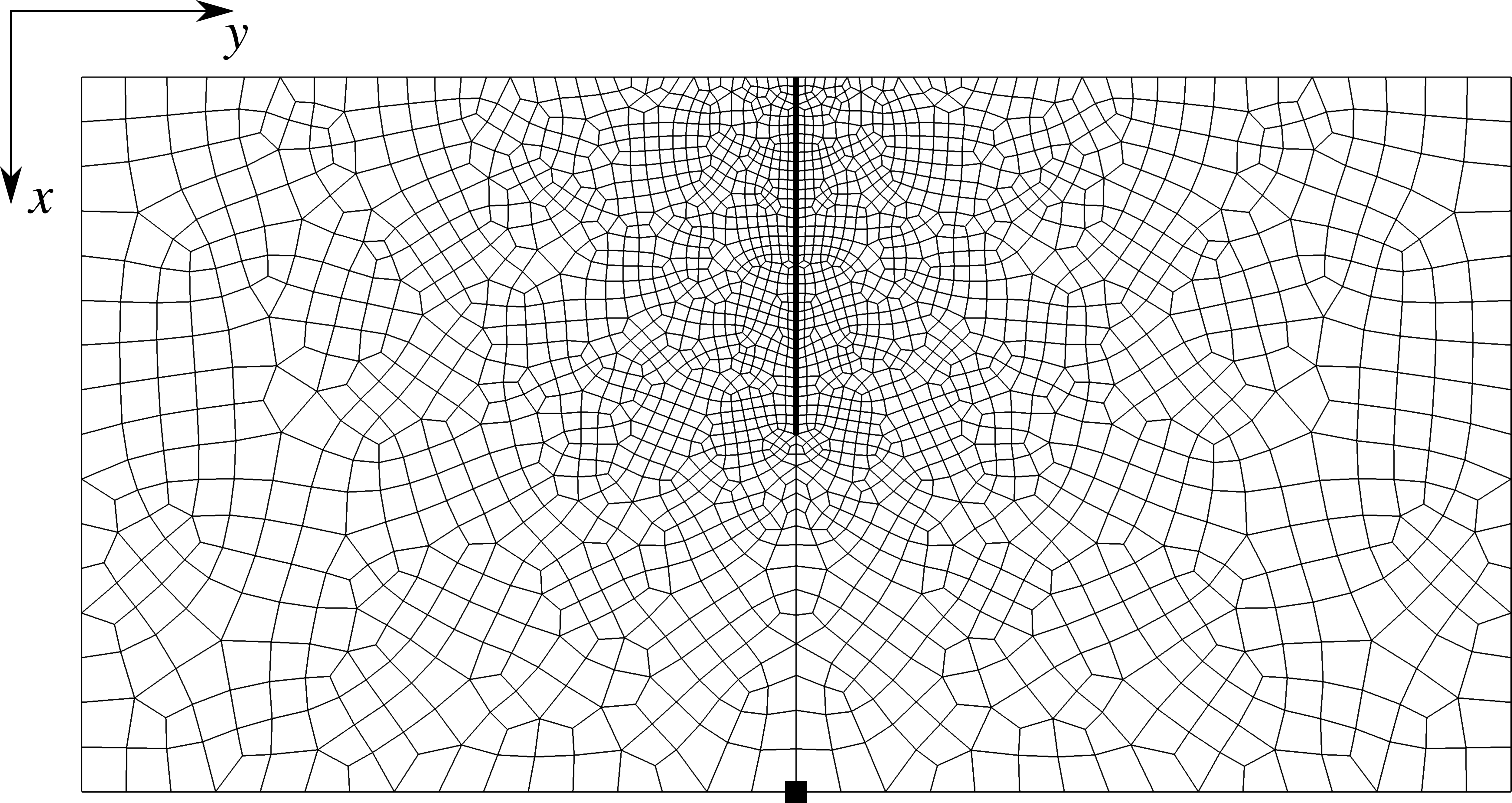}}
    \caption{}  
    \label{fig:KGD_mesh}
  \end{subfigure}
  \hfill\null
  \caption{KGD problem: (a) sketch of the setup; (b) computational domain with the fracture highlighted as a gray segment. The square (\protect\tikz \protect\draw[black,fill=black,line width=.25ex] (0,0) rectangle (1ex,1ex);) represents the projection on the $x$-$y$ plane of the $y$-constrained nodes.}
\end{figure}

We consider a 2D hydraulic fracture propagation assuming plane strain conditions. The
medium is isotropic, homogeneous, impermeable,
and is fully described using a linear elastic model. An incompressible Newtonian fluid is
injected
from a fixed point, at a constant rate $Q_0$. The fracture propagates in the direction
that is orthogonal to the maximum principal direction of the stress tensor in the
surrounding medium. In Fig.\ \ref{fig:KGD_sketch}, we represent the set-up of the problem
and introduce the quantities of interest:
\begin{itemize}
  \item $g_N(x,t)$: the fracture opening for any time and position $x \le l$;
  \item $p(x,t)$: the net fluid pressure inside the fracture for any time and position
    $x \le l$;
  \item $l(t)$: the fracture half-length.
\end{itemize}

The analytical solution is provided in terms of $g_N(x,t)$, $p(x,t)$ and $l(t)$, but the
complete expressions require the definition of some dimensionless quantities
\cite{detournay2004propagation}: (i) an opening $\Omega$, (ii) a net pressure $\Pi$ and
(iii) a fracture length $\gamma$. Given these dimensionless variables, the quantities of
interest become:
\begin{subequations}
\begin{align}
  l(t) &= \gamma L(t), &
  g_N(x,t) &= \varepsilon(t) L(t) \Omega(\xi), &
  p(x,t) &= \varepsilon(t) E^\prime \Pi(\xi), &
\end{align}
\label{eq:KGD_sol1}\null
\end{subequations}
where $\xi = x/l$ is the similarity variable, i.e., a dimensionless fracture coordinate.
For the zero toughness case, we have \cite{adachi2002self}:
\begin{subequations}
\begin{align}
  \varepsilon(t) &= \left(\frac{\mu^\prime}{E^\prime t}\right)^{\frac{1}{3}}, &
  L(t) &= Q_0^{1/2} \left(\frac{E^\prime}{\mu^\prime}\right)^{\frac{1}{6}} t^{2/3}, &
  \gamma &= \frac{1}{\left(2 \int_0^1 \overline{\Omega} \, \mathrm{d} \xi\right)^{1/2}}.
\end{align}
\label{eq:KGD_sol2}\null
\end{subequations}
In \eqref{eq:KGD_sol2}, we introduced $\overline{\Omega} = \frac{\Omega}{\gamma}$,
$\mu^\prime = 12 \mu_l$, with $\mu_l$ the fluid viscosity, and $E^\prime =
\frac{E}{1-\nu^2}$, i.e., the plane strain modulus. The two functions
$\overline{\Omega}(\xi)$ and $\Pi(\xi)$, called self-similar fracture opening and
self-similar net fluid pressure, respectively, are approximated through polynomial series,
based on the Gegenbauer polynomials \cite{abramowitz1988handbook} for the former one,
while the latter uses Euler's beta functions and Gauss's hyper-geometric functions
\cite{abramowitz1988handbook}. The full expression of these two dimensionless functions,
with the numerical coefficients of the series expansion, are provided in
\cite{adachi2002self}.
This analytical solution is based on zero-toughness and zero-lag assumptions. The
first hypothesis implies that the energy dissipated by the fracture propagation is negligible
compared to the energy dissipated in the fluid by viscous
flow. In our current implementation, the fracture surface is pre-defined and no energy is dissipated in fracturing. Regarding the second assumption, for real fractures a fluid lag or
non-wetted zone may appear. From physical considerations, the fluid pressure at the
tip has to be finite, as well as the stress in the rock surrounding the fracture tip. In the regime of interest, however, the influence of this region on the global pressure and displacement
solution is confined in a small region near the tip (see numerical
results \cite{desroches1994crack,carbonell1999comparison}). Given these considerations,
the analytical solution used herein is a good approximation of the main physical
behavior.
We emphasize the fact that $\lim_{\xi \rightarrow 1} \Pi = -\infty$, thus, the
analytical solution predicts an infinite pressure at the fracture tip. This nonphysical
result is due to the assumption that the fluid reaches the tip and fill every empty space
at the same density, without allowing for cavitation. Note that both the fluid density and the initial
stress regime do not affect any quantity of interest.

In Fig. \ref{fig:KGD_mesh}, we show the mesh used for the simulation. It is composed by
10209 nodes, 6600 hexahedra and 76 interface elements. The domain is $300 \times 600
\times 8 \text{ m}$, with a fracture of $150 \text{ m}$ and an average element size of
$h_x \approx 4 \text{ m}$. We highlight that currently our model does not handle 
fracture propagation, but we know a priori the fracture trajectory and can pre-discretize a surface of sufficient length.  In some sense this is then a fracture ``reactivation'' problem. To minimize the boundary effects, we imposed the
symmetric condition on the face parallel to the $y$ axis and intersecting the fracture and
$z$-constrained the two surfaces parallel to the $x$-$y$ plane, to reflect the plane
strain assumption. The 2D node highlighted in Fig. \ref{fig:KGD_mesh}, that corresponds to
a column of nodes in 3D, is the only one that is $y$-constrained. The material parameters are $E = 30
\text{ GPa}$ and $\nu = 0.25$. The fracture frictional behavior is governed by Coulomb's
criterion, characterized by $\theta = 30^\circ$ and zero cohesion. The fluid viscosity is
$\mu_l = 10^{-9} \text{ MPa}\cdot\text{s}$. The injection rate and the confining stress
are $Q_0 = 2 \cdot 10^{-3} \text{ (m}^3/\text{s)/m}$ and $\sigma_0 = 10 \text{ MPa}$,
respectively. Finally, according to \cite{kamenov2013laboratory}, the conductivity initial
value (see Eq. \eqref{eq:cond_def}) is $C_{f,0} = 10 \text{ mD}\cdot\text{m} = 9.87 \cdot
10^{-15} \text{ m}^2\cdot\text{m}$. The simulated time is $100 \text{ s}$, with $\dt =
0.01 \text{ s}$ from the beginning to $t = 0.2 \text{ s}$, then $\dt = 0.1 \text{ s}$ up
to $t = 2 \text{ s}$ and finally $\dt = 1 \text{ s}$ up to the end of the simulation.
We highlight that for the current flow rate, the Reynolds number is about $2000$, thus the
regime is laminar and the cubic law assumption is reasonable.

\begin{figure}
\hfill
\begin{subfigure}{0.24\linewidth}
  \centering
  \includegraphics[width=\linewidth]{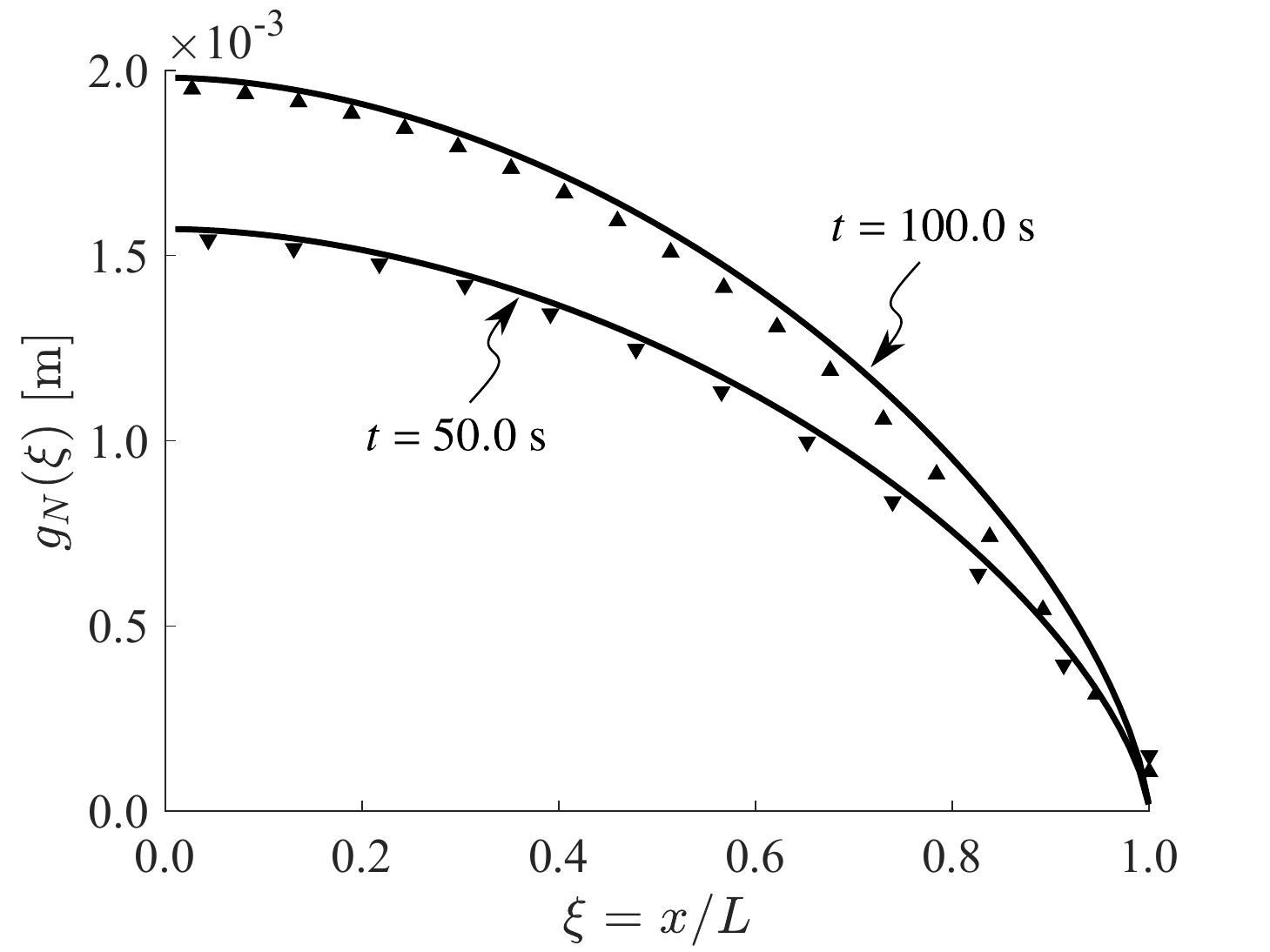}
  \caption{}
  \label{fig:KGD_opening}
\end{subfigure}
\begin{subfigure}{0.24\linewidth}
  \centering
  \includegraphics[width=\linewidth]{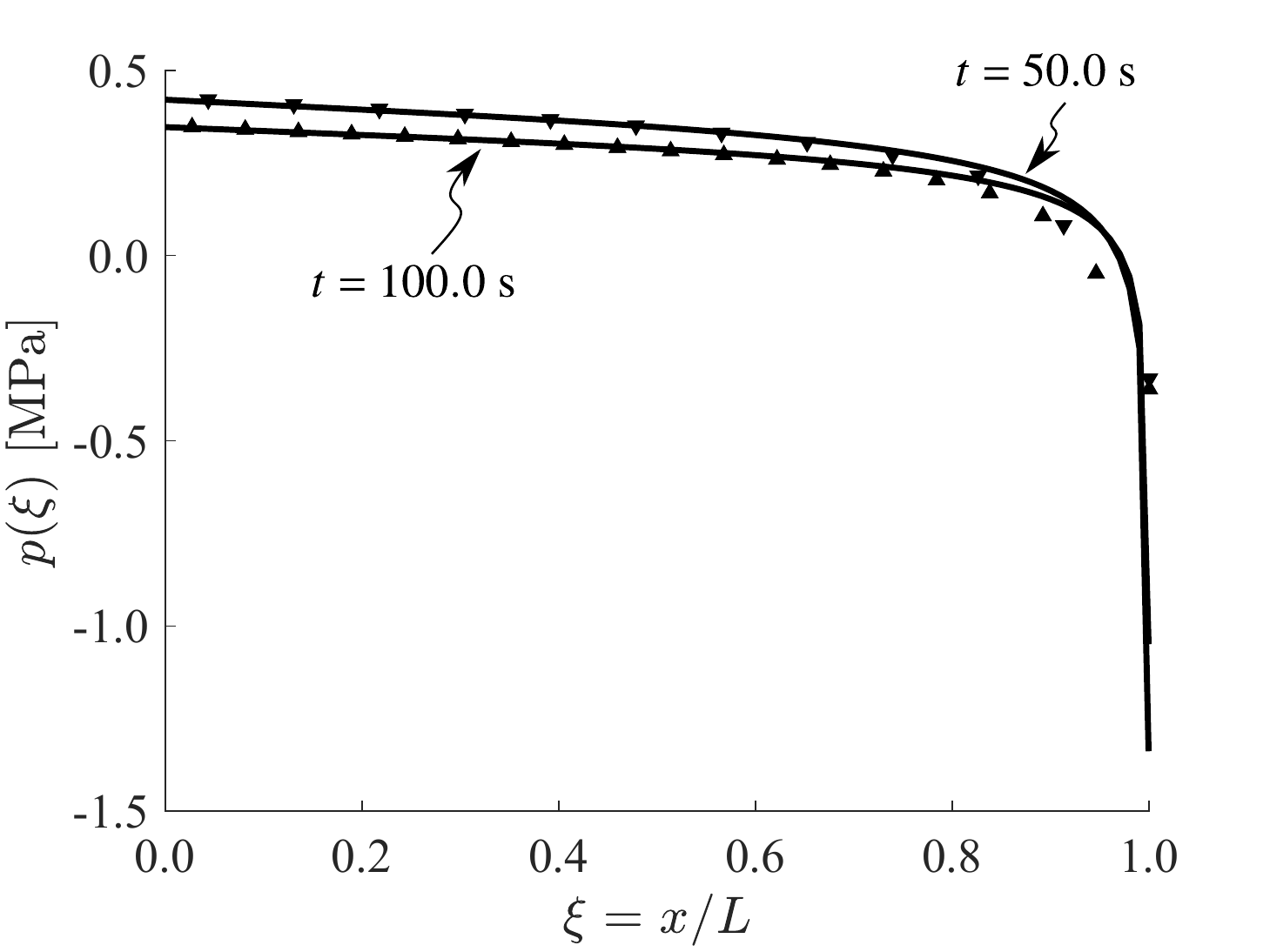}
  \caption{}
  \label{fig:KGD_pressure}
\end{subfigure}
\begin{subfigure}{0.24\linewidth}
  \centering
  \includegraphics[width=\linewidth]{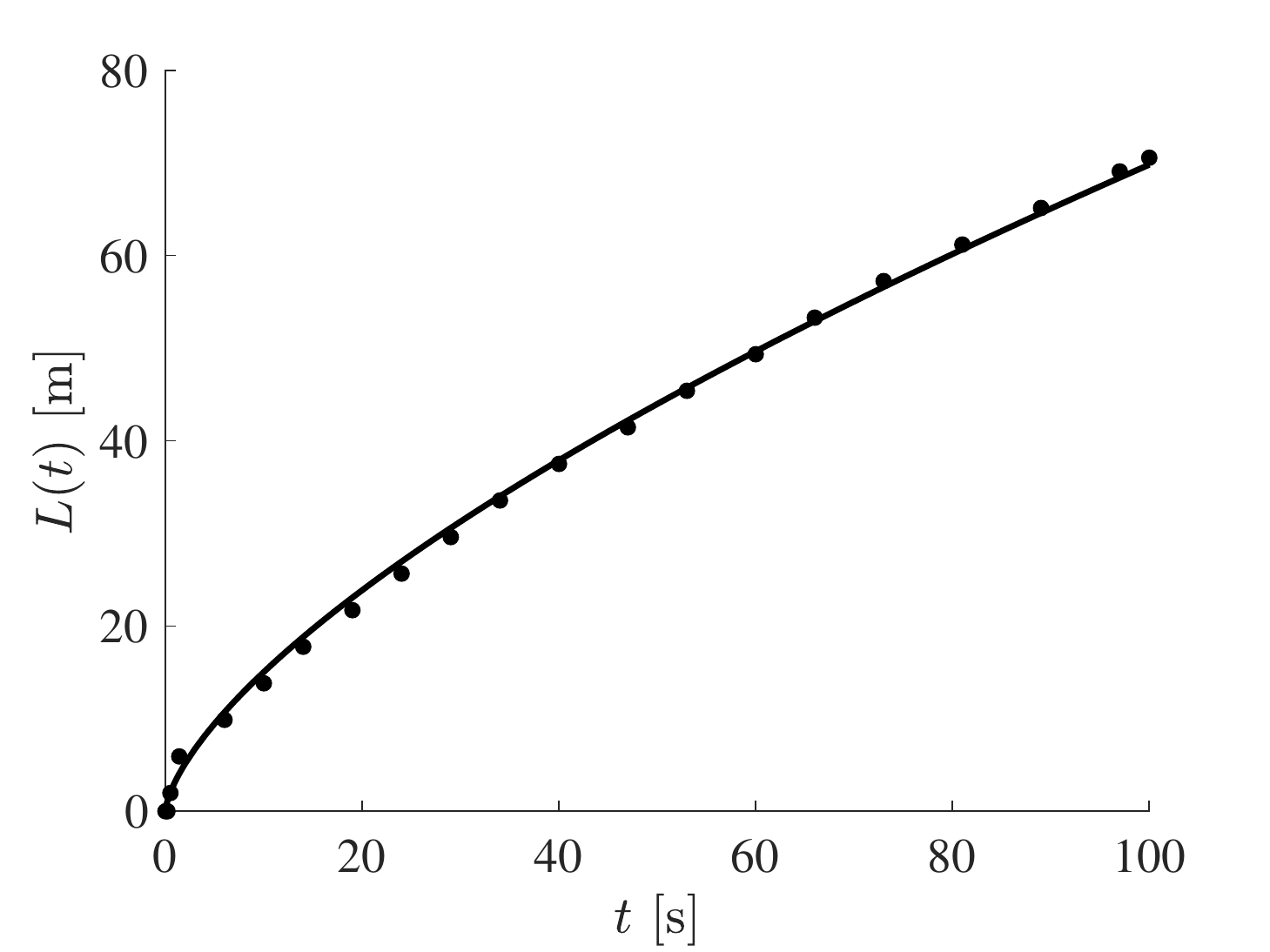}
  \caption{}
  \label{fig:KGD_length}
\end{subfigure}
\begin{subfigure}{0.24\linewidth}
  \centering
  \includegraphics[width=\linewidth]{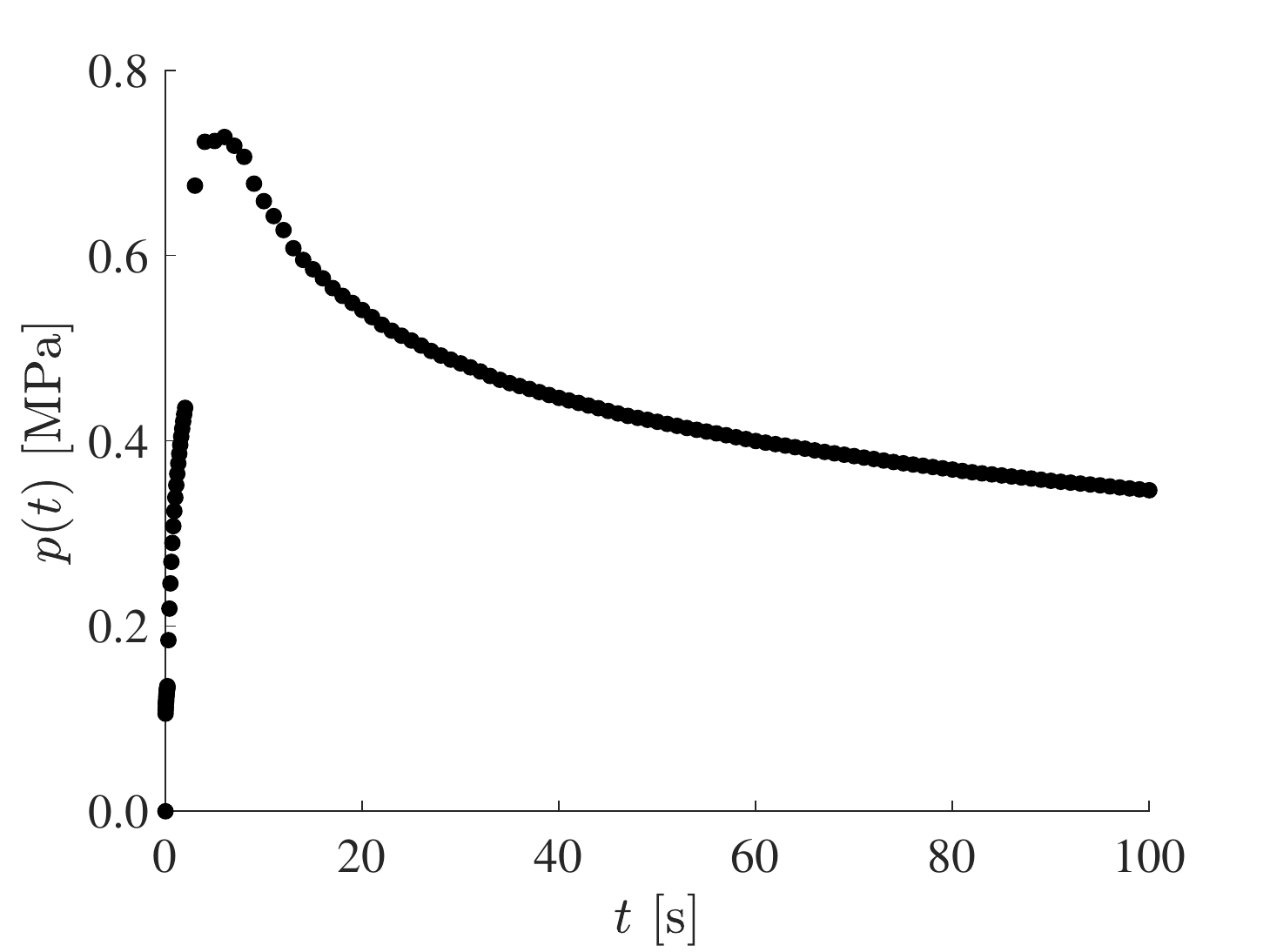}
  \caption{}
  \label{fig:KGD_pressure_x0}
\end{subfigure}
\hfill\null
  \caption{Results for the KGD simulation. From the left to the right: fracture opening,
    fluid pressure, fracture length and pressure behavior at the injection location.}
  \label{fig:KGD_res}
\end{figure}

Figs. \ref{fig:KGD_opening}-\ref{fig:KGD_pressure} show the outcomes of the model in terms
of opening and pressure for two different time steps, i.e., at half ($t = 50 \text{ s}$)
and at the end of the simulation ($t = 100 \text{ s}$). The continuous line is the
analytical solution. Overall, there is a good agreement, for both the aperture (with an
integral relative error of $4.4 \%$ for $t = 100 \text{ s}$) and the pressure (with an
integral relative error of $3.8 \%$ at the same time step). At time $t = 100 \text{ s}$,
the pressure is slightly different close to the tip, where the theoretical behavior tends
to $-\infty$. We emphasize that the analytical solution for the pressure is constrained in
an integral sense by the propagation criterion, being
\begin{equation}
  \int_0^1 \frac{\Pi}{\sqrt{1-\xi^2}} \, \mathrm{d} \xi = 0.
\end{equation}
In our model, we use a Dirichlet boundary condition on the pressure value at the fracture
end, that is not the fracture ``tip'', where $p = 0$ is imposed. To compare the two behaviors,
our outcome in terms of pressure is shifted by a constant value. In Fig.
\ref{fig:KGD_length}, we report the fracture length as function of time. The model is able
to predict quite accurately the fracture length, with an average relative error of $2.1
\%$. Finally, in Fig. \ref{fig:KGD_pressure_x0}, the pressure at the injection location
is shown for all time steps. Because of the different strategy used to impose the pressure
boundary condition, we cannot directly compare this profile with the analytical solution. We
limit ourselves to a visual comparison with similar works, e.g.
\cite{carrier2012numerical,Set_etal17,vahab2017numerical}, observing profiles that are in excellent
agreement.

\subsection{Penny-shaped crack}
\begin{figure}
  \hfill
  \begin{subfigure}[b]{0.45\textwidth}
    \centerline{\includegraphics[width=\linewidth]{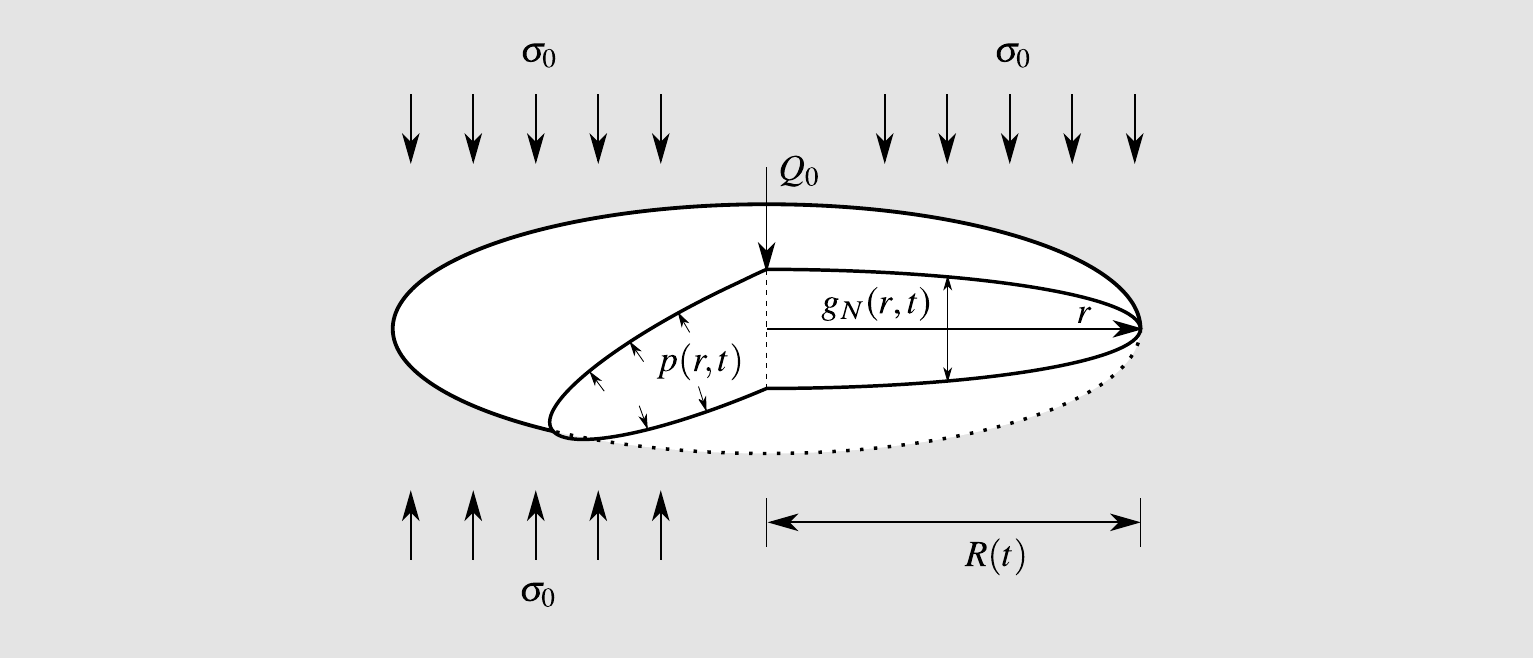}}
    \caption{}  
    \label{fig:penny_sketch}
  \end{subfigure}
  \hfill
  \begin{subfigure}[b]{0.22\textwidth}
    \centerline{\includegraphics[width=\linewidth]{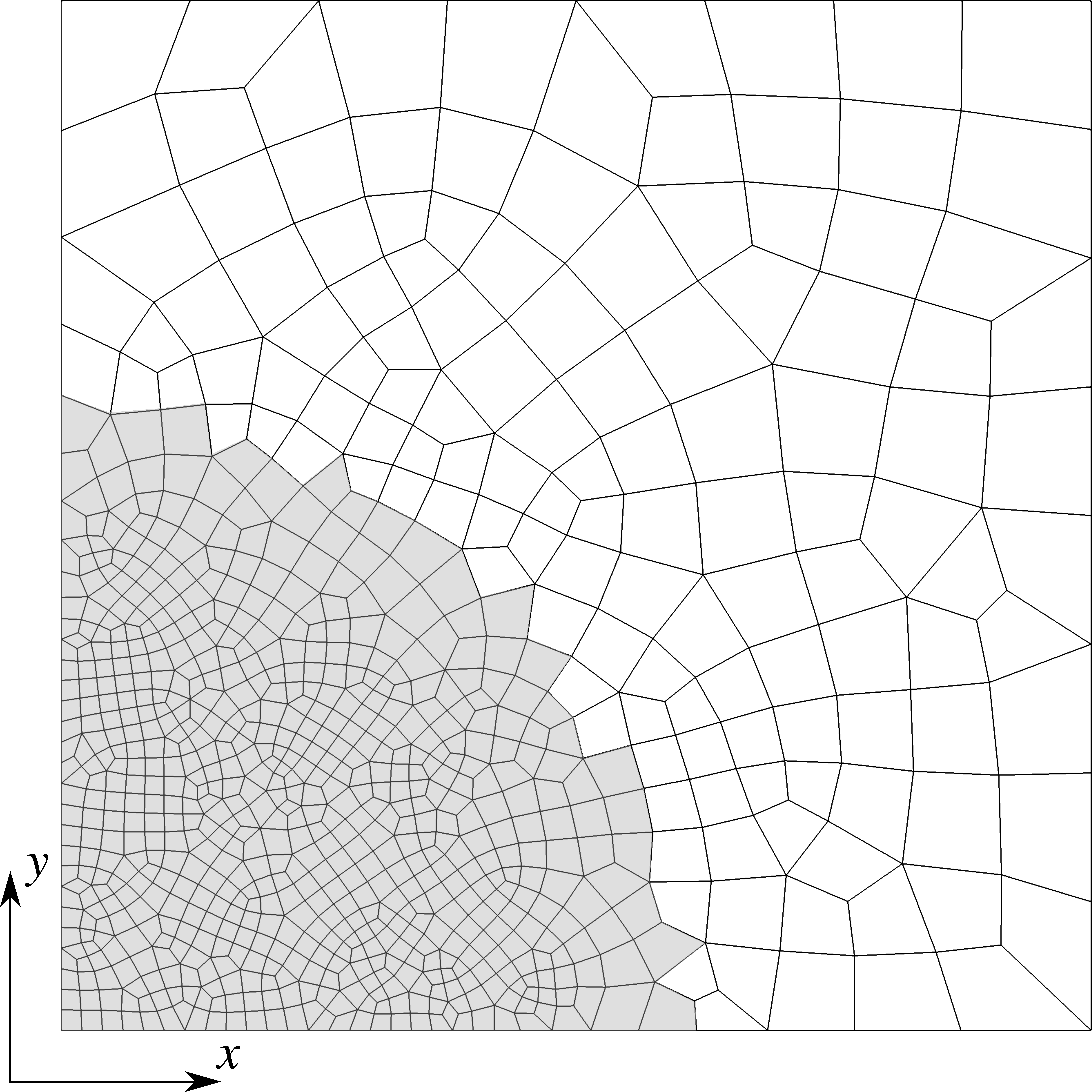}}
    \caption{}  
    \label{fig:penny_mesh_2d}
  \end{subfigure}
  \hfill
  \begin{subfigure}[b]{0.22\textwidth}
    \centerline{\includegraphics[width=\linewidth]{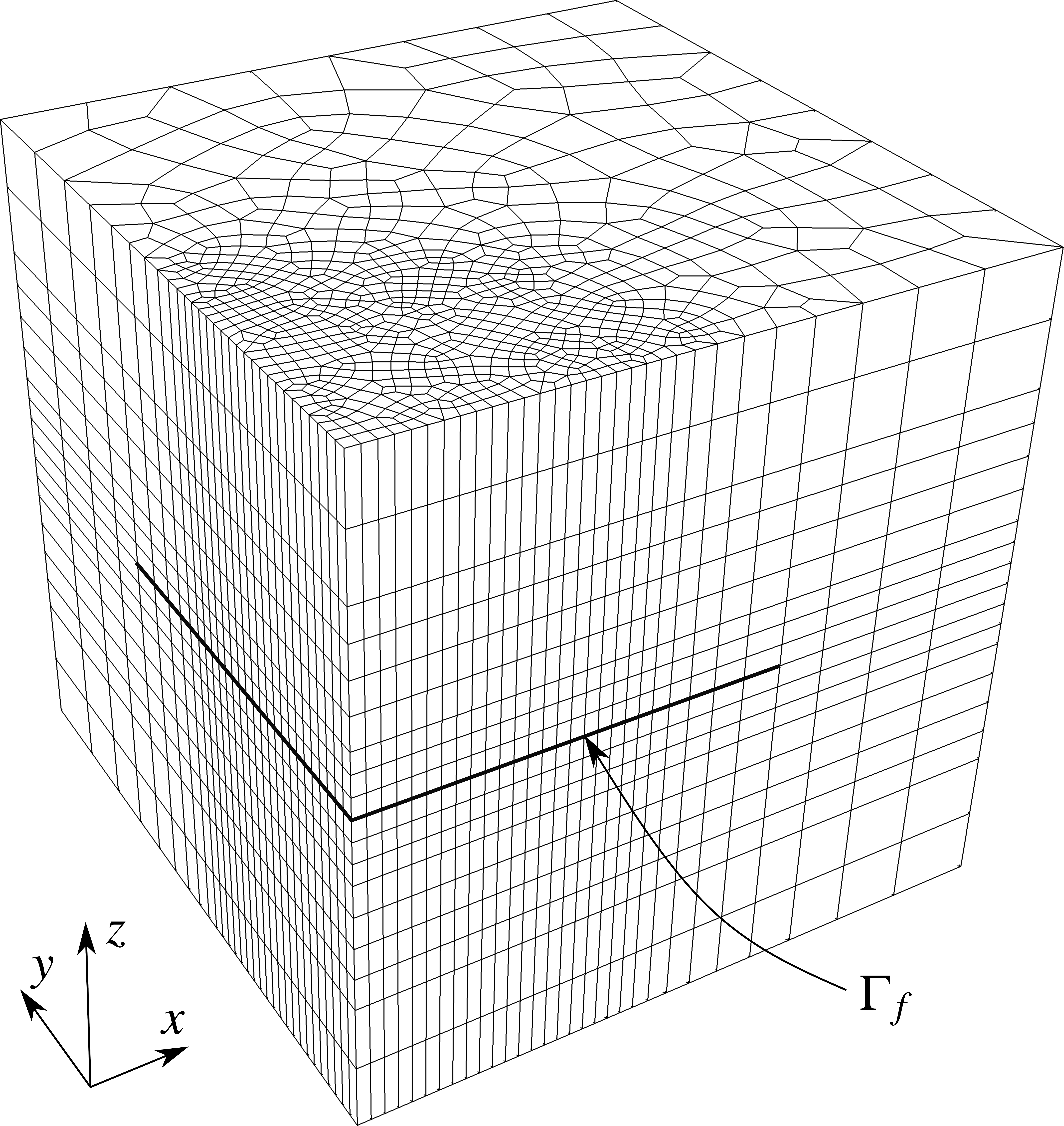}}
    \caption{}  
    \label{fig:penny_mesh_3d}
  \end{subfigure}
  \hfill\null
  \caption{Penny-shaped crack problem: (a) sketch of the setup; (b) $x$-$y$ view of computational domain; (c) and full 3D mesh. In (b), the maximum area of the fracture is highlighted, while in (c) the black line indicating the fracture trace.}
\end{figure}

The problem consists on an axisymmetric hydraulic fracture in an infinite medium, that is
isotropic, homogeneous, impermeable and with a linear elastic behavior. From the center of
the fracture, an incompressible Newtonian fluid is injected, at a constant rate $Q_0$. The
fracture propagation does not depend on the far-field stress status and, as in the KGD
example, it is enough to solve for the net fluid pressure. In Fig. \ref{fig:penny_sketch},
we represent the setup for the problem. The quantities of interest are the same as before,
except the fracture half-length, that is now substitute by $R(t)$, i.e., the fracture
radius.

The behavior in space and time, represented by $g_N(x,t)$, $p(x,t)$ and $L(t)$, is provided
by the analytical solution through some dimensionless quantities
\cite{detournay2004propagation}: (i) an opening $\Omega$, (ii) a net pressure $\Pi$ and
(iii) a fracture radius $\gamma$. The main quantities can be expressed from these
dimensionless functions as:
\begin{subequations}
\begin{align}
  R(t) &= \gamma L(t), &
  g_N(x,t) &= \varepsilon(t) L(t) \Omega(\rho), &
  p(x,t) &= \varepsilon(t) E^\prime \Pi(\rho), &
\end{align}
\label{eq:penny_sol1}\null
\end{subequations}
where the similarity variable $\rho = r/R$ is the dimensionless fracture coordinate. For
the zero toughness case, we have \cite{savitski2002propagation}:
\begin{subequations}
\begin{align}
  \varepsilon(t) &= \left(\frac{\mu^\prime}{E^\prime t}\right)^{\frac{1}{3}},
&
  L(t) &= Q_0^{1/3} \left(\frac{E^\prime}{\mu^\prime}\right)^{\frac{1}{9}} t^{4/9},
&
  \gamma &= \frac{1}{2\pi \left(\int_0^1 \overline{\Omega}\rho \, \mathrm{d}
\rho\right)^{1/3}}.
&
\end{align}
\label{eq:penny_sol2}\null
\end{subequations}
In \eqref{eq:penny_sol2}, we introduced $\overline{\Omega} = \frac{\Omega}{\gamma}$,
$\mu^\prime = 12 \mu_l$ and $E^\prime = \frac{E}{1-\nu^2}$ as in the KDG study. The
self-similar fracture opening $\overline{\Omega}(\rho)$ and the self-similar net fluid
pressure $\Pi(\rho)$ are expressed as sum of a general and a particular solution. The
first respects the governing equations, while the second represents the inlet asymptotic
behavior. According to \cite{savitski2002propagation}, the chosen basis function is a
combination of Jacobi polynomials of arbitrary order in the interval $[0,1]$. In this
reference, there is the complete expression for the two self-similarity functions. In the
current work, we use just the first order expansion for both $\overline{\Omega}(\rho)$ and
$\Pi(\rho)$.
The
pressure solution is unbounded, being $\lim_{\rho \rightarrow 0} \Pi = +\infty$ and
$\lim_{\rho \rightarrow 1} \Pi = -\infty$, and the model can struggle in the approximation
of this nonphysical values happening at the extremes of the domain. As for the KGD
solution, we highlight that the analytical solutions are not affected by neither the fluid
density nor the initial stress regime.

The computational domain simulates one fourth of the problem, composed by 16061 nodes,
13896 hexahedra and 640 interface elements, as shown in Figs. \ref{fig:penny_mesh_2d}-\ref{fig:penny_mesh_3d}. The
global size is $100 \times 100 \times 100 \text{ m}$, with a fracture radius of about $60
\text{ m}$ and an average element area of $2.4 \text{ m}^2$. The discretized fracture
surface is large enough to allow the propagation of the fracture without geometric
constraints. The symmetric boundary conditions are imposed on the two symmetry plane,
i.e., the boundaries parallel to $x$ and $y$ axes containing the well. The other two faces
parallel to $x$ and $y$ axes are $z$-constrained. The material parameters are $E = 30
\text{ GPa}$ and $\nu = 0.25$. The fracture frictional behavior is governed by Coulomb's
criterion, characterized by $\theta = 30^\circ$ and zero cohesion. The fluid viscosity is
$\mu_l = 10^{-9} \text{ MPa}\cdot\text{s}$. The injection rate and the confining stress
are $Q_0 = 1 \cdot 10^{-2} \text{ m}^3/\text{s}$ and $\sigma_0 = 10 \text{ MPa}$,
respectively. Finally, we use the same value of initial conductivity as in the KGD
example, i.e., $C_{f,0} = 10 \text{ mD}\cdot\text{m} = 9.87 \cdot 10^{-15}
\text{ m}^2\cdot\text{m}$. As for the previous case, the simulated time is $100 \text{
s}$, with $\dt = 0.01 \text{ s}$ from the beginning to $t = 0.2 \text{ s}$, then $\dt =
0.1 \text{ s}$ up to $t = 2 \text{ s}$ and finally $\dt = 1 \text{ s}$ up to the end of
the simulation. According to \cite{lai2015experimental}, the effective Reynolds number is
$1.1 \cdot 10^{-1}$ after 1 s and $3.2 \cdot 10^{-3}$ at the end of the simulation. Thus,
the laminar flow assumption required by the cubic law holds. We emphasize that
the proposed mesh is not really suited for a TPFA finite volumes solution scheme,
nevertheless, the results prove to be accurate enough to verify the analytical solution.

\begin{figure}
\hfill
\begin{subfigure}{0.24\linewidth}
  \centering
  \includegraphics[width=\linewidth]{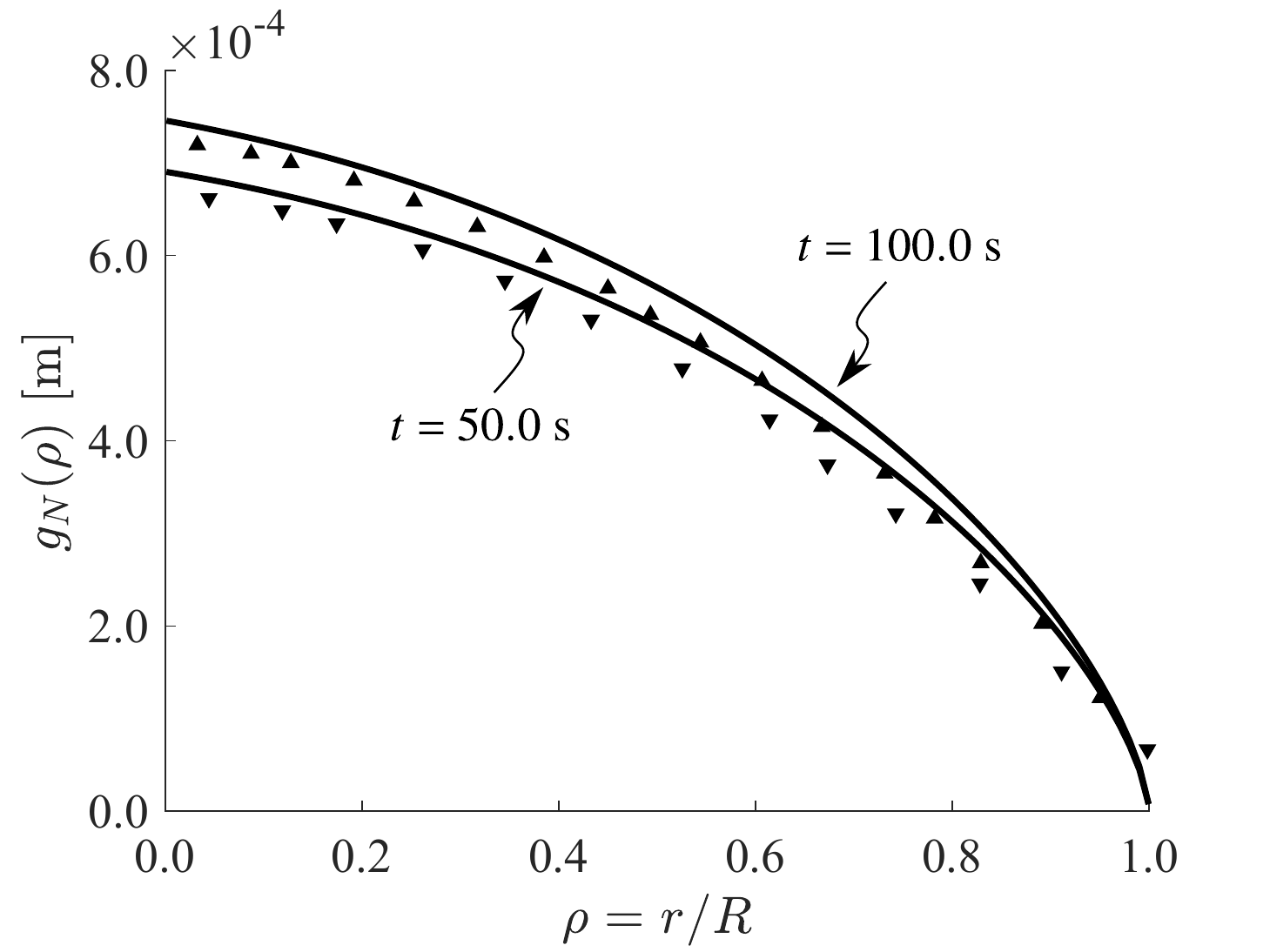}
  \caption{}
  \label{fig:penny_opening}
\end{subfigure}
\begin{subfigure}{0.24\linewidth}
  \centering
  \includegraphics[width=\linewidth]{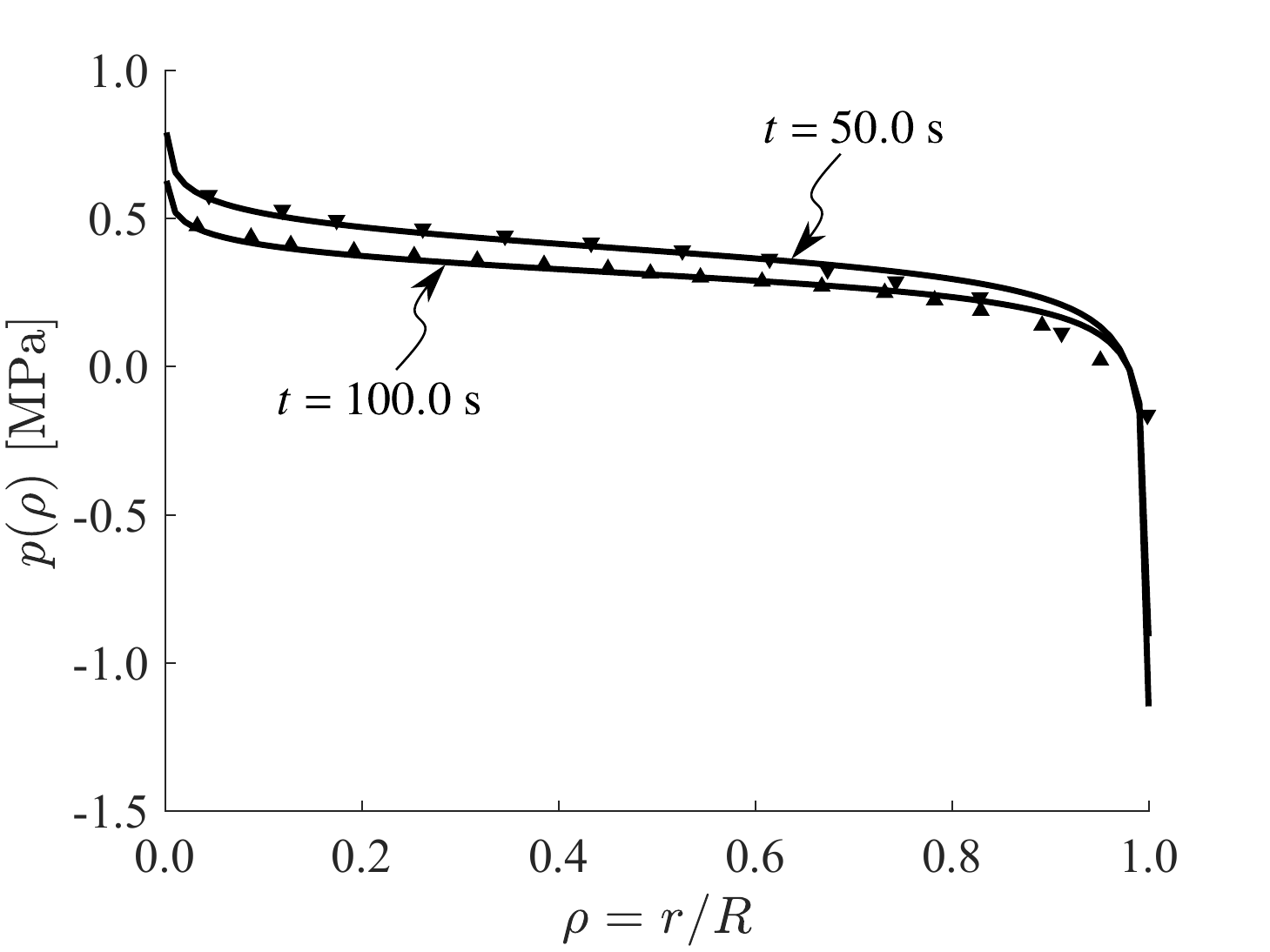}
  \caption{}
  \label{fig:penny_pressure}
\end{subfigure}
\begin{subfigure}{0.24\linewidth}
  \centering
  \includegraphics[width=\linewidth]{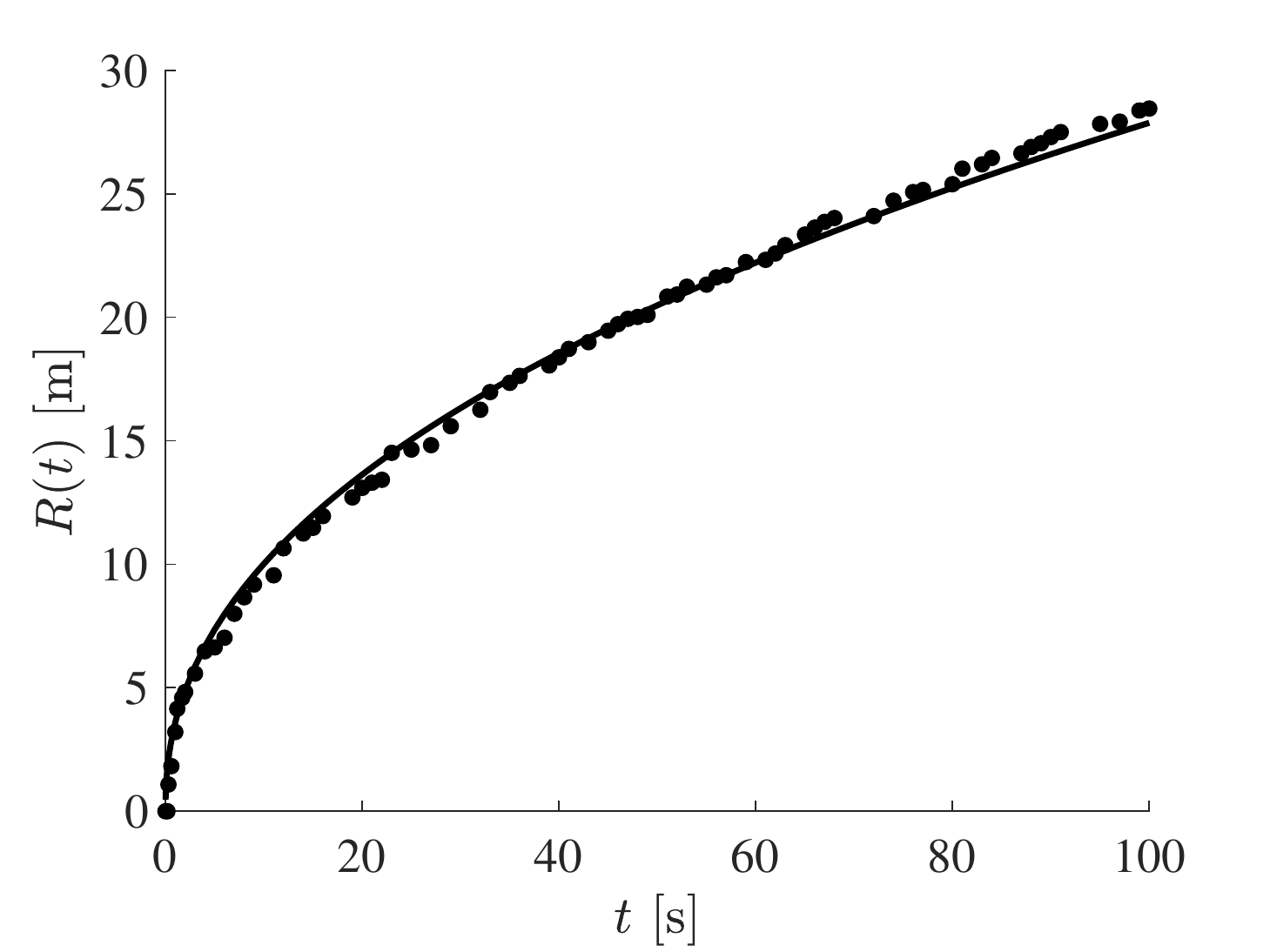}
  \caption{}
  \label{fig:penny_length}
\end{subfigure}
\begin{subfigure}{0.24\linewidth}
  \centering
  \includegraphics[width=\linewidth]{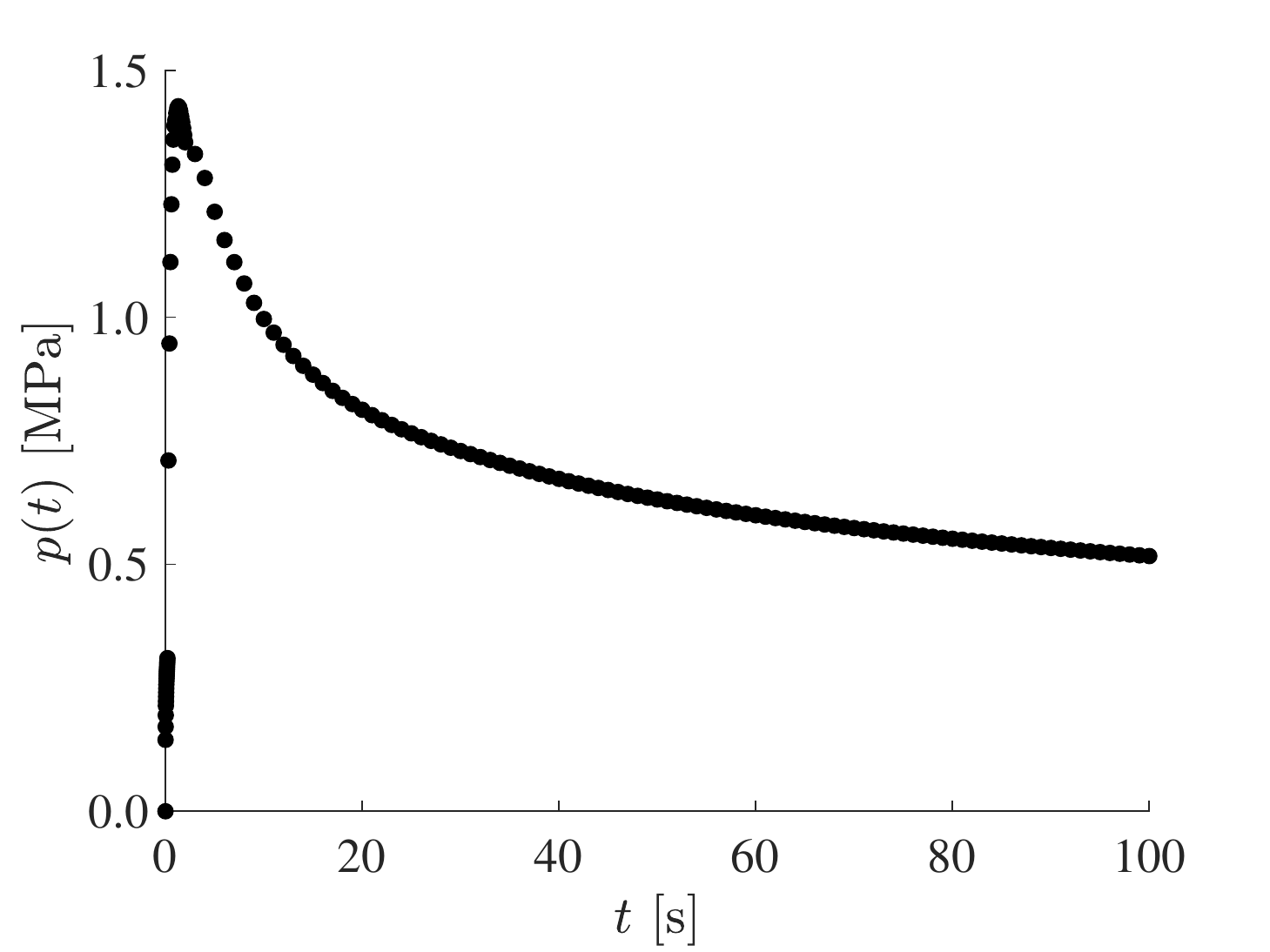}
  \caption{}
  \label{fig:penny_pressure_x0}
\end{subfigure}
\hfill\null
  \caption{Results for the penny-shaped crack simulation. From the left to the right:
    fracture opening, fluid pressure, fracture radius and pressure behavior at the
    injection location.} \label{fig:penny_res}
\end{figure}

For two different simulated times, i.e., at half ($t = 50 \text{ s}$) and at the end of
the simulation ($t = 100 \text{ s}$), we show the opening and pressure profiles for the $y
= 0$ fracture trace in Figs. \ref{fig:penny_opening}-\ref{fig:penny_pressure}, where the
continuous line is the analytical solution. Overall, there is a quite good agreement, for
both the aperture (with an integral relative error of $5.2 \%$ for $t = 100 \text{ s}$)
and the pressure (with an integral relative error of $4.2 \%$ at the same time step). The
pressure is slightly different close to the extremes of the domain, i.e., $\rho = 0$ and
$\rho = 1$, because the theoretical behavior diverges, tending to $\pm\infty$. As in the
KGD case, the propagation criterion provides an integral constraint for the pressure
analytical solution, being
\begin{equation}
  \int_0^1 \frac{\Pi}{\sqrt{1-\rho^2}}\rho \, \mathrm{d} \rho = 0.
\end{equation}
Using a simple Dirichlet boundary condition on the pressure value at the fracture maximum
radius, where $p = 0$ is imposed, we need to shift our outcome for the sake of comparison.
Fig. \ref{fig:penny_length} represents the fracture radius as function of time. The model
is able to predict quite accurately the fracture length, with an average relative error of
$2.2 \%$. Finally, in Fig. \ref{fig:penny_pressure_x0}, the pressure behavior
for the entire simulation at the injection point is shown. As for the KGD example, we rely on a
visual comparison with other solutions available in the literature, see for instance
\cite{Set_etal17}.

\section{Conclusions}
\label{sec:concl}
In this work, we have presented a stabilized displacement-Lagrange multiplier-pressure formulation for quasi-static contact mechanics coupled with fracture fluid flow.
Our discretization is based on a finite element method for the contact mechanics subproblem combined with a finite volume scheme for the flow subproblem.
The global nonlinear problem is solved using an active set strategy.
We employ lowest-order continuous finite elements for the displacement field, piecewise constant functions for both Lagrange multiplier (traction) and pressure field, and a linear two-point flux approximation for intercell numerical fluxes on the discontinuity surface.

The mixed space we adopt does not automatically satisfy the discrete inf-sup stability condition.
To stabilize the formulation, we began by revisiting the well-known macroelement approach, originally proposed for the Stokes equation.
We then developed two new algebraic strategies that work under less restrictive assumptions while providing good performance.
The resulting approach, both without and with fluid flow in the fracture, has been benchmarked against analytical solutions to validate its behavior in case of normal and frictional activation, as well as in two- and three-dimensional fracture reactivation.

Future developments will deal with the simulation of fluid flow in the matrix coupled with the structural mechanics problem through porosity variation and fracture flow through leak-off.
Networks of fractures, for which the computation of the hydraulic conductivity requires particular care, will also be investigated.
Finally, a key step to enable high-performance computing applications is the design of an efficient preconditioning strategy for the particular Jacobian linear systems encountered here.

\section*{Acknowledgements}
\label{sec::acknow}
Funding was provided by TOTAL S.A. through the FC-MAELSTROM project.
The authors wish to thank Randolph Settgast for helpful discussions.
Portions of this work were performed under the auspices of the U.S. Department of
Energy by Lawrence Livermore National Laboratory under Contract DE-AC52-07-NA27344.

\appendix
\section{Finite element and finite volume vectors and matrices}\label{app:FEFV_vec_mat_blocks}

The linearized system \eqref{eq:jac_sys} is assembled in the standard way from elementary contributions.
The global expressions for the residual block vectors read:
\begin{subequations}
\begin{align}
  [ \Vec{r}_u ]_i =
  & ( \nabla^s \tensorOne{\eta}_i, \tensorTwo{\sigma}_n )_{\Omega}
  + ( \llbracket \tensorOne{\eta}_{i} \rrbracket, \tensorOne{t}^h_{n} - p^h_n \tensorOne{n}_f )_{\Gamma_f}
  - ( \tensorOne{\eta}_i, \bar{\vec{t}}_n )_{\partial \Omega_{\sigma}}
  && \forall i \in \mathcal{N}_u,
  \label{eq:r_disp_new}\\
  [ \Vec{r}_S ]_i =
  & ( \mu_{i,N},  g_{N,n} )_{\Gamma_{f,n}^{\ell,\text{stick}} }
  + ( \tensorOne{\mu}_{i,T}, \Delta_n \tensorOne{g}_T )_{\Gamma_{f,n}^{\ell,\text{stick}}}
  && \forall i \in \mathcal{N}_t^S,
  \label{eq:r_lam_stick}\\
  [ \Vec{r}_N ]_i =
  & ( \mu_{i,N}, g_{N,n} )_{ \Gamma_{f,n}^{\ell,\text{slip}}}
  && \forall i \in \mathcal{N}_{t}^N,
  \label{eq:r_lam_slip_N}\\
  [ \Vec{r}_T ]_i =
  & ( \tensorOne{\mu}_{i,T}, \tensorOne{t}^h_{T,n} - \tensorOne{t}^*_n )_{\Gamma_{f,n}^{\ell,\text{slip}}}
  && \forall i \in \mathcal{N}_{t}^T,
  \label{eq:r_lam_slip_T}\\
  [ \Vec{r}_O ]_i =
  & ( \tensorOne{\mu}_{i}, \tensorOne{t}^h_n )_{\Gamma_{f,n}^{\ell,\text{open}}}
  && \forall i \in \mathcal{N}_t^O,
  \label{eq:r_lam_open}\\
  [ \Vec{r}_p ]_i =
  & \left( \chi_i, \frac{\Delta_n g_{N} }{\Delta_n t} \right)_{\Gamma_f^{\ell,\text{open}}}
  + [ \chi_i, p_n^h ]_{\mathcal{F}_f} 
  - F_{\mathcal{F}_f}(\chi_i)
  + G_{\mathcal{F}_f}(\chi_i)
  - ( \chi_i, q_{s,n} )_{\Gamma_f},
  && \forall i \in \mathcal{N}_p.
  \label{eq:r_pres_new}
\end{align}
\label{eq:residual}\null
\end{subequations}
The global expressions for the sub-matrices appearing in the Jacobian matrix read:
\begin{subequations}
\begin{align}
  [ A_{uu} ]_{ij} =
  & ( \nabla^s \tensorOne{\eta}_i, \tensorFour{C} : \nabla^s \tensorOne{\eta}_j )_{\Omega}
  && \forall(i,j) \in \mathcal{N}_u \times \mathcal{N}_u,
  \label{eq:jac_blk_Auu} \\
  [ A_{u\beta} ]_{ij} =
  & ( \llbracket \eta_{i,N} \rrbracket, \mu_{j,N} )_{\Gamma_{f,n}}
  + ( \llbracket \tensorOne{\eta}_{i,T} \rrbracket, \tensorOne{\mu}_{j,T} )_{\Gamma_{f,n}}
  && \forall(i,j) \in \mathcal{N}_u \times \mathcal{N}_t^{\beta = \{ S,N,T,O \} },
  \label{eq:jac_blk_AuSNTO} \\
  [ A_{up} ]_{ij} =
  & - ( \llbracket \eta_{i,N} \rrbracket, \chi_j )_{\Gamma_{f,n}}
  && \forall(i,j) \in \mathcal{N}_u \times \mathcal{N}_p,
  \label{eq:jac_blk_Aup} \\
  [ A_{Su} ]_{ij} =
  & ( \mu_{i,N}, \llbracket \eta_{j,N} \rrbracket )_{\Gamma_{f,n}^{\ell,\text{stick}} }
  + ( \tensorOne{\mu}_{i,T}, \llbracket \tensorOne{\eta}_{j,T} \rrbracket )_{\Gamma_{f,n}^{\ell,\text{stick}}}
  &&  \forall(i,j) \in \mathcal{N}_t^S \times \mathcal{N}_u,
  \label{eq:jac_blk_ASu} \\
  [ A_{Nu} ]_{ij} =
  & ( \mu_{i,N}, \llbracket \eta_{j,N} \rrbracket )_{\Gamma_{f,n}^{\ell,\text{slip}} }
  && \forall(i,j) \in \mathcal{N}_t^N \times \mathcal{N}_u,
  \label{eq:jac_blk_ANu} \\
  [ A_{Tu} ]_{ij} =
  & - \left( \tensorOne{\mu}_{i,T}, \left.\frac{\partial \vec{t}_T^* }{\partial \Delta_n \tensorTwo{g}_T } \right|_n^{\ell,(k)} \cdot \llbracket \tensorOne{\eta}_{j,T} \rrbracket \right)_{\Gamma_{f,n}^{\ell,\text{slip}}}
  && \forall(i,j) \in \mathcal{N}_t^T \times \mathcal{N}_u,
  \label{eq:jac_blk_ATu} \\
  [ A_{TN} ]_{ij} =
  & - \left( \tensorOne{\mu}_{i,T}, \left.\frac{\partial \vec{t}_T^* }{\partial t_N } \right|_n^{\ell,(k)} \mu_{j,N} \right)_{ \Gamma_{f,n}^{\ell,\text{slip}} }
  && \forall(i,j) \in \mathcal{N}_t^T \times \mathcal{N}_t^N,
  \label{eq:jac_blk_ATN} \\
  [ A_{TT} ]_{ij} =
  & ( \tensorOne{\mu}_{i,T}, \tensorOne{\mu}_{j,T} )_{ \Gamma_{f,n}^{\ell,\text{slip}} }
  && \forall(i,j) \in \mathcal{N}_t^T \times \mathcal{N}_t^T,
  \label{eq:jac_blk_ATT} \\
  [ A_{OO} ]_{ij} =
  & ( \mu_{i,N}, \mu_{j,N} )_{\Gamma_{f,n}^{\ell,\text{open}} }
  + ( \tensorOne{\mu}_{i,T}, \tensorOne{\mu}_{j,T} )_{\Gamma_{f,n}^{\ell,\text{open}} }
  && \forall(i,j) \in \mathcal{N}_t^O \times \mathcal{N}_t^O,
  \label{eq:jac_blk_AOO} \\
  [ A_{pu} ]_{ij} =
  & \frac{1}{\Delta_n t} ( \chi_i, \llbracket \eta_{j,N} \rrbracket )_{\Gamma_{f,n}^{\ell,\text{open}} }
  + \left.\frac{\partial ( [ \chi_i, p^h ]_{\mathcal{F}_f} ) }{\partial u_j } \right|_n^{\ell,(k)}
  - \left.\frac{\partial ( F_{\mathcal{F}_f}( \chi_i ) ) }{\partial u_j } \right|_n^{\ell,(k)}
  && \forall(i,j) \in \mathcal{N}_p \times \mathcal{N}_u,
  \label{eq:jac_blk_Apu} \\
  [ A_{pp} ]_{ij} =
  &
  \left.\frac{\partial ( [ \chi_i, p^h ]_{\mathcal{F}_f} ) }{\partial p_j } \right|_n^{\ell,(k)} 
  && \forall(i,j) \in \mathcal{N}_p \times \mathcal{N}_p,
  \label{eq:jac_blk_App}
\end{align}
\label{eq:jac_blk_new}\null
\end{subequations}
with the partial derivatives expanded as:
\begin{subequations}
\begin{align}
  \left. \frac{\partial \vec{t}_T^* }{\partial \Delta_n \vec{g}_T }
\right|_n^{\ell,(k)} &
  = \left. \tau_{\text{max}}(t^{h}_{N}) \frac{|| \Delta_n \vec{g}_T ||_2^2 \tensorTwo{1}
  - \Delta_n \vec{g}_T \otimes \Delta_n \vec{g}_T }{|| \Delta_n \vec{g}_T ||_2^3} \right|_{n}^k,
  \label{eq:jac_der1_new} \\
  \left. \frac{\partial \vec{t}_T^* }{\partial t_N } \right|_n^{\ell,(k)} &=
\left. - \tan(\theta) \frac{ \Delta_n \vec{g}_T }{|| \Delta_n \vec{g}_T ||_2} \right|_{n}^k, \label{eq:jac_der2_new} \\
  \left.\frac{\partial ( [ \chi_i, p^h ]_{\mathcal{F}_f} ) }{\partial u_j } \right|_n^{\ell,(k)}  &=
  \smashoperator{\sum_{\varepsilon = \varepsilon_{K,L} \in \mathcal{E}_{\text{int}} }}
  ( \chi_i\vert_{ \varphi_L } - \chi_i\vert_{ \varphi_K } ) \left ( \left.\frac{\partial \Upsilon_{KL} }{\partial u_j } \right|_n^{\ell,(k)} \right) ( p^h\vert_{ \varphi_L } - p^h\vert_{ \varphi_K } )
  + \smashoperator{\sum_{\varepsilon = \varepsilon_{K} \in \mathcal{E}_p }}
  \chi\vert_{ \varphi_K } \left ( \left.\frac{\partial \Upsilon_{K} }{\partial u_j } \right|_n^{\ell,(k)} \right) \, p^h\vert_{ \varphi_K },
  \label{eq:jac_der3d_new} \\
  \left.\frac{\partial ( F_{\mathcal{F}_f}( \chi_i ) ) }{\partial u_j } \right|_n^{\ell,(k)} &=
  \smashoperator{\sum_{\varepsilon = \varepsilon_{K} \in \mathcal{E}_p }}
  \chi_i\vert_{ \varphi_K } \left ( \left.\frac{\partial \Upsilon_{K} }{\partial u_j } \right|_n^{\ell,(k)} \right) \, \frac{1}{|\varepsilon|}\int_{\varepsilon} \bar{p}_{n} \, \mathrm{d}l,
  \label{eq:jac_der4d_new} \\
  \left.\frac{\partial \Upsilon_{KL} }{\partial u_j } \right|_n^{\ell,(k)} &=
  \frac{\Upsilon_{KL}^2 \overline{\Upsilon}_K}{\Upsilon_K^2 \, |\varphi_K|} \int_{\varphi_{K} } \frac{ g_N^2 }{4} \, \llbracket \eta_{j,N} \rrbracket \, \mathrm{d}A
  + \frac{\Upsilon_{KL}^2 \overline{\Upsilon}_L}{\Upsilon_L^2 \, |\varphi_L|}
  \int_{\varphi_{L} } \frac{ g_N^2 }{4} \, \llbracket \eta_{j,N} \rrbracket \, \mathrm{d}A
  \label{eq:jac_der5d_new} \\
  \left.\frac{\partial \Upsilon_K }{\partial u_j } \right|_n^{\ell,(k)} &=
  \frac{\overline{\Upsilon}_K }{| \varphi_{K} |} \int_{\varphi_{K} } \frac{ g_N^2 }{4} \, \llbracket \eta_{j,N} \rrbracket \, \mathrm{d}A,
  \label{eq:jac_der6d_new} \\
  \left.\frac{\partial ( [ \chi_i, p^h ]_{\mathcal{F}_f} ) }{\partial p_j } \right|_n^{\ell,(k)}  &=
  \smashoperator{\sum_{\varepsilon = \varepsilon_{K,L} \in \mathcal{E}_{\text{int}} }}
  ( \chi_i\vert_{ \varphi_L } - \chi_i\vert_{ \varphi_K } ) \Upsilon_{KL} ( \chi_j\vert_{ \varphi_L } - \chi_j\vert_{ \varphi_K } )
  + \smashoperator{\sum_{\varepsilon = \varepsilon_{K} \in \mathcal{E}_p }}
  \chi\vert_{ \varphi_K } \Upsilon_K \, \chi_j\vert_{ \varphi_K },
  \label{eq:jac_der7d_new}
\end{align}
\label{eq:jac_der_new}\null
\end{subequations}
Note that Eqs. \eqref{eq:jac_der3d_new}-\eqref{eq:jac_der4d_new}, provide nonzero contributions only for faces in the open mode, i.e., $\varphi \in \Gamma_{f,n}^{\ell,\text{open}}$ .

\end{document}